\documentclass[11pt,twoside,leqno]{aomamlt2e} 
  \pageno{711}
 
\received{December 20, 2002}
  \begingroup \makeatletter
\@ifundefined{theorem@style}{\input{theorem.sty}}{}
             [\FMithmInfo]
\gdef\th@plain{\normalfont\itshape
  \def\@begintheorem##1##2{%
        \item[\hskip\labelsep \theorem@headerfont \hskip22pt ##2\ ##1.]}%
\def\@opargbegintheorem##1##2##3{%
   \item[\hskip\labelsep \theorem@headerfont\hskip22pt  ##2\ ##1 {\rm (##3)}.\ ]}}
\endgroup

\newcommand\sdemo[1]{\demo{\scshape #1}}
   \newtheorem{theorem}{Theorem}[section]
   
   \newtheorem{lemma}[theorem]{Lemma}

   \newtheorem{corollary}[theorem]{Corollary}
   \newtheorem{proposition}[theorem]{Proposition}
   \newtheorem{main-theorem}[theorem]{Main Theorem}
   \newtheorem{key-estimates}[theorem]{Key Estimates}
   
   \newtheorem{key-lemma}[theorem]{Key-Lemma}
\theorembodyfont{\upshape}
   \newtheorem{definition}[theorem]{{\it Definition}}
   \newtheorem{remark}[theorem]{{\it Remark}}

\numberwithin{equation}{section} 

\def\C{{\mathbb C}}
\def\E{{\mathbb E}}

\def\N{{\mathbb N}}
\def\R{{\mathbb R}}
\def\TT{{\mathbb T}}
\def\Z{{\mathbb Z}}
\def\X{{\mathbb X}}
\def\Y{{\mathbb Y}}
\def\V{{\mathbb V}}
\def\Q{{\mathbb Q}}

\def\CA{{\mathcal A}}
\def\CB{{\mathcal B}}
\def\CC{{\mathcal C}}
\def\CD{{\mathcal D}}
\def\CE{{\mathcal E}}
\def\CF{{\mathcal F}}

\def\CH{{\mathcal H}}
\def\CI{{\mathcal I}}
\def\CJ{{\mathcal J}}
\def\CK{{\mathcal K}}
\def\CL{{\mathcal L}}

\def\CO{{\mathcal O}}
\def\CP{{\mathcal P}}

\def\CR{{\mathcal R}}

\def\CT{{\mathcal T}}
\def\CU{{\mathcal U}}

\def\1{{{\ae}}}
\def\2{{{\o}}}
\def\3{{{\aa}}}
\def\6{\, {\rm d}}
\def\i{{\rm i}}
\def\e{{\rm e}}
\def\<{{\langle}}
\def\>{{\rangle}}
\def\brinv{{\langle -1\rangle}}

\def\unit{{{\pmb 1}}}
\def\spe{{\text{\rm sp}}}

\def\ep{{\varepsilon}}

\def\lmax{{\lambda_{\max}}}
\def\lmin{{\lambda_{\min}}}

\def\id{{\text{\rm id}}}
\def\tr{{\text{\rm tr}}}
\def\Tr{{\text{\rm Tr}}}
\def\lmax{{\lambda_{\max}}}
\def\GRM{{\text{\rm GRM}}}
\def\SGRM{{\text{\rm SGRM}}}
\def\GL{{\rm GL}}
\def\Mnsa{{M_n(\C)_{\rm sa}}}
\def\Mmsa{{M_m(\C)_{\rm sa}}}

\def\Ext{{\rm Ext}}
\def\Cred{C^*_{\rm red}}
\def\grad{{\rm grad}}
\def\diff{\frac{\rm d}{{\rm d}t}_{\bigm |_{t=0}}}

\def\im{{\rm Im}\, }
\def\re{{\rm Re}}
\def\supp{{\rm supp}}

\begin{document}
\currannalsline{162}{2005} 

 \title{A new application of random matrices:\\ $\Ext(\Cred(F_2))$ is  
not a
group}

 \acknowledgements{This work was carried out, while the first author was a member of the MaPhySto  --
Centre for Mathematical Physics and Stochastics, funded by the Danish National Research Foundation. 
The second author was supported by The Danish National Science Research Council.}
\twoauthors{Uffe Haagerup}{Steen
Thorbj{\o}rnsen}

 \institution{University of Southern Denmark, Odense, Denmark
\\
\emails{haagerup@imada.sdu.dk}{steenth@imada.sdu.dk}}

 \shorttitle{A new application of random matrices}
\shortname{Uffe Haagerup and Steen
Thorbj\hskip1pt\raise1pt\hbox{$\slash$}{\hskip-5.5pt\hbox{o}}rnsen}

 \vglue-12pt
\centerline{\it Dedicated to the memory of Gert Kj{\ae}rg{\aa}rd Pedersen}
 
\vglue10pt

\centerline{\bf Abstract}
\vglue6pt

In the process of developing the theory of free probability and free
entropy, Voiculescu introduced in 1991 a random matrix model for a free
semicircular system. Since then, random matrices have played a key role
in von Neumann algebra theory (cf.\ \cite{V8}, \cite{V9}). The main
result of this paper is the following extension of Voiculescu's random
matrix result: Let $(X_1^{(n)},\dots,X_r^{(n)})$ be a system of $r$
stochastically independent $n\times n$ Gaussian self-adjoint random
matrices as in Voiculescu's random matrix paper \cite{V3}, and let
$(x_1,\dots,x_r)$ be a semi-circular system in a $C^*$-probability
space. Then for every polynomial $p$ in $r$ noncommuting variables
\[
\lim_{n\to\infty}
\big\|p\big(X_1^{(n)}(\omega),\dots,X_r^{(n)}(\omega)\big)\big\|
=\|p(x_1,\dots,x_r)\|,
\]
for almost all $\omega$ in the underlying probability space. We use the
result to show that the $\Ext$-invariant for the reduced $C^*$-algebra
of the free group on 2 generators is not a group but only a
semi-group. This problem has been open since Anderson in 1978 found the
first example of a $C^*$-algebra $\CA$ for which $\Ext(\CA)$ is not a  
group.
 
\vglue-20pt
\phantom{goone}
\section{Introduction}\label{sec0}
\vglue-6pt

A random matrix $X$ is a matrix whose entries are real or complex random
variables on a probability space $(\Omega,\CF,P)$. As in
\cite{T}, we denote by $\SGRM(n,\sigma^2)$ the class of complex
self-adjoint $n\times n$ random matrices
\[
X = (X_{ij})^n_{i,j=1},
\]
for which $(X_{ii})_i$, $(\sqrt{2} \re X_{ij})_{i<j}$, $(\sqrt{2} \im
X_{ij})_{i<j}$ are $n^2$ independent identically distributed (i.i.d.)
Gaussian random variables with mean value 0 and variance~$\sigma^2$. In  
the
terminology of Mehta's book \cite{Me}, $X$ is a Gaussian unitary
ensemble (GUE). In the following we put $\sigma^2=\frac 1n$ which is the
normalization used in Voiculescu's random matrix paper \cite{V3}. We
shall need the following basic definitions from free probability theory
(cf.~\cite{V1}, \cite{vdn}):
\begin{itemize}
\item[a)] A $C^*$-probability space is a pair $(\CB,\tau)$ consisting  
of a
unital $C^*$-algebra $\CB$ and a state $\tau$ on $\CB$.
\item[b)] A family of elements $(a_i)_{i\in I}$ in a $C^*$-probability
space $(\CB,\tau)$ is free if for all $n\in\N$ and all polynomials
$p_1,\dots,p_n\in\C[X]$, one has
\[
\tau(p_1(a_{i_1})\cdots p_n(a_{i_n})) = 0,
\]
whenever $i_1\ne i_2, i_2\ne i_3,\dots,i_{n-1}\ne i_n$
and $\varphi(p_k(a_{i_k}))=0$ for $k=1,\dots,n$.
\item[c)] A family $(x_i)_{i\in I}$ of elements in a $C^*$-probability
space $(\CB,\tau)$ is a semicircular family, if $(x_i)_{i\in I}$ is a  
free
family, $x_i=x_i^*$ for all $i\in I$ and
\[
\tau(x_i^k) = \frac{1}{2\pi} \int^2_{-2}
t^k\sqrt{4-t^2}\6t=
\begin{cases}
\frac{1}{k/2+1}\binom{k}{k/2}, &\mbox{if $k$ is even}, \\
0, &\mbox{if $k$ is odd},
\end{cases}
\]
for all $k\in\N$ and $i\in I$.
\end{itemize}

We can now formulate Voiculescu's random matrix result from \cite{V4}:
Let, for each $n\in\N$, $(X_i^{(n)})_{i\in I}$ be a family of
independent random matrices from the class $\SGRM(n,\frac 1n)$, and let
$(x_i)_{i\in I}$ be a semicircular family in a $C^*$-probability
space $(\CB,\tau)$. Then for all $p\in \N$ and all $i_1,\dots,i_p\in  
I$, we
have
\begin{equation} \label{eq0-1}
\lim_{n\to\infty}\E\big\{\tr_n\big(X_{i_1}^{(n)}\cdots
X_{i_p}^{(n)}\big)\big\}=\tau(x_{i_1}\cdots x_{i_p}),
\end{equation}
where $\tr_n$ is the normalized trace on $M_n(\C)$, i.e.,
$\tr_n=\frac{1}{n}\Tr_n$, where $\Tr_n(A)$ is the sum of the diagonal
elements of $A$. Furthermore, $\E$ denotes expectation (or integration)
with respect to  the probability measure $P$.

The special case $|I|=1$ is Wigner's semi-circle law (cf.~\cite{Wi},
\cite{Me}). The strong law corresponding to \eqref{eq0-1} also holds,  
i.e.,
\begin{equation}
\label{eq0-2}
\lim_{n\to\infty} \tr_n\big(X_{i_1}^{(n)}(\omega)\cdots
X_{i_p}^{(n)}(\omega)\big)=\tau(x_{i_1}\cdots x_{i_p}),
\end{equation}
for almost all $\omega\in\Omega$ (cf.~\cite{Ar} for the case $|I|=1$ and
\cite{HP2}, \cite[Cor.~3.9]{T} for the general case). Voiculescu's  
result is
actually more general than the one quoted above. It also involves
sequences of non random diagonal matrices. We will, however, only
consider the case, where there are no diagonal matrices. The main result
of this paper is that the strong version \eqref{eq0-2} of Voiculescu's
random matrix result also holds for the operator norm in the following
sense:

\demo{\scshape Theorem A}   {\it Let $r\in\N$ and{\rm ,} for each $n\in\N${\rm ,} let
   $(X_1^{(n)},\dots,X_r^{(n)})$ be
a set of $r$ independent random matrices from the class $\SGRM(n,\frac
1n)$. Let further $(x_1,\dots,x_r)$ be a semicircular system in a
$C^*$\/{\rm -}\/probability space $(\CB,\tau)$ with a faithful state $\tau$. Then
there is a $P$\/{\rm -}\/null set $N\subseteq\Omega$ such that for all
$\omega\in\Omega\backslash N$ and all polynomials $p$ in $r$
noncommuting variables{\rm ,} we have}
\begin{equation}
\label{eq0-3}
\lim_{n\to\infty}\big\|p\big(X_1^{(n)}(\omega),
\dots,X_r^{(n)}(\omega)\big)\big\|=\|p(x_1,\dots,x_r)\|.
\end{equation}

\Enddemo

The proof of Theorem A is given in Section~\ref{sec6}. The special case
\[
\lim_{n\to\infty}\big\|X_1^{(n)}(\omega)\big\|=\|x_1\|=2
\]
is well known (cf.\ \cite{BY}, \cite[Thm.~2.12]{Ba}
or \cite[Thm.~3.1]{HT1}).

 From Theorem A above, it is not hard to obtain the following
result (cf.~\S \ref{ext-resultat}).

\demo{\scshape Theorem B}  {\it Let $r\in\N\cup\{\infty\}${\rm ,} let $F_r$
denote the free group on $r$ generators{\rm ,} and let $\lambda\colon
F_r\to\CB(\ell^2(F_r))$ be the left regular representation of $F_r$.  
Then
there exists a sequence of unitary representations $\pi_n\colon F_r\to
M_n(\C)$ such that for all $h_1,\dots,h_m\in F_r$ and}
$c_1,\dots,c_m\in\C$\/{\rm :}\/
\[
\lim_{n\to\infty} \Big\| \sum^m_{j=1}
c_j\pi_n(h_j)\Big\|=\Big\|\sum^m_{j=1} c_j\lambda(h_j)\Big\|.
\]
\Enddemo

The invariant $\Ext(\CA)$ for separable unital $C^*$-algebras $\CA$ was
introduced by
Brown, Douglas and Fillmore in 1973 (cf.~\cite{BDF1},
\cite{BDF2}). $\Ext(\CA)$ is the set of equivalence classes $[\pi]$ of
one-to-one $*$-homomorphisms $\pi\colon \CA\to \CC(\CH)$, where
$\CC(\CH)=\CB(\CH)/\CK(\CH)$ is the Calkin algebra for the Hilbert space
$\CH=\ell^2(\N)$. The equivalence relation is defined as follows:
\[
\pi_1\sim\pi_2\iff \exists u\in\CU(\CB(\CH))\ \forall a\in\CA\colon
\pi_2(a)=\rho(u)\pi_1(a)\rho(u)^*,
\]
where $\CU(\CB(\CH))$ denotes the unitary group of $\CB(\CH)$ and
$\rho\colon \CB(\CH)\to \CC(\CH)$ is the quotient map. Since
$\CH\oplus \CH\simeq \CH$, the map $(\pi_1,\pi_2)\to\pi_1\oplus\pi_2$  
defines
a natural semi-group structure on $\Ext(\CA)$. By Choi and Effros
\cite{CE}, $\Ext(\CA)$ is a group for every separable unital nuclear
$C^*$-algebra and by Voiculescu \cite{V0}, $\Ext(\CA)$ is a unital
semi-group for all separable
unital $C^*$-algebras $\CA$. Anderson \cite{An} provided in 1978 the  
first
example of a unital $C^*$-algebra $\CA$ for which $\Ext(\CA)$ is not a
group. The $C^*$-algebra $\CA$ in \cite{An} is generated by the reduced
$C^*$-algebra $\Cred(F_2)$ of the free group $F_2$ on 2 generators
and a projection $p\in \CB(\ell^2(F_2))$. Since then, it has been an  
open
problem whether $\Ext(\Cred(F_2))$ is a group. In
\cite[Sect.~5.14]{V5}, Voiculescu shows that if one could prove  
Theorem~B,
then it would follow that $\Ext(\Cred(F_r))$ is not a group for any  
$r\ge
2$. Hence we have

\demo{\scshape Corollary 1}  {\it Let $r\in\N\cup\{\infty\}$, $r\ge 2$. Then
   $\Ext(\Cred(F_r))$ is not a group.}
\Enddemo

The problem of proving Corollary 1 has
been considered by a number of mathematicians; see \cite[\S
5.11]{V5} for
a more detailed discussion.

In Section 9 we extend Theorem A (resp.~Theorem B)
to polynomials (resp.~linear combinations) with coefficients in an
arbitrary unital exact $C^*$-algebra. The first of these two results is
used to provide new proofs of two key results from our previous paper
\cite{HT2}: ``Random matrices and $K$-theory for exact
$C^*$-algebras''. Moreover, we use the second result to make an exact
computation of the constants $C(r)$, $r\in\N$, introduced by Junge and
Pisier \cite{JP} in connection with their proof of
\[
\CB(\CH)\mathop{\otimes}_{\max} \CB(\CH)\ne
\CB(\CH)\mathop{\otimes}_{\min} \CB(\CH).
\]
Specifically, we prove the following:

\demo{\scshape Corollary 2}    {\it Let $r\in\N${\rm ,} $r\ge 2${\rm ,} and let $C(r)$ be the
infimum of all real numbers
$C>0$ with the following property\/{\rm :}\/ There exists a sequence of natural
numbers $(n(m))_{m\in\N}$ and a sequence of $r$\/{\rm -}\/tuples
$(u_1^{(m)},\dots,u_r^{(m)})_{m\in\N}$ of $n(m)\times n(m)$ unitary
matrices{\rm ,} such that
\[
\Big\|\sum^r_{i=1} u_i^{(m)}\otimes\bar{u}_i^{(m')}\Big\|\le C,
\]
whenever $m,m'\in\N$ and $m\ne m'$. Then $C(r)=2\sqrt{r-1}$. }
\Enddemo

Pisier proved in \cite{P3} that $C(r)\ge 2\sqrt{r-1}$ and Valette
proved subsequently in \cite{V} that $C(r)=2\sqrt{r-1}$, when $r$ is of
the form $r=p+1$ for an odd prime number $p$.

We end Section~9 by using Theorem A to prove the following result on
powers of ``circular'' random matrices (cf.\ \S 9):

\sdemo{Corollary 3}  {\it Let $Y$ be a random matrix in the class
   $\GRM(n,\frac 1n)${\rm ,} i.e.{\rm ,} the
entries of $Y$ are independent and identically distributed  complex Gaussian random variables
with density $z\mapsto\frac{n}{\pi}\e^{-n|z|^2}${\rm ,} $z\in\C$. Then for  
every
$p\in\N$ and almost all} $\omega\in\Omega$,
\[
\lim_{n\to\infty}\big\|Y(\omega)^p\big\|
=\bigg(\frac{(p+1)^{p+1}}{p^p}\bigg)^\frac12.
\]
\Enddemo

Note that for $p=1$, Corollary 3 follows from Geman's result
\cite{Ge}.

In the remainder of this introduction, we sketch the main steps in the
proof of Theorem A. Throughout the paper, we denote by $\CA_{\rm sa}$
the real vector space of self-adjoint elements in a $C^*$-algebra $\CA$.
In Section 2 we prove the following ``linearization trick'':

\medbreak
{\it Let $\CA,\CB$ be unital $C^*$\/{\rm -}\/algebras{\rm ,} and let $x_1,\dots,x_r$ and
$y_1,\dots,y_r$ be operators in $\CA_{\rm sa}$ and $\CB_{\rm sa}${\rm ,}
respectively. Assume that for all $m\in\N$ and all matrices  
$a_0,\dots,a_r$
in $\Mmsa${\rm ,} we have
\[
\spe\big(a_0\otimes\unit_\CB + \textstyle{\sum^r_{i=1} a_i\otimes
   y_i}\big)\subseteq
\spe\big(a_0\otimes\unit_\CA + \sum^r_{i=1} a_i\otimes x_i\big),
\]
where $\spe(T)$ denotes the spectrum of an operator $T${\rm ,} and where
$\unit_{\CA}$ and $\unit_{\CB}$ denote the units of $\CA$ and $\CB${\rm ,}
respectively. Then there
exists a unital $*$\/{\rm -}\/homo\-morphism
\[
\Phi\colon C^*(x_1,\dots,x_r,\unit_{\CA})\to  
C^*(y_1,\dots,y_r,\unit_\CB),
\]
such that $\Phi(x_i)=y_i${\rm ,} $i=1,\dots,r$. In particular{\rm ,}
\[
\|p(y_1,\dots,y_r)\|\le\|p(x_1,\dots,x_r)\|,
\]
for every polynomial $p$ in $r$ noncommuting variables}.
\Enddemo

The linearization trick allows us to conclude (see \S\ref{sec6}):
\sdemo{Lemma 1} {\it In order to prove Theorem {\rm A,} it is sufficient to
   prove the following\/{\rm :}\/ With
$(X_1^{(n)},\dots,X_r^{(n)})$ and $(x_1,\dots,x_r)$ as in Theorem
{\rm A,} one has for all $m\in\N${\rm ,} all matrices $a_0,\dots,a_r$ in $\Mmsa$
and all $\varepsilon >0$ that
\[
\spe\big(a_0\otimes\unit_n + \textstyle{\sum^r_{i=1} a_i\otimes
X_i^{(n)}(\omega)}\big)\subseteq\spe(a_0\otimes\unit_\CB +
\sum^r_{i=1}a_i\otimes x_i\big)+{}]-\ep,\ep[,
\]
eventually as $n\to\infty${\rm ,} for almost all $\omega\in\Omega${\rm ,} and where
$\unit_n$ denotes the unit of $M_n(\C)$.}
\Enddemo

In the rest of this section, $(X_1^{(n)},\dots,X_r^{(n)})$ and
$(x_1,\dots,x_r)$ are defined as in Theorem A. Moreover we let
$a_0,\dots,a_r\in \Mmsa$ and put
\begin{eqnarray*}
s &=& a_0\otimes\unit_{\CB} + \sum^r_{i=1} a_i\otimes x_i\\
S_n &=& a_0\otimes\unit_n+\sum^r_{i=1} a_i\otimes X_i^{(n)},\quad  
n\in\N.
\end{eqnarray*}
It was proved by Lehner in \cite{Le} that Voiculescu's $R$-transform of  
$s$
with amalgamation over $M_m(\C)$ is given by
\begin{equation}
\label{eq0-4}
\CR_s(z) = a_0+\sum^r_{i=1} a_i z a_i, \quad z\in M_m(\C).
\end{equation}
For $\lambda\in M_m(\C)$, we let $\im\lambda$ denote the self-adjoint
matrix $\im\lambda = \frac{1}{2\i}(\lambda-\lambda^*)$, and we put
\[
\CO = \big\{\lambda\in M_m(\C)\mid\im \lambda\ \mbox{is positive
definite}\big\} .
\]
 From \eqref{eq0-4} one gets (cf.~\S 6) that the matrix-valued
Stieltjes transform of $s$,
\[
G(\lambda) =  
(\id_m\otimes\tau)\big[(\lambda\otimes\unit_\CB-s)^{-1}\big]\in
M_m(\C),
\]
is defined for all $\lambda\in\CO$, and satisfies the matrix equation
\begin{equation}
\label{eq0-5}
\sum^r_{i=1} a_iG(\lambda)a_i  
G(\lambda)+(a_0-\lambda)G(\lambda)+\unit_m=0 .
\end{equation}
For $\lambda\in\CO$, we let $H_n(\lambda)$ denote the $M_m(\C)$-valued
random variable
\[
H_n(\lambda) = (\id_m\otimes\tr_n)\big[(\lambda\otimes
\unit_n-S_n)^{-1}\big],
\]
and we put
\[
G_n(\lambda) = \E\big\{H_n(\lambda)\big\}\in M_m(\C).
\]
Then the following analogy to \eqref{eq0-5} holds
(cf.~\S 3):

\sdemo{Lemma 2 {\rm (Master equation)}}  {\it For all $\lambda\in\CO$ and}  
$n\in\N$\/{\rm :}\/
\begin{equation}
\label{eq0-6}
\E\Big\{\sum^r_{i=1} a_iH_n(\lambda)a_iH_n(\lambda)
+(a_0-\lambda)H_n(\lambda)+\unit_m\Big\}=0.
\end{equation}
\Enddemo

The proof of \eqref{eq0-6} is completely different from the proof of
\eqref{eq0-5}. It is based on the simple observation that the density of
the standard Gaussian distribution,
$\varphi(x)=\frac{1}{\sqrt{2\pi}}\e^{-x^2/2}$ satisfies the first order
differential equation $\varphi'(x)+x\varphi(x)=0$. In the special case  
of a
single $\SGRM(n,\frac{1}{n})$ random matrix (i.e., $r=m=1$ and
$a_0=0,a_1=1$), equation \eqref{eq0-6} occurs in a recent paper by  
Pastur
(cf.\ \cite[Formula (2.25)]{pas}). Next we use the
so-called ``Gaussian Poincar\'e inequality'' (cf.\ \S 4) to estimate
the norm of the difference
\[
\E\Big\{\sum^r_{i=1} a_iH_n(\lambda)a_iH_n(\lambda)\Big\} - \sum^r_{i=1}
a_i\E\{H_n(\lambda)\}a_i\E\{H_n(\lambda)\},
\]
and we obtain thereby (cf.~\S 4):

\sdemo{Lemma 3  {\rm (Master inequality)}}  {\it For all $\lambda\in\CO$ and  
all
   $n\in\N${\rm ,} we have
\begin{equation}
\label{eq0-7}
\Big\|\sum^r_{i=1}
a_iG_n(\lambda)a_iG_n(\lambda)-(a_0-\lambda)G_n(\lambda)+\unit_m\Big\|
\le\frac{C}{n^2}\big\|(\im\lambda)^{-1}\big\|^4,
\end{equation}
where $C=m^3\big\|\sum^r_{i=1} a^2_i\big\|^2$.}
\Enddemo

In Section 5, we deduce from \eqref{eq0-5} and \eqref{eq0-7} that
\begin{equation}
\label{eq0-8}
\|G_n(\lambda)-G(\lambda)\|\le\frac{4C}{n^2}
\big(K+\|\lambda\|\big)^2\big\|(\im\lambda)^{-1}\big\|^7,
\end{equation}
where $C$ is as above and $K=\|a_0\|+4\sum^r_{i=1} \|a_i\|$. The
estimate \eqref{eq0-8} implies that for every $\varphi\in
C^\infty_c(\R,\R)$:
\begin{equation}
\label{eq0-9}
\E\big\{(\tr_m\otimes\tr_n)\varphi(S_n)\big\} =
(\tr_m\otimes\tau)(\varphi(s))+ O\big(\textstyle{\frac{1}{n^2}}\big),
\end{equation}
for $n\to\infty$ (cf.\ \S 6). Moreover, a second application of the
Gaussian Poincar\'e inequality yields that
\begin{equation}
\label{eq0-10}
\V\big\{(\tr_m\otimes\tr_n)\varphi(S_n)\big\}
\le\frac{1}{n^2}\E\big\{(\tr_m\otimes\tr_n)(\varphi'(S_n)^2)\big\},
\end{equation}
where $\V$ denotes the variance. Let now $\psi$ be a $C^\infty$-function
with values in $[0,1]$, such that $\psi$ vanishes on a neighbourhood of
the spectrum $\spe(s)$ of $s$, and such that $\psi$ is 1 on the
complement of $\spe(s)+{}]-\varepsilon,\varepsilon[$.

By applying \eqref{eq0-9} and \eqref{eq0-10} to $\varphi=\psi-1$, one  
gets
\begin{eqnarray*}
\E\big\{(\tr_m\otimes\tr_n)\psi(S_n)\big\} &=&
O(n^{-2}),\\[.2cm]
\V\big\{(\tr_m\otimes\tr_n)\psi(S_n)\big\} &=&
O(n^{-4}),
\end{eqnarray*}
and by a standard application of the Borel-Cantelli lemma, this implies
that
\[
(\tr_m\otimes\tr_n)\psi(S_n(\omega))=O(n^{-4/3}),
\]
for almost all $\omega\in\Omega$. But the number of eigenvalues of
$S_n(\omega)$ outside $\spe(s)+{}]-\varepsilon,\varepsilon[$ is
dominated by $mn(\tr_m\otimes\tr_n)\psi(S_n(\omega))$, which is
$O(n^{-1/3})$ for $n\to\infty$. Being an integer, this number must
therefore vanish eventually as $n\to\infty$, which shows that for almost
all $\omega\in\Omega$,
\[
\spe(S_n(\omega))\subseteq\spe(s)+{}]-\varepsilon,\varepsilon[,
\]
eventually as $n\to\infty$, and Theorem A now follows from
Lemma 1.

\section{A linearization trick}\label{sec1}

Throughout this section we consider two unital $C^*$-algebras $\CA$ and
$\CB$ and self-adjoint elements $x_1,\dots,x_r\in\CA$,
$y_1,\dots,y_r\in\CB$. We put
\[
\CA_0=C^*(\unit_{\CA},x_1,\dots,x_r) \quad \textrm{and} \quad
\CB_0=C^*(\unit_{\CB},y_1,\dots,y_r).
\]
Note that since $x_1,\dots,x_r$ and
$y_1,\dots,y_r$ are self-adjoint, the complex linear spaces
\[
E=\textrm{span}_{\C}\{\unit_{\CA},x_1,\dots,x_r,\textstyle{\sum_{i=1}^r 
x_i^2}\}
\enspace
\textrm{and}
\enspace  
F=\textrm{span}_{\C}\{\unit_{\CB},y_1,\dots,y_r,\sum_{i=1}^ry_i^2\}
\]
are both operator systems.

\begin{lemma}\label{lin1}
Assume that
$u_0\colon E\to F$ is a unital
completely positive \/{\rm (}\/linear\/{\rm )}\/ mapping{\rm ,} such that
\[
u_0(x_i)=y_i, \qquad i=1,2,\dots,r,
\]
and
\[
u_0\big(\textstyle{\sum_{i=1}^rx_i^2}\big)=\sum_{i=1}^ry_i^2.
\]
Then there exists a surjective $*$\/{\rm -}\/homomorphism $u\colon\CA_0\to\CB_0${\rm ,}
such that
\[
u_0=u_{\mid E}.
\]
\end{lemma}

\Proof   The proof is inspired by Pisier's proof of \cite[Prop.~1.7]{P2}.  
We
may assume that $\CB$ is a unital sub-algebra of $\CB(\CH)$ for
some Hilbert space~$\CH$. Combining Stinespring's theorem
(\cite[Thm.~4.1]{pa}) with Arveson's extension theorem
(\cite[Cor.~6.6]{pa}), it follows then that there exists a Hilbert
space $\CK$ containing $\CH$, and a unital $*$-homomorphism
$\pi\colon\CA\to\CB(\CK)$, such that
\[
u_0(x)=p\pi(x)p \qquad (x\in E),
\]
where $p$ is the orthogonal projection of $\CK$ onto $\CH$. Note in
particular that

\begin{itemize}

\item[(a)] $u_0(\unit_{\CA})=p\pi(\unit_{\CA})p=p=\unit_{\CB(\CH)}$,

\item[(b)] $y_i=u_0(x_i)=p\pi(x_i)p$, \ $i=1,\dots,r$,

\item[(c)]
   $\sum_{i=1}^ry_i^2=u_0\big(\sum_{i=1}^rx_i^2\big)
=\sum_{i=1}^rp\pi(x_i)^2p$.

\end{itemize}

 From (b) and (c), it follows that $p$ commutes with $\pi(x_i)$
for all $i$ in $\{1,2,\dots,r\}$. Indeed, using (b) and (c), we find  
that
\[
\sum_{i=1}^rp\pi(x_i)p\pi(x_i)p=\sum_{i=1}^ry_i^2=\sum_{i=1}^rp\pi(x_i)^ 
2p,
\]
so that
\[
\sum_{i=1}^rp\pi(x_i)\big(\unit_{\CB(\CK)}-p\big)\pi(x_i)p=0.
\]
Thus, putting $b_i=(\unit_{\CB(\CK)}-p)\pi(x_i)p$, $i=1,\dots,r$, we  
have
that $\sum_{i=1}^rb_i^*b_i=0$, so that $b_1=\cdots=b_r=0$. Hence, for  
each
$i$ in $\{1,2,\dots,r\}$, we have
\begin{eqnarray*}
[p,\pi(x_i)] 
&=&p\pi(x_i)-\pi(x_i)p\\
&=&
p\pi(x_i)(\unit_{\CB(\CK)}-p)-(\unit_{\CB(\CK)}-p)\pi(x_i)p= b_i^*-b_i  
= 0,
\end{eqnarray*}
as desired. Since $\pi$ is a unital $*$-homomorphism, we may conclude
further that $p$ commutes with all elements of the $C^*$-algebra
$\pi(\CA_0)$.

Now define the mapping $u\colon\CA_0\to\CB(\CH)$ by
\[
u(a)=p\pi(a)p, \quad (a\in\CA_0).
\]
Clearly $u(a^*)=u(a)^*$ for all $a$ in $\CA_0$, and, using (a) above,
$u(\unit_{\CA})=u_0(\unit_{\CA})\break=\unit_{\CB}$. Furthermore, since $p$
commutes with $\pi(\CA_0)$, we find for any $a,b$ in $\CA_0$ that
\[
u(ab)=p\pi(ab)p=p\pi(a)\pi(b)p=p\pi(a)p\pi(b)p=u(a)u(b).
\]
Thus, $u\colon\CA_0\to\CB(\CH)$ is a unital $*$-homomorphism, which  
extends
$u_0$, and $u(\CA_0)$ is a $C^*$-sub-algebra of $\CB(\CH)$. It remains  
to note
that $u(\CA_0)$ is generated, as a $C^*$-algebra, by the set
$u(\{\unit_{\CA},x_1,\dots,x_r\})=\{\unit_{\CB},y_1,\dots,y_r\}$, so  
that
$u(\CA_0)=C^*(\unit_{\CB},y_1,\dots,y_r)=\CB_0$, as desired.
\Endproof\vskip4pt 

For any element $c$ of a $C^*$-algebra $\CC$, we denote by $\spe(c)$ the
spectrum of~$c$, i.e.,
\[
\spe(c)=\{\lambda\in\C\mid c-\lambda\unit_{\CC} \ \textrm{is not
   invertible}\}.
\]

\begin{theorem}\label{lin2}
Assume that the self\/{\rm -}\/adjoint elements $x_1,\dots,x_r\in\CA$ and
$y_1,\dots,y_r\in\CB$ satisfy the property\/{\rm :}\/
\begin{multline}\label{e1.1}
\forall m\in\N \ \forall a_0,a_1,\dots,a_r\in\Mmsa\colon \\
 \spe\big(a_0\otimes\unit_{\CA}+\textstyle{\sum_{i=1}^ra_i\otimes  
x_i}\big)
\supseteq  
\spe\big(a_0\otimes\unit_{\CB}+\textstyle{\sum_{i=1}^ra_i\otimes
   y_i}\big).
\end{multline}
Then there exists a unique surjective unital $*$\/{\rm -}\/homomorphism
$\varphi\colon\CA_0\to\CB_0${\rm ,} such that
\[
\varphi(x_i)=y_i, \qquad i=1,2,\dots,r.
\]
\end{theorem}

Before the proof of Theorem~\ref{lin2}, we make a few observations:

\begin{remark}\label{obs}
(1) \ In connection with condition \eqref{e1.1} above, let $V$ be a  
subspace
     of $M_m(\C)$ containing the unit $\unit_m$. Then the condition:
\begin{multline}
\label{e1.3}
\forall a_0, a_1,\dots,a_r\in V\colon \\ 
 \spe\big(a_0\otimes\unit_{\CA}+\textstyle{\sum_{i=1}^ra_i\otimes  
x_i}\big)
\supseteq  
\spe\big(a_0\otimes\unit_{\CB}+\textstyle{\sum_{i=1}^ra_i\otimes
   y_i}\big)
\end{multline}
is equivalent to the condition:
\begin{eqnarray}
\label{e1.4}
\forall a_0, a_1,\dots,a_r\in V\colon &&
 a_0\otimes\unit_{\CA}+\textstyle{\sum_{i=1}^ra_i\otimes x_i} \  
\textrm{is
   invertible}\\
&& \Longrightarrow
a_0\otimes\unit_{\CB}+\textstyle{\sum_{i=1}^ra_i\otimes y_i} \  
\textrm{is
   invertible}.\nonumber
\end{eqnarray}
Indeed, it is clear that \eqref{e1.3} implies \eqref{e1.4}, and the  
reverse
implication follows by replacing, for any complex number $\lambda$, the
matrix $a_0\in V$ by $a_0-\lambda\unit_m\in V$.

(2) \ Let $\CH_1$ and $\CH_2$ be Hilbert spaces and consider the Hilbert
     space direct sum $\CH=\CH_1\oplus\CH_2$. Consider further the  
operator
     $R$ in $\CB(\CH)$ given in matrix form as
\[
R=
\begin{pmatrix} x &  y \\
          z  & \unit_{\CB(\CH_2)},
\end{pmatrix}
\]
where $x\in\CB(\CH_1), y\in\CB(\CH_2,\CH_1)$ and $z\in\CB(\CH_1,\CH_2)$.
Then $R$ is invertible in $\CB(\CH)$ if and only if $x-yz$ is  
invertible in
$\CB(\CH_1)$.

This follows immediately by writing
\[
\begin{pmatrix} x &  y \\
          z  & \unit_{\CB(\CH_2)}
\end{pmatrix}
=
\begin{pmatrix} \unit_{\CB(\CH_1)} &  y \\
          0  & \unit_{\CB(\CH_2)}
\end{pmatrix}
\cdot
\begin{pmatrix} x-yz &  0 \\
          0  & \unit_{\CB(\CH_2)}
\end{pmatrix}
\cdot
\begin{pmatrix} \unit_{\CB(\CH_1)} & 0 \\
          z  & \unit_{\CB(\CH_2)}
\end{pmatrix}, 
\]\pagebreak
where the first and last matrix on the right-hand side are invertible  
with
inverses given by:
\begin{eqnarray*}
\begin{pmatrix} \unit_{\CB(\CH_1)} &  y \\
          0  & \unit_{\CB(\CH_2)}
\end{pmatrix}^{-1}
&=&
\begin{pmatrix} \unit_{\CB(\CH_1)} & -y \\
           0  & \unit_{\CB(\CH_2)}
\end{pmatrix}\\ 
\noalign{\noindent 
 and} 
\begin{pmatrix} \unit_{\CB(\CH_1)} & 0 \\
          z  & \unit_{\CB(\CH_2)}
\end{pmatrix}^{-1}
&=&
\begin{pmatrix} \unit_{\CB(\CH_1)} & 0 \\
           -z & \unit_{\CB(\CH_2)}
\end{pmatrix}.
\end{eqnarray*}
\end{remark}

\demo{Proof of Theorem~{\rm \ref{lin2}}}  By Lemma~\ref{lin1}, our objective  
is to
   prove the existence of a
unital completely positive map $u_0\colon E\to F$, satisfying that
$u_0(x_i)=y_i$, $i=1,2,\dots,r$ and
$u_0(\sum_{i=1}^rx_i^2)=\sum_{i=1}^ry_i^2$.

\demo{Step {\rm I}} We show first that the assumption \eqref{e1.1} is  
equivalent
to the seemingly stronger condition:
\begin{multline}
\label{e1.2}
\forall m\in\N \ \forall a_0, a_1,\dots,a_r\in M_m(\C)\colon \\ 
 \spe\big(a_0\otimes\unit_{\CA}+\textstyle{\sum_{i=1}^ra_i\otimes  
x_i}\big)
\supseteq  
\spe\big(a_0\otimes\unit_{\CB}+\textstyle{\sum_{i=1}^ra_i\otimes
   y_i}\big).
\end{multline}
Indeed, let $a_0,a_1,\dots,a_r$ be arbitrary matrices in $M_m(\C)$ and
consider then the self-adjoint matrices
$\tilde{a}_0,\tilde{a}_1,\dots,\tilde{a}_r$ in $M_{2m}(\C)$ given by:
\[
\tilde{a}_i=
\begin{pmatrix} 0 &  a_i^* \\
          a_i  & 0
\end{pmatrix},
\qquad i=0,1,\dots,r.
\]
Note then that
\begin{eqnarray*}
&&\hskip-6pt \tilde{a}_0\otimes\unit_{\CA}+\sum_{i=1}^r\tilde{a}_i\otimes x_i \\
&&\qquad =
\begin{pmatrix} 0 &  a_0^*\otimes\unit_{\CA}+\sum_{i=1}^ra_i^*\otimes  
x_i \\
          a_0\otimes\unit_{\CA}+\sum_{i=1}^ra_i\otimes x_i  & 0
\end{pmatrix} \\
&&\qquad =
\begin{pmatrix} 0 &  \unit_{\CA} \\
          \unit_{\CA}  & 0
\end{pmatrix}
\cdot
\begin{pmatrix} a_0\otimes\unit_{\CA}+\sum_{i=1}^ra_i\otimes x_i & 0\\
          0 & a_0^*\otimes\unit_{\CA}+\sum_{i=1}^ra_i^*\otimes x_i
\end{pmatrix}.
\end{eqnarray*}
Therefore, $\tilde{a}_0\otimes\unit_{\CA}+\sum_{i=1}^r\tilde{a}_i\otimes
x_i$ is invertible in $M_{2m}(\CA)$ if and only if
$a_0\otimes\unit_{\CA}+\sum_{i=1}^ra_i\otimes x_i$
is invertible in $M_m(\CA)$, and similarly,
of course, $\tilde{a}_0\otimes\unit_{\CB}+\sum_{i=1}^r\tilde{a}_i\otimes
y_i$ is invertible in $M_{2m}(\CB)$ if and only if
$a_0\otimes\unit_{\CB}+\sum_{i=1}^ra_i\otimes y_i$
is invertible in $M_m(\CB)$. It follows that
\[
\begin{split}
a_0\otimes\unit_{\CA}+\sum_{i=1}^ra_i\otimes x_i \ \textrm{is  
invertible}
&\iff
\tilde{a}_0\otimes\unit_{\CA}+\sum_{i=1}^r\tilde{a}_i\otimes x_i \
\textrm{is invertible}
\\[.2cm]
&\implies
\tilde{a}_0\otimes\unit_{\CB}+\sum_{i=1}^r\tilde{a}_i\otimes y_i \
\textrm{is invertible}
\\[.2cm]
&\iff
a_0\otimes\unit_{\CB}+\sum_{i=1}^ra_i\otimes y_i \ \textrm{is  
invertible},
\end{split}
\]
where the second implication follows from the assumption
\eqref{e1.1}. Since the argument above holds for arbitrary matrices
$a_0,a_1,\dots,a_r$ in $M_m(\C)$, it follows from Remark~\ref{obs}(1)  
that
condition \eqref{e1.2} is satisfied.

\demo{Step {\rm II}} We prove next that the assumption \eqref{e1.1} implies  
the
condition:
\begin{equation}
\begin{split}
\forall m\in\N \ \forall a_0,a_1,\dots,a_r,a_{r+1}&\in M_m(\C)\colon
\\[.2cm]
\spe\big(a_0\otimes\unit_{\CA}+&\textstyle{\sum_{i=1}^ra_i\otimes
   x_i}+a_{r+1}\otimes\sum_{i=1}^rx_i^2\big) \\[.2cm]
&\supseteq  
\spe\big(a_0\otimes\unit_{\CB}+\textstyle{\sum_{i=1}^ra_i\otimes
   y_i} + a_{r+1}\otimes\sum_{i=1}^ry_i^2\big).
\label{e1.5}
\end{split}
\end{equation}
Using Remark~\ref{obs}(1), we have to show, given $m$ in $\N$ and
$a_0,a_1,\dots,a_{r+1}$ in $M_m(\C)$, that invertibility of
$a_0\otimes\unit_{\CA}+\textstyle{\sum_{i=1}^ra_i\otimes
   x_i}+a_{r+1}\otimes\sum_{i=1}^rx_i^2$ in $M_m(\CA)$ implies  
invertibility
of $a_0\otimes\unit_{\CA}+\textstyle{\sum_{i=1}^ra_i\otimes
   y_i}+a_{r+1}\otimes\sum_{i=1}^ry_i^2$ in $M_m(\CB)$.
For this, consider the matrices:

 \begin{small}
\[
S=
\begin{pmatrix} a_0\otimes\unit_{\CA} &  -\unit_m\otimes x_1 &  
-\unit_m\otimes x_2
& \cdots & -\unit_m\otimes
x_r \\
a_1\otimes\unit_{\CA}+a_{r+1}\otimes x_1 & \unit_m\otimes\unit_{\CA} &  
\ &
\ & O \\
a_2\otimes\unit_{\CA}+a_{r+1}\otimes x_2 & \ &
\unit_m\otimes\unit_{\CA} & \ & \ \\
\vdots & \ & \ & \ddots & \ \\
a_r\otimes\unit_{\CA}+a_{r+1}\otimes x_r & O & \ & \ &
\unit_m\otimes\unit_{\CA}
\end{pmatrix}
\in M_{(r+1)m}(\CA)
\] \end{small}

\noindent 
and 

\begin{small}
\[
T=
\begin{pmatrix} a_0\otimes\unit_{\CB} &  -\unit_m\otimes y_1 &  
-\unit_m\otimes y_2
& \cdots & -\unit_m\otimes
y_r \\
a_1\otimes\unit_{\CB}+a_{r+1}\otimes y_1 & \unit_m\otimes\unit_{\CB} &  
\ &
\ & O \\
a_2\otimes\unit_{\CB}+a_{r+1}\otimes y_2 & \ &
\unit_m\otimes\unit_{\CB} & \ & \ \\
\vdots & \ & \ & \ddots & \ \\
a_r\otimes\unit_{\CB}+a_{r+1}\otimes y_r & O & \ & \ &
\unit_m\otimes\unit_{\CB}
\end{pmatrix}
\in M_{(r+1)m}(\CB).
\]
\end{small}

By Remark~\ref{obs}(2), invertibility of $S$ in $M_{(r+1)m}(\CA)$ is
equivalent to invertibility of
\[
\begin{split}
a_0\otimes\unit_{\CA}+\textstyle{\sum_{i=1}^r(\unit_m\otimes
   x_i)}\cdot(&a_i\otimes\unit_{\CA} + a_{r+1}\otimes x_i) \\[.2cm]
&=a_0\otimes\unit_{\CA}+\textstyle{\sum_{i=1}^ra_i\otimes
   x_i}+a_{r+1}\otimes\sum_{i=1}^rx_i^2
\end{split}
\]
in $M_m(\CA)$. Similarly, $T$ is invertible in $M_{(r+1)m}(\CB)$ if and
only if
\[
a_0\otimes\unit_{\CB}+\textstyle{\sum_{i=1}^ra_i\otimes
   y_i}+a_{r+1}\otimes\sum_{i=1}^ry_i^2
\]
is invertible in $M_m(\CB)$. It remains thus to
show that invertibility of $S$ implies that of $T$. This, however,  
follows
immediately from Step~I, since we may write $S$ and $T$ in the form:
\[
S=b_0\otimes\unit_{\CA}+\sum_{i=1}^rb_i\otimes x_i \quad \textrm{and}  
\quad
T=b_0\otimes\unit_{\CB}+\sum_{i=1}^rb_i\otimes y_i,
\]
for suitable matrices $b_0,b_1,\dots,b_r$ in $M_{(r+1)m}(\C)$; namely
\[
b_0=
\begin{pmatrix} a_0 & 0 & 0 & \cdots & 0 \\
a_1 & \unit_m & \ & \ & \textrm{\large O} \\
a_2 & \ & \unit_m & \ & \ \\
\vdots & \ & \ & \ddots & \ \\
a_r & \textrm{\large O} & \ & \ & \unit_m
\end{pmatrix}
\]
and
\[
b_i=
\begin{pmatrix} 0 & \cdots & 0 & -\unit_m & 0 & \cdots & 0 \\
\vdots & \ & \ & \ & \ & \ & \  \\
0 & \ & \ & \ & \ & \ & \ \\
a_{r+1} & \ & \ & \textrm{\Huge O}  &  & \ & \ \\
0 & \ & \ & \ & \ & \ & \ \\
\vdots & \ & \ & \ & \ & \ & \  \\
0 & \ & \ & \ & \ & \ & \
\end{pmatrix},
\qquad i=1,2,\dots,r.
\]
For $i$ in $\{1,2,\dots,r\}$, the (possible) nonzero entries in
$b_i$ are at positions\break $(1,i+1)$ and $(i+1,1)$. This concludes
Step~II.

\demo{Step~{\rm III}} We show, finally, the existence of a unital completely
positive mapping $u_0\colon E\to F$, satisfying that
$u_0(x_i)=y_i$, $i=1,2,\dots,r$ and
$u_0(\sum_{i=1}^rx_i^2)=\sum_{i=1}^ry_i^2$.
\vglue3pt

Using Step~II in the case $m=1$, it follows that for any complex numbers
$a_0,a_1,\dots,a_{r+1}$, we have that
\begin{multline}
\spe\big(a_0\unit_{\CA}+\textstyle{\sum_{i=1}^ra_ix_i} +
a_{r+1}\sum_{i=1}^rx_i^2\big) \\ \supseteq
\spe\big(a_0\unit_{\CB}+\textstyle{\sum_{i=1}^ra_iy_i} +
a_{r+1}\sum_{i=1}^ry_i^2\big).
\label{e1.6}
\end{multline}
If $a_0,a_1,\dots,a_{r+1}$ are real numbers, then the operators
\[
a_0\unit_{\CA}+\textstyle{\sum_{i=1}^ra_ix_i} + a_{r+1}\sum_{i=1}^rx_i^2
\quad \textrm{and} \quad
a_0\unit_{\CB}+\textstyle{\sum_{i=1}^ra_iy_i} + a_{r+1}\sum_{i=1}^ry_i^2
\]
are self-adjoint, since $x_1,\dots,x_r$ and $y_1,\dots,y_r$ are  
self-adjoint.
Hence \eqref{e1.6} implies that
\begin{equation}
\begin{split}
\forall a_0,&\ldots,a_{r+1}\in\R\colon \\[.2cm]
&\big\|a_0\unit_{\CA}+\textstyle{\sum_{i=1}^ra_ix_i} +
a_{r+1}\sum_{i=1}^rx_i^2\big\| \ge
\big\|a_0\unit_{\CB}+\textstyle{\sum_{i=1}^ra_iy_i} +
a_{r+1}\sum_{i=1}^ry_i^2\big\|.
\label{e1.7}
\end{split}
\end{equation}
Let \pagebreak $E'$ and $F'$ denote, respectively, the $\R$-linear span of
$\{\unit_{\CA},x_1,\dots,x_r,\sum_{i=1}^rx_i^2\}$ and
$\{\unit_{\CB},y_1,\dots,y_r,\sum_{i=1}^ry_i^2\}$:
\begin{eqnarray*}
E'&=&\textrm{span}_{\R}\{\unit_{\CA},x_1,\dots,x_r,\textstyle{\sum_{i=1}^ 
rx_i^2}\}
\\
\noalign{\noindent and} 
F'&=&\textrm{span}_{\R}\{\unit_{\CB},y_1,\dots,y_r,\sum_{i=1}^ry_i^2\}.
\end{eqnarray*}
It follows then from \eqref{e1.7} that there is a (well-defined)
$\R$-linear mapping $u_0'\colon E'\to F'$ satisfying that
$u_0'(\unit_{\CA})=\unit_{\CB}$, $u_0'(x_i)=y_i$,
$i=1,2,\dots,r$ and $u_0'(\sum_{i=1}^rx_i^2)=\sum_{i=1}^ry_i^2$.
For an arbitrary element $x$ in $E$, note that
$\re(x)=\frac{1}{2}(x+x^*)\in E'$ and ${\rm Im}(x)=\frac{1}{2\i}(x-x^*)\in
E'$. Hence, we may define a mapping $u_0\colon E\to F$ by setting:
\[
u_0(x)=u'_0(\re(x))+\i u_0'({\rm Im}(x)), \qquad (x\in E).
\]
It is straightforward, then, to check that $u_0$ is a $\C$-linear
mapping from $E$ onto~$F$, which extends $u_0'$.

Finally, it follows immediately from Step~II that for all $m$
in $\N$, the mapping $\id_{M_m(\C)}\otimes u_0$ preserves
positivity. In other words, $u_0$ is a completely positive
mapping. This concludes the proof.
\Endproof\vskip4pt 

In Section \ref{sec6}, we shall need the following strengthening of
Theorem~\ref{lin2}:

\begin{theorem}\label{thm1-4}
Assume that the self adjoint elements $x_1,\dots,x_r\in\CA${\rm ,}
$y_1,\dots,y_r\in\CB$ satisfy the property
\begin{multline}
\label{eq1-8}
\forall m\in\N \ \forall a_0,\dots,a_r\in M_m(\Q + \i\Q)_{\rm sa} :\\
  \spe\big(a_0\otimes\unit_A+\sum^r_{i=1} a_i\otimes  
x_i\big)\supseteq
\spe\big(a_0\otimes 1_B+\sum_{i=1}^ra_i\otimes y_i\big).
\end{multline}
Then there exists a unique surjective unital $*$\/{\rm -}\/homomorphism
$\varphi\colon A_0\to B_0$ such that $\varphi(x_i)=y_i${\rm ,} $i=1,\dots,r$.
\end{theorem}

\Proof   By Theorem~\ref{lin2}, it suffices to prove that condition
\eqref{eq1-8} is equivalent to condition \eqref{e1.1} of that theorem.
Clearly \eqref{e1.1} $\Rightarrow$ \eqref{eq1-8}. It remains to be  
proved
that \eqref{eq1-8}
$\Rightarrow$ \eqref{e1.1}. Let $d_H(K,L)$ denote the Hausdorff distance
between two subsets $K$, $L$ of $\C$:
\begin{equation}
\label{eq1-9}
d_H(K,L)=\max\Big\{\sup_{x\in K} d(x,L), \ \sup_{y\in L}
d(y,K)\Big\}.
\end{equation}
For normal operators $A,B$ in $M_m(\C)$ or $\CB(\CH)$ ($\CH$ a Hilbert
space) one has
\begin{equation}
\label{eq1-10}
d_H(\spe(A),\spe(B))\le\|A-B\|
\end{equation}
(cf. \cite[Prop.~2.1]{Da}). Assume now that \eqref{eq1-8} is satisfied,  
let
$m\in\N$, $b_0,\dots,b_r \in M_m(\C)$ and let $\varepsilon >0$.

Since $M_m(\Q+ \i\Q)_{\rm sa}$ is dense in $M_m(\C)_{\rm sa}$, we can  
choose
$a_0,\dots,a_r\in M_m(\Q+\i\Q)_{\rm sa}$ such that
\[
\|a_0-b_0\|+\sum^r_{i=1} \|a_i-b_i\| \| x_i\| <\varepsilon
\]
and
\[
\|a_0-b_0\|+\sum^r_{i=1} \|a_i-b_i\| \| y_i\| <\varepsilon.
\]
Hence, by \eqref{eq1-10},
\[
d_H\big(\spe\big(a_0\otimes 1 + \textstyle{\sum^r_{i=1}} a_i\otimes
x_i\big), \spe\big(b_0\otimes 1 + \sum^r_{i=1} b_i\otimes
x_i\big)\big)<\varepsilon
\]
and
\[
d_H\big(\spe\big(a_0\otimes 1 + \textstyle{\sum^r_{i=1}} a_i\otimes  
y_i\big),
\spe\big(b_0\otimes 1 + \sum^r_{i=1} b_i\otimes  
y_i\big)\big)<\varepsilon.
\]
By these two inequalities and \eqref{eq1-8} we get
\begin{eqnarray*}
\spe\big(b_0\otimes 1 + \textstyle{\sum^r_{i=1}}b_i\otimes y_i\big)
&\subseteq & \spe\big(a_0\otimes 1 + \textstyle{\sum^r_{i=1}}a_i\otimes
y_i\big) \ + \ ]-\varepsilon,\varepsilon[\\
&\subseteq &
\spe\big(a_0\otimes 1 + \textstyle{\sum^r_{i=1}} a_i\otimes
x_i\big) \ + \ ]-\varepsilon,\varepsilon[\\
&\subseteq &
\spe\big(b_0\otimes 1 + \textstyle{\sum^r_{i=1}}b_i\otimes
x_i) \ + \ ]-2\varepsilon,2\varepsilon[.
\end{eqnarray*}
Since $\spe(b_0\otimes 1 + \sum^r_{i=1} b_i\otimes y_i)$ is compact and
$\varepsilon >0$ is arbitrary, it follows that
\[
\spe\big(b_0\otimes 1 + \textstyle{\sum^r_{i=1}}b_i\otimes
y_i\big)\subseteq\spe\big(b_0\otimes 1 + \sum^r_{i=1} b_i\otimes  
x_i\big),
\]
for all $m\in\N$ and all $b_0,\dots,b_r\in M_m(\C)_{\rm sa}$, i.e.  
\eqref{e1.1}
holds. This completes the proof of Theorem~\ref{thm1-4}.
\hfill\qed

\section{The master equation} \label{sec2}

Let $\CH$ be a Hilbert space. For $T\in \CB(\CH)$ we let $\im T$ denote  
the
self adjoint operator $\im T = \frac{1}{2\i}(T-T^*)$. We say that a
matrix $T$ in $M_m(\C)_{\rm sa}$ is positive definite if all its  
eigenvalues
are strictly positive, and we denote by $\lmax(T)$ and $\lmin(T)$ the
largest and smallest eigenvalues of $T$, respectively.

\begin{lemma}\label{imaginaerdel}
\begin{itemize}
\item[{\rm (i)}] Let $\CH$ be a Hilbert space and let $T$ be an operator in
   $\CB(\CH)${\rm ,} such that the imaginary part $\im T$ satisfies one of
   the two conditions\/{\rm :}\/
\[
\im T\ge\varepsilon\unit_{\CB(\CH)} \qquad \textrm{or} \qquad
\im T\le-\varepsilon\unit_{\CB(\CH)},
\]
for some $\varepsilon$ in $]0,\infty[$. Then $T$ is invertible and
$\|T^{-1}\|\le\frac{1}{\varepsilon}$.

\item[{\rm (ii)}]
Let $T$ be a matrix in $M_m(\C)$ and assume that $\im T$ is positive
definite. Then $T$ is invertible and $\|T^{-1}\|\le \|(\im T)^{-1}\|$.
\end{itemize}
\end{lemma}

\Proof   Note first that (ii) is a special case of (i). Indeed, since
$\im T$ is self-adjoint, we have that $\im T\ge\lmin(\im T)\unit_m$.
Since $\im T$ is positive definite, $\lmin(\im T)>0$, and hence (i)
applies. Thus, $T$ is invertible and furthermore
\[
\|T^{-1}\|\le\frac{1}{\lmin(\im T)}=\lmax\big((\im T)^{-1}\big)=\|(\im
T)^{-1}\|,
\]
since $(\im T)^{-1}$ is positive.

To prove (i), note first that by replacing, if necessary, $T$ by
$-T$, it suffices to consider the case where
$\im T\ge\varepsilon\unit_{\CB(\CH)}$.
Let $\|\cdot\|$ and $\<\cdot,\cdot\>$ denote, respectively, the norm
and the inner product on $\CH$. Then, for any unit vector $\xi$ in
$\CH$, we have
\begin{eqnarray*}
\|T\xi\|^2&=&\|T\xi\|^2\|\xi\|^2\ge |\<T\xi,\xi\>|^2
\\
&=&\big|\<\re(T)\xi,\xi\>+\i\<\im T\xi,\xi\>\big|^2
\ge \<\im T\xi,\xi\>^2 \ge \varepsilon^2\|\xi\|^2,
\end{eqnarray*}
where we used that $\<\re(T)\xi,\xi\>,\<\im T\xi,\xi\>\in\R$. Note
further, for any unit vector $\xi$ in $\CH$, that
\[
\|T^*\xi\|^2\ge|\<T^*\xi,\xi\>|^2=|\<T\xi,\xi\>|^2\ge
\varepsilon^2\|\xi\|^2.
\]
Altogether, we have verified that $\|T\xi\|\ge\varepsilon\|\xi\|$ and
that $\|T^*\xi\|\ge\varepsilon\|\xi\|$ for any (unit) vector $\xi$ in
$\CH$, and by \cite[Prop.~3.2.6]{pe} this implies that $T$ is  
invertible and
that $\|T^{-1}\|\le\frac{1}{\varepsilon}$.
\hfill\qed

\begin{lemma}\label{Kaplansky trick}
Let $\CA$ be a unital $C^*$\/{\rm -}\/algebra and denote by $\GL(\CA)$ the group
of invertible elements of $\CA$. Let further $A\colon I\to\GL(\CA)$
be a mapping from an open interval $I$ in $\R$ into $\GL(\CA)${\rm ,} and
assume that $A$ is differentiable{\rm ,} in the sense that
\[
A'(t_0):=\lim_{t\to t_0}\frac{1}{t-t_0}\big(A(t)-A(t_0)\big)
\]
exists in the operator norm{\rm ,} for any $t_0$ in $I$. Then the mapping
$t\mapsto A(t)^{-1}$ is also differentiable and
\[
\frac{\rm d}{{\rm d}t}A(t)^{-1}=-A(t)^{-1}A'(t)A(t)^{-1}, \qquad (t\in  
I).
\]
\end{lemma}

\Proof   The lemma is well known. For the reader's convenience we include
a proof. For any $t,t_0$ in $I$, we have
\[
\begin{split}
\frac{1}{t-t_0}\big(A(t)^{-1}-A(t_0)^{-1}\big) &=
\frac{1}{t-t_0}A(t)^{-1}\big(A(t_0)-A(t)\big)A(t_0)^{-1} \\[.2cm]
&=-A(t)^{-1}\Big(\frac{1}{t-t_0}\big(A(t)-A(t_0)\big)\Big)A(t_0)^{-1}
\\[.2cm]
&\underset{t\to t_0}{\longrightarrow} -A(t_0)^{-1}A'(t_0)A(t_0)^{-1},
\end{split}
\]
where the limit is taken in the operator norm, and we use that the
mapping $B\mapsto B^{-1}$ is a homeomorphism of $\GL(\CA)$ with respect to  the
operator norm.
\phantom{endofline}\hfill\qed

\begin{lemma}\label{Gauss trick}
Let $\sigma$ be a positive number{\rm ,} let $N$ be a positive integer and
let $\gamma_1,\dots,\gamma_N$ be
$N$ independent identically distributed real valued random variables
with distribution $N(0,\sigma^2)${\rm ,} defined on the same probability
space $(\Omega,\CF,P)$.
Consider further a finite dimensional vector space $E$ and
a $C^1$\/{\rm -}\/mapping\/{\rm :}\/
\[
(x_1,\dots,x_N)\mapsto F(x_1,\dots,x_N)\colon\R^N\to E,
\]
satisfying that $F$ and all its first order partial derivatives
$\frac{\partial F}{\partial x_1},\dots,\frac{\partial F}{\partial
   x_N}$ are polynomially bounded.
For any $j$ in $\{1,2,\dots,N\}${\rm ,} we then have
\[
\E\big\{\gamma_jF(\gamma_1,\dots,\gamma_N)\big\} =
\sigma^2\E\big\{\textstyle{\frac{\partial F}{\partial
   x_j}}(\gamma_1,\dots,\gamma_N)\big\},
\]
where $\E$ denotes expectation with respect to  $P$.
\end{lemma}

\Proof  
Clearly it is sufficient to treat the case $E=\C$. The joint
distribution of $\gamma_1,\dots,\gamma_N$ is given by the density
function
\[
\varphi(x_1,\dots,x_N) = (2\pi\sigma^2)^{-\frac n2}
\exp\big(\textstyle{-\frac{1}{2\sigma^2}\sum^N_{i=1}x_i^2}\big), \qquad
(x_1,\dots,x_N)\in\R^N.
\]
Since
\[
\frac{\partial\varphi}{\partial x_j} (x_1,\dots,x_N) =
-\frac{1}{\sigma^2} x_j \varphi(x_1,\dots,x_N),
\]
we get by partial integration in the variable $x_j$,
\begin{eqnarray*}
\E\big\{\gamma_j F(\gamma_1,\dots,\gamma_N)\big\} &=& \int_{\R^N}
F(x_1,\dots,x_N)x_j\varphi(x_1,\dots,x_N)\6x_1,\dots,\6x_N\\
&=& -\sigma^2 \int_{\R^N}
F(x_1,\dots,x_N)\frac{\partial\varphi}{\partial x_j}
(x_1,\dots,x_N)\6x_1,\dots,\6x_N\\
&=& \sigma^2 \int_{\R^N}
\frac{\partial F}{\partial x_j}(x_1,\dots,x_N) \varphi
(x_1,\dots,x_N)\6x_1,\dots,\6x_N\\
&=& \sigma^2\E\bigg\{\frac{\partial F}{\partial x_j}
(\gamma_1,\dots,\gamma_N)\bigg\}.  
\end{eqnarray*}
\vglue-24pt \Endproof
\vglue16pt

Let $r$ and $n$ be positive integers. In the following we denote by
$\CE_{r,n}$
the {\it real} vector space $(\Mnsa)^r$. Note that $\CE_{r,n}$ is a
Euclidean space with inner product $\<\cdot,\cdot\>_e$ given by
\begin{multline*}
\<(A_1,\dots,A_r),(B_1,\dots,B_r)\>_{e}
\\=\Tr_n\Big(\sum_{j=1}^rA_jB_j\Big), \qquad
((A_1,\dots,A_r),(B_1,\dots,B_r)\in\CE_{r,n}),
\end{multline*}
and with norm given by
\[
\|(A_1,\dots,A_r)\|^2_{e} = \Tr_n\Big(\sum_{j=1}^rA_j^2\Big)
=\sum_{j=1}^r\|A_j\|^2_{2,\Tr_n},
\qquad ((A_,\dots,A_r)\in\CE_{r,n}).
\]
Finally, we shall denote by $S_1(\CE_{r,n})$ the unit sphere of  
$\CE_{r,n}$
with respect to  $\|\cdot\|_e$.

\begin{remark}\label{kanonisk iso}
Let $r,n$ be positive integers, and consider the linear isomorphism
$\Psi_0$ between $\Mnsa$ and $\R^{n^2}$ given by
\begin{equation}
\Psi_0((a_{kl})_{1\le k,l\le n})=\big((a_{kk})_{1\le k\le
n},(\sqrt{2}\re(a_{kl}))_{1\le k<l\le n},
(\sqrt{2}{\rm Im}(a_{kl}))_{1\le k<l\le n}\big),
\label{e2.0}
\end{equation}
for $(a_{kl})_{1\le k,l\le n}$ in $\Mnsa$.
We denote further by $\Psi$ the natural extension of $\Psi_0$ to a  
linear
isomorphism between $\CE_{r,n}$ and $\R^{rn^2}$:
\[
\Psi(A_1,\dots,A_r)=(\Psi_0(A_1),\dots,\Psi_0(A_r)), \qquad
(A_1,\dots,A_r\in\Mnsa).
\]
We shall identify $\CE_{r,n}$ with $\R^{rn^2}$ via the isomorphism
$\Psi$. Note that under this identification, the norm $\|\cdot\|_e$ on
$\CE_{r,n}$ corresponds to the usual Euclidean norm on $\R^{rn^2}$. In
other words, $\Psi$ is an isometry.

Consider next independent random matrices $X_1^{(n)},\dots,X_r^{(n)}$  
from\break
$\SGRM(n,\frac{1}{n})$ as defined in the introduction. Then  
$\X=(X_1^{(n)},\dots,X_r^{(n)})$ is a random
variable taking values in $\CE_{r,n}$, so that $\Y=\Psi(\X)$ is a random
variable taking values in $\R^{rn^2}$. From the definition of
$\SGRM(n,\frac{1}{n})$ and the fact that $X_1^{(n)},\dots,X_r^{(n)}$  
are
independent, it is easily seen that the distribution of $\Y$ on  
$\R^{rn^2}$
is the product measure $\mu=\nu\otimes\nu\otimes\cdots\otimes\nu$  
($rn^2$
terms), where $\nu$ is the Gaussian distribution with mean $0$ and  
variance
$\frac{1}{n}$.
\end{remark}

In the following, we consider a given family $a_0,\dots,a_r$ of
matrices in $\Mmsa$, and, for each $n$ in $\N$, a family
$X_1^{(n)},\dots,X_r^{(n)}$ of independent random matrices in
$\SGRM(n,\frac{1}{n})$. Furthermore, we consider the following random  
variable with values in $M_m(\C)\otimes M_n(\C)$:
\begin{equation}
S_n= a_0\otimes \unit_n + \sum_{i=1}^ra_i\otimes X_i^{(n)}.
\label{e2.0a}
\end{equation}

\begin{lemma}\label{lemma til master eq}
For each $n$ in $\N${\rm ,} let $S_n$ be as above. For any matrix
$\lambda$ in $M_m(\C)${\rm ,} for which $\im\lambda$ is positive definite{\rm ,}
we define a random variable with values in $M_m(\C)$ by \/{\rm (}\/cf.\
Lemma~{\rm \ref{imaginaerdel}),}
\[
H_n(\lambda)=(\id_m\otimes\tr_n)
\big[(\lambda\otimes\unit_n-S_n)^{-1}\big].
\]
Then{\rm ,} for any $j$ in $\{1,2,\dots,r\}${\rm ,} we have
\[
\E\big\{H_n(\lambda)a_jH_n(\lambda)\big\}=
\E\big\{ (\id_m\otimes\tr_n)\big[(\unit_m\otimes
X_j^{(n)})\cdot(\lambda\otimes\unit_n-S_n)^{-1}\big]\big\}.
\]
\end{lemma}

\Proof   Let $\lambda$ be a fixed matrix in $M_m(\C)$, such that
$\im\lambda$ is positive definite. Consider the canonical isomorphism
$\Psi\colon\CE_{r,n}\to\R^{rn^2}$, introduced in Remark~\ref{kanonisk
   iso}, and then define the mappings $\tilde{F}\colon\CE_{r,n}\to
M_m(\C)\otimes M_n(\C)$ and $F\colon\R^{rn^2}\to M_m(\C)\otimes
M_n(\C)$ by (cf.\ Lemma~\ref{imaginaerdel})
\begin{multline*}
\tilde{F}(v_1,\dots,v_r)\\
=\big(\lambda\otimes\unit_n - a_0\otimes\unit_n
   -\textstyle{\sum_{i=1}^ra_i\otimes
   v_i}\big)^{-1}, \qquad (v_1,\dots,v_r\in\Mnsa),
\end{multline*}
and
\[
F=\tilde{F}\circ\Psi^{-1}.
\]
Note then that
\[
\big(\lambda\otimes\unit_n-S_n\big)^{-1} =
F(\Psi(X_1^{(n)},\dots,X_r^{(n)})),
\]
where $\Y=\Psi(X_1^{(n)},\dots,X_r^{(n)})$ is a random variable taking
values in $\R^{rn^2}$, and the distribution of $\Y$ equals that of a
tuple $(\gamma_1,\dots,\gamma_{rn^2})$ of $rn^2$ independent
identically $N(0,\frac{1}{n})$-distributed real-valued random
variables.

Now, let $j$ in $\{1,2,\dots,r\}$ be fixed, and then define
\[
\begin{split}
X_{j,k,k}^{(n)}&=(X_j^{(n)})_{kk}, \qquad (1\le k\le n), \\[.2cm]
Y_{j,k,l}^{(n)}&=\sqrt{2}\re(X_j^{(n)})_{k,l}, \qquad (1\le k<l\le n),
\\[.2cm]
Z_{j,k,l}^{(n)}&=\sqrt{2}{\rm Im}(X_j^{(n)})_{k,l}, \qquad (1\le k<l\le n).
\end{split}
\]
Note that
$\big((X_{j,k,k}^{(n)})_{1\le k\le n},(Y_{j,k,l}^{(n)})_{1\le k<l\le n},
(Z_{j,k,l}^{(n)})_{1\le k<l\le n}\big)=\Psi_0(X_j^{(n)})$, where
$\Psi_0$ is the mapping defined in \eqref{e2.0} of
Remark~\ref{kanonisk iso}. Note also that the standard orthonormal
basis for $\R^{n^2}$ corresponds, via $\Psi_0$, to the following
orthonormal basis for $\Mnsa$:
\begin{eqnarray}
&&e_{k,k}^{(n)}, \qquad (1\le k\le n)\label{e2.1}
 \\ 
&&f_{k,l}^{(n)}=\textstyle{\frac{1}{\sqrt{2}}}\big(e_{k,l}^{(n)}+e_{l,k}^ 
{(n)}\big)
\qquad (1\le k<l\le n), \nonumber\\ 
&&g_{k,l}^{(n)}=\textstyle{\frac{\i}{\sqrt{2}}}\big(e_{k,l}^{(n)}- 
e_{l,k}^{(n)}\big)
\qquad (1\le k<l\le n).
\nonumber 
\end{eqnarray}
In other words, $\big((X_{j,k,k}^{(n)})_{1\le k\le  
n},(Y_{j,k,l}^{(n)})_{1\le k<l\le n},
(Z_{j,k,l}^{(n)})_{1\le k<l\le n}\big)$ are the coefficients of
$X_j^{(n)}$ with respect to  the orthonormal basis set out in \eqref{e2.1}.

Combining now the above observations with Lemma~\ref{Gauss trick},
it follows that
\[
\begin{split}
\frac{1}{n}\E\Big\{\diff\big(\lambda\otimes\unit_n-S_n-ta_j\otimes
e_{k,k}^{(n)}\big)^{-1}\Big\}
&=\E\big\{X_{j,k,k}^{(n)}\cdot\big(\lambda\otimes\unit_n-S_n\big)^{ 
-1}\big\},
\\[.2cm]
\frac{1}{n}\E\Big\{\diff\big(\lambda\otimes\unit_n-S_n-ta_j\otimes
f_{k,l}^{(n)}\big)^{-1}\Big\}
&=\E\big\{Y_{j,k,l}^{(n)}\cdot\big(\lambda\otimes\unit_n-S_n\big)^{ 
-1}\big\},
\\[.2cm]
\frac{1}{n}\E\Big\{\diff\big(\lambda\otimes\unit_n-S_n-ta_j\otimes
g_{k,l}^{(n)}\big)^{-1}\Big\}
&=\E\big\{Z_{j,k,l}^{(n)}\cdot\big(\lambda\otimes\unit_n-S_n\big)^{ 
-1}\big\},
\end{split}
\]
for all values of $k,l$ in $\{1,2,\dots,n\}$ such that $k<l$.
On the other hand, it follows from Lemma~\ref{Kaplansky trick} that
for any vector $v$ in $\Mnsa$,
\[
\diff\big(\lambda\otimes\unit_n-S_n-ta_j\otimes v\big)^{-1}
=(\lambda\otimes\unit_n-S_n)^{-1}(a_j\otimes
v)(\lambda\otimes\unit_n-S_n)^{-1},
\]
and we obtain thus the identities:
\begin{eqnarray}
&&\E\big\{X_{j,k,k}^{(n)}\cdot\big(\lambda\otimes\unit_n-S_n\big)^{ 
-1}\big\}
\label{e2.2}
\\
&&\qquad\qquad =\frac{1}{n}\E\big\{(\lambda\otimes\unit_n-S_n)^{-1}(a_j\otimes
e_{k,k}^{(n)})(\lambda\otimes\unit_n-S_n)^{-1}\big\}
\nonumber\\
&& \E\big\{Y_{j,k,l}^{(n)}\cdot\big(\lambda\otimes\unit_n-S_n\big)^{ 
-1}\big\}
\label{e2.3}\\
&&\qquad\qquad =\frac{1}{n}\E\big\{(\lambda\otimes\unit_n-S_n)^{-1}(a_j\otimes
f_{k,l}^{(n)})(\lambda\otimes\unit_n-S_n)^{-1}\big\}
\nonumber
\\ 
&&\E\big\{Z_{j,k,l}^{(n)}\cdot\big(\lambda\otimes\unit_n-S_n\big)^{ 
-1}\big\}
\label{e2.4}\\
&&\qquad\qquad =\frac{1}{n}\E\big\{(\lambda\otimes\unit_n-S_n)^{-1}(a_j\otimes
g_{k,l}^{(n)})(\lambda\otimes\unit_n-S_n)^{-1}\big\}
\nonumber
\end{eqnarray}
for all relevant values of $k,l$, $k<l$. Note next that for $k<l$, we
have
\[
\begin{split}
(X_j^{(n)})_{k,l}
&=\textstyle{\frac{1}{\sqrt{2}}}\big(Y_{j,k,l}^{(n)}+\i  
Z_{j,k,l}^{(n)}\big),
\\[.2cm]
(X_j^{(n)})_{l,k}
&=\textstyle{\frac{1}{\sqrt{2}}}\big(Y_{j,k,l}^{(n)}-\i  
Z_{j,k,l}^{(n)}\big),
\\[.2cm]
e_{k,l}^{(n)}
&=\textstyle{\frac{1}{\sqrt{2}}}\big(f_{k,l}^{(n)}-\i  
g_{k,l}^{(n)}\big),
\\[.2cm]
e_{l,k}^{(n)}
&=\textstyle{\frac{1}{\sqrt{2}}}\big(f_{k,l}^{(n)}+\i  
g_{k,l}^{(n)}\big),
\\[.2cm]
\end{split}
\]
and combining this with \eqref{e2.3}--\eqref{e2.4}, it follows that
\begin{multline}
\E\big\{(X_j^{(n)})_{k,l}\cdot\big(\lambda\otimes\unit_n-S_n\big)^{ 
-1}\big\}
\\
=\frac{1}{n}\E\big\{(\lambda\otimes\unit_n-S_n)^{-1}(a_j\otimes
e_{l,k}^{(n)})(\lambda\otimes\unit_n-S_n)^{-1}\big\},
\label{e2.5}
\end{multline}
and that
\begin{multline}
\E\big\{(X_j^{(n)})_{l,k}\cdot\big(\lambda\otimes\unit_n-S_n\big)^{ 
-1}\big\}
\\
=\frac{1}{n}\E\big\{(\lambda\otimes\unit_n-S_n)^{-1}(a_j\otimes
e_{k,l}^{(n)})(\lambda\otimes\unit_n-S_n)^{-1}\big\},
\label{e2.6}
\end{multline}
for all $k,l$, $k<l$. Taking also \eqref{e2.2} into account, it follows
that \eqref{e2.5} actually holds for all $k,l$ in
$\{1,2,\dots,n\}$. By adding the equation \eqref{e2.5} for
all values of $k,l$ and by recalling that
\[
X_j^{(n)}=\sum_{1\le
k,l\le n}(X_j^{(n)})_{k,l}e_{k,l}^{(n)},
\]
we conclude that
\begin{multline}
\E\big\{(\unit_m\otimes
X_j^{(n)}) (\lambda\otimes\unit_n-S_n)^{-1}\big\} \\ 
 =\frac{1}{n}\sum_{1\le k,l\le n}\E\big\{(\unit_m\otimes
e_{k,l}^{(n)})(\lambda\otimes\unit_n-S_n)^{-1}
(a_j\otimes e_{l,k}^{(n)})(\lambda\otimes\unit_n-S_n)^{-1}\big\}.
\label{e2.7}
\end{multline}
To calculate the right-hand side of \eqref{e2.7}, we write
\[
\big(\lambda\otimes\unit_n-S_n\big)^{-1}=\sum_{1\le u,v\le
   n}F_{u,v}\otimes e_{u,v},
\]
where, for all $u,v$ in $\{1,2,\dots,n\}$, $F_{u,v}\colon\Omega\to
M_m(\C)$ is an $M_m(\C)$-valued random variable. Recall then that
for any $k,l,u,v$ in $\{1,2,\dots,n\}$,
\[
e_{k,l}^{(n)}\cdot e_{u,v}^{(n)} =
\begin{cases}
e_{k,v}, &\textrm{if} \ l=u, \\
0, &\textrm{if} \ l\ne u.
\end{cases}
\]
For any fixed $u,v$ in $\{1,2,\dots,n\}$, it follows thus that
\begin{equation}
\sum_{1\le k,l\le n}(\unit_m\otimes e_{k,l}^{(n)})(F_{u,v}\otimes
e_{u,v}^{(n)})(a_j\otimes e_{l,k}^{(n)})=
\begin{cases}
(F_{u,u}\cdot a_j)\otimes\unit_n, &\textrm{if} \ u=v, \\
0, &\textrm{if} \ u\ne v.
\end{cases}
\label{e2.8}
\end{equation}

Adding the equation \eqref{e2.8} for all values of $u,v$ in
$\{1,2,\dots,n\}$, it follows that
\[
\sum_{1\le k,l\le n}(\unit_m\otimes
e_{k,l}^{(n)})(\lambda\otimes\unit_n-S_n)^{-1}
(a_j\otimes e_{l,k}^{(n)})=
\big(\textstyle{\sum_{u=1}^nF_{u,u}a_j}\big)\otimes\unit_n.
\]
Note here that
\[
\sum_{u=1}^nF_{u,u}=n\cdot\id_m\otimes\tr_n\big[(\lambda\otimes\unit_n- 
S_n)^{-1}\big]
=n\cdot H_n(\lambda),
\]
so that
\[
\sum_{1\le k,l\le n}(\unit_m\otimes
e_{k,l}^{(n)})(\lambda\otimes\unit_n-S_n)^{-1}
(a_j\otimes e_{l,k}^{(n)})=
nH_n(\lambda)a_j\otimes\unit_n.
\]

Combining this with \eqref{e2.7}, we find that
\begin{equation}
\E\big\{(\unit_m\otimes
X_j^{(n)})(\lambda\otimes\unit_n-S_n)^{-1}\big\}
=\E\big\{(H_n(\lambda)a_j\otimes\unit_n)
(\lambda\otimes\unit_n-S_n)^{-1}\big\}.
\label{e2.9}
\end{equation}
Applying finally $\id_m\otimes\tr_n$ to both sides of \eqref{e2.9},
we conclude that
\begin{eqnarray*}
&&\hskip-36pt\E\big\{\id_m\otimes\tr_n\big[(\unit_m\otimes
X_j^{(n)})(\lambda\otimes\unit_n-S_n)^{-1}\big]\big\}\\
&&\qquad =\E\big\{H_n(\lambda)a_j\cdot
\id_m\otimes\tr_n\big[(\lambda\otimes\unit_n-S_n)^{-1}\big]\big\}
\\ 
&&\qquad=\E\big\{H_n(\lambda)a_jH_n(\lambda)\big\},
\end{eqnarray*}
which is the desired formula.
\hfill\qed

\begin{theorem}[Master equation]\label{master eq}
Let{\rm ,} for each $n$ in $\N${\rm ,} $S_n$ be the random matrix introduced in
{\rm \eqref{e2.0a},} and let $\lambda$ be a matrix in $M_m(\C)$ such that
${\rm Im}(\lambda)$ is positive definite. Then with
\[
H_n(\lambda)=(\id_m\otimes\tr_n)
\big[(\lambda\otimes\unit_n-S_n)^{-1}\big]
\]
\/{\rm (}\/cf.\ Lemma~{\rm \ref{imaginaerdel}),} we have the formula
\begin{equation}
\E\Big\{\sum_{i=1}^r a_iH_n(\lambda)a_iH_n(\lambda) +  
(a_0-\lambda)H_n(\lambda)
+\unit_m\Big\}=0,
\label{e2.10}
\end{equation}
as an $M_m(\C)$\/{\rm -}\/valued expectation.
\end{theorem}

\Proof   By application of Lemma~\ref{lemma til master eq}, we find that
\begin{eqnarray*}
&&\hskip-36pt \E\Big\{\sum_{j=1}^ra_jH_n(\lambda)a_jH_n(\lambda)\Big\}
\\[.2cm]
&&\qquad=\sum_{j=1}^ra_j\E\big\{H_n(\lambda)a_jH_n(\lambda)\big\} \\[.2cm]
&&\qquad= \sum_{j=1}^ra_j\E\big\{\id_m\otimes\tr_n\big[(\unit_m\otimes
X_j^{(n)})(\lambda\otimes\unit_n-S_n)^{-1}\big]\big\} \\[.2cm]
&&\qquad=\sum_{j=1}^r\E\big\{\id_m\otimes\tr_n\big[(a_j\otimes\unit_n)(\unit_m\otimes
X_j^{(n)})(\lambda\otimes\unit_n-S_n)^{-1}\big]\big\} \\[.2cm]
 &&\qquad=\sum_{j=1}^r\E\big\{\id_m\otimes\tr_n\big[(a_j\otimes
X_j^{(n)})(\lambda\otimes\unit_n-S_n)^{-1}\big]\big\}.
\end{eqnarray*}
Moreover,
\begin{eqnarray*}
\E\{a_0 H_n(\lambda)\} &=& \E\{a_0(\id_m\otimes\tr_n)((\lambda\otimes
\unit_n-S_n)^{-1})\}\\
&=&
\E\{(\id_m\otimes\tr_n)((a_0\otimes\unit_n)(\lambda\otimes\unit_n- 
S_n)^{-1}\}.
\end{eqnarray*}
Hence,
\[
\begin{split}
\E\Big\{a_0H_n(\lambda)+\sum^r_{i=1}  
&a_jH_n(\lambda)a_jH_n(\lambda)\Big\} \\[.2cm]
&=\E\big\{\id_m\otimes\tr_n\big[S_n(\lambda\otimes\unit_n-S_n)^{ 
-1}\big]\big\}
\\[.2cm]
&=\E\big\{\id_m\otimes\tr_n\big[\big(\lambda\otimes\unit_n
-(\lambda\otimes\unit_n-S_n)\big)(\lambda\otimes\unit_n-S_n)^{ 
-1}\big]\big\}
\\[.2cm]
&=\E\big\{\id_m\otimes\tr_n\big[(\lambda\otimes\unit_n)(\lambda\otimes\unit_n-S_n)^{-1}
-\unit_m\otimes\unit_n\big]\big\}
\\[.2cm]
&=\E\big\{\lambda H_n(\lambda)-\unit_m\big\},
\end{split}
\]
from which \eqref{e2.10} follows readily.
\hfill\qed

\section{Variance estimates} \label{sec3}

Let $K$ be a positive integer. Then we denote by $\|\cdot\|$ the
usual Euclidean norm $\C^K$; i.e.,
\[
\|(\zeta_1,\dots,\zeta_K)\|
=\big(|\zeta_1|^2+\cdots+|\zeta_K|^2\big)^{1/2},
\qquad (\zeta_1,\dots,\zeta_K\in\C).
\]
Furthermore, we denote by $\|\cdot\|_{2,\Tr_K}$ the Hilbert-Schmidt  
norm on
$M_K(\C)$, i.e.,
\[
\|T\|_{2,\Tr_K}=\big(\Tr_K(T^*T)\big)^{1/2}, \qquad (T\in M_K(\C)).
\]
We shall also, occasionally, consider the norm $\|\cdot\|_{2,\tr_k}$  
given
by:
\[
\|T\|_{2,\tr_K}=\big(\tr_K(T^*T)\big)^{1/2}=K^{-1/2}\|T\|_{2,\Tr_K},  
\qquad
(T\in M_K(\C)).
\]

\begin{proposition}[Gaussian Poincar\'e inequality]\label{concentration}
Let $N$ be a positive integer and equip $\R^N$ with the probability  
measure
$\mu=\nu\otimes\nu\otimes\cdots\otimes\nu$ \/{\rm (}\/$N$~terms\/{\rm ),}\/ where $\nu$ is  
the
Gaussian distribution on $\R$ with mean $0$ and variance~$1$. Let
$f\colon\R^N\to\C$ be a $C^1$-function{\rm ,} such that $\E\{|f|^2\}<\infty$.
Then with $\V\{f\}=\E\{|f-\E\{f\}|^2\}${\rm ,} we have
\[
\V\{f\}\le \E\big\{\|\grad(f)\|^2\big\}.
\]
\end{proposition}

\Proof   See \cite[Thm.~2.1]{Cn}.
\Endproof\vskip4pt 

The Gaussian Poincar\'e inequality is a folklore result which goes back
to the 30's (cf.~Beckner \cite{Be}). It was rediscovered by Chernoff
\cite{Cf} in 1981 in the case $N=1$ and by Chen \cite{Cn} in 1982 for
general $N$. The original proof as well as Chernoff's proof is based on
an expansion of $f$ in Hermite polynomials (or tensor products of
Hermite polynomials in the case $N\ge 2$). Chen gives in \cite{Cn} a
self-contained proof which does not rely on Hermite polynomials. In a
preliminary version of this paper, we proved the slightly weaker
inequality: $\V\{f\}\le\frac{\pi^2}{8}\E\{\| \grad f \|^2\}$ using the  
method
of proof of \cite[Lemma 4.7]{P1}. We wish to thank Gilles Pisier for  
bringing
the papers by Bechner, Chernoff and Chen to our attention.

\begin{corollary} \label{cor3-2}
Let $N\in\N${\rm ,} and let $Z_1,\dots,Z_N$ be $N$ independent and identically distributed  real
Gaussian random variables with mean zero and variance $\sigma^2$ and let
$f\colon\R^N\to\C$ be a $C^1$\/{\rm -}\/function{\rm ,} such that $f$ and $\grad(f)$  
are both
polynomially bounded. Then
\[
\V\big\{ f(Z_1,\dots,Z_N)\big\}\le\sigma^2\E\big\{\|(\grad  
f)(Z_1,\dots,Z_N)\|^2\big\}.
\]
\end{corollary}

\Proof  
In the case $\sigma =1$, this is an immediate consequence of
Proposition~\ref{concentration}.
In the general case, put $Y_j=\frac{1}{\sigma} Z_j$, $j=1,\dots,N$, and
define $g\in C^1(\R^N)$ by
\begin{equation}
\label{eq31}
g(y) = f(\sigma y),\qquad (y\in\R^N).
\end{equation}
Then
\begin{equation}
\label{eq32}
(\grad g)(y) = \sigma(\grad f)(\sigma y),\qquad (y\in\R^N).
\end{equation}
Since $Y_1,\dots,Y_N$ are independent standard Gaussian distributed  
random
variables, we have from Proposition~\ref{concentration} that
\begin{equation}
\label{eq33}
\V\big\{g(Y_1,\dots,Y_N)\big\}\le\E\big\{\|(\grad  
g)(Y_1,\dots,Y_N)\|^2\big\}.
\end{equation}
Since $Z_j=\sigma Y_j$, $j=1,\dots,N$, it follows from \eqref{eq31},
\eqref{eq32}, and \eqref{eq33} that
\vglue12pt
\hfill $
\displaystyle{\V\big\{f(Z_1,\dots,Z_N)\big\}
\le\sigma^2\E\big\{\|(\grad f)(Z_1,\dots,Z_N)\|^2\big\}}. $\hfill\qed

\begin{remark}\label{matrix concentration}
Consider the canonical isomorphism $\Psi\colon\CE_{r,n}\to\R^{rn^2}$
introduced in Remark~\ref{kanonisk iso}.
Consider further independent random matrices  
$X_1^{(n)},\dots\break \dots,X_r^{(n)}$ from
$\SGRM(n,\frac{1}{n})$. Then $\X=(X_1^{(n)},\dots,X_r^{(n)})$ is a  
random
variable taking values in $\CE_{r,n}$, so that $\Y=\Psi(\X)$ is a random
variable taking values in $\R^{rn^2}$. As mentioned in
Remark~\ref{kanonisk iso}, it is easily seen that the distribution of
$\Y$ on $\R^{rn^2}$
is the product measure $\mu=\nu\otimes\nu\otimes\cdots\otimes\nu$  
($rn^2$
terms), where $\nu$ is the Gaussian distribution with mean $0$ and  
variance
$\frac{1}{n}$.
Now, let
$\tilde{f}\colon\R^{rn^2}\to\C$ be a $C^1$-function, such that
$\tilde{f}$ and $\grad\tilde{f}$ are both polynomially bounded, and  
consider
further the $C^1$-function $f\colon\CE_{r,n}\to\C$ given by
$f=\tilde{f}\circ\Psi$.
Since $\Psi$ is a linear isometry (i.e., an orthogonal transformation),  
it
follows from Corollary~\ref{cor3-2} that
\begin{equation}
\V\big\{f(\X)\big\}\le
\frac 1n \E\big\{\big\|\grad f(\X)\big\|_e^2\big\}.
\label{e3.0a}
\end{equation}
\end{remark}

\begin{lemma}\label{tensor estimat}
Let $m,n$ be positive integers{\rm ,} and assume that $a_1,\dots,a_r\in  
\Mmsa$
and $w_1,\dots,w_r\in M_n(\C)$. Then
\[
\Big\|\sum_{i=1}^ra_i\otimes w_i\Big\|_{2,\Tr_m\otimes\Tr_n}
\le m^{1/2}\Big\|\sum_{i=1}^ra_i^2\Big\|^{1/2}
\Big(\sum_{i=1}^r\|w_i\|_{2,\Tr_n}^2\Big)^{1/2}.
\]
\end{lemma}\vglue8pt

\Proof   We find that
\[
\begin{split}
\big\|\textstyle{\sum_{i=1}^ra_i\otimes w_i}\big\|_{2,\Tr_m\otimes\Tr_n}
&\le\textstyle{\sum_{i=1}^r\|a_i\otimes w_i\|_{2,\Tr_m\otimes\Tr_n}}  
\\[.2cm]
&=\textstyle{\sum_{i=1}^r\|a_i\|_{2,\Tr_m}\cdot\|w_i\|_{2,\Tr_n}}  
\\[.2cm]
&\le \big(\textstyle{\sum_{i=1}^r\|a_i\|_{2,\Tr_m}^2}\big)^{1/2}
\big(\textstyle{\sum_{i=1}^r\|w_i\|_{2,\Tr_n}^2}\big)^{1/2} \\[.2cm]
&=\big(\Tr_m\big(\textstyle{\sum_{i=1}^ra_i^2}\big)\big)^{1/2}\cdot
\big(\textstyle{\sum_{i=1}^r\|w_i\|_{2,\Tr_n}^2}\big)^{1/2} \\[.2cm]
&\le m^{1/2}\big\|\textstyle{\sum_{i=1}^ra_i^2}\big\|^{1/2}\cdot
\big(\textstyle{\sum_{i=1}^r\|w_i\|_{2,\Tr_n}^2}\big)^{1/2}.
\end{split}
\]
\vglue-20pt\Endproof 
\vskip12pt

Note, in particular, that if $w_1,\dots,w_r\in\Mnsa$, then
Lemma~\ref{tensor estimat} provides the estimate:
\[
\big\|\textstyle{\sum_{i=1}^ra_i\otimes w_i}\big\|_{2,\Tr_m\otimes\Tr_n}
\le m^{1/2}\big(\textstyle{\sum_{i=1}^r}\|a_i\|^2\big)^{1/2}\cdot
\big\|(w_1,\dots,w_r)\big\|_e.
\]

\begin{theorem}[Master inequality]\label{master ineq}
Let $\lambda$ be a matrix in $M_m(\C)$ such that ${\rm Im}(\lambda)$ is  
positive
definite. Consider further the random matrix $H_n(\lambda)$ introduced  
in
Theorem~{\rm \ref{master eq}} and put
\[
G_n(\lambda)=\E\big\{H_n(\lambda)\big\}\in M_m(\C).
\]
Then
\[
\Big\|\sum_{i=1}^ra_iG_n(\lambda)a_iG_n(\lambda) +(a_0-\lambda)  
G_n(\lambda) +
\unit_m\Big\| \le \frac{C}{n^2}\big\|({\rm Im}(\lambda))^{-1}\big\|^4,
\]
where $C= m^3\|\sum_{i=1}^ra_i^2\|^2$.
\end{theorem}

\Proof   We put
\[
K_n(\lambda)=H_n(\lambda)-G_n(\lambda)
=H_n(\lambda)-\E\big\{H_n(\lambda)\big\}.
\]
Then, by Theorem~\ref{master eq}, we have
\[
\begin{split}
\E\Big\{\sum_{i=1}^ra_i&K_n(\lambda)a_iK_n(\lambda)\Big\} \\
&=\E\Big\{\sum_{i=1}^ra_i\big(H_n(\lambda)-G_n(\lambda)\big)
a_i\big(H_n(\lambda)-G_n(\lambda)\big)\Big\} \\
&=\E\Big\{\sum_{i=1}^ra_iH_n(\lambda)a_iH_n(\lambda)\Big\}
-\sum_{i=1}^ra_iG_n(\lambda)a_iG_n(\lambda) \\ &=\Big( -(a_0-\lambda)
\E\big\{H_n(\lambda)\big\}-\unit_m\Big)-
\sum_{i=1}^ra_iG_n(\lambda)a_iG_n(\lambda) \\ 
&=-\Big(\sum_{i=1}^ra_iG_n(\lambda)a_iG_n(\lambda)
+ (a_0- \lambda)   G_n(\lambda) +
\unit_m\Big).
\end{split}
\]
Hence, we can make the following estimates
\begin{eqnarray*}
&&\hskip-.75in\Big\|\sum_{i=1}^ra_iG_n(\lambda)a_iG_n(\lambda) +  
(a_0-\lambda)G_n(\lambda) +
\unit_m\Big\|
\\ &&= \Big\|\E\Big\{\sum_{i=1}^ra_iK_n(\lambda)a_iK_n(\lambda)\Big\}\Big\|
\\
&&\le  
\E\Big\{\Big\|\sum_{i=1}^ra_iK_n(\lambda)a_iK_n(\lambda)\Big\|\Big\}
\\
&&\le \E\Big\{\Big\|\sum_{i=1}^r a_iK_n(\lambda)a_i\Big\|\cdot
\big\|K_n(\lambda)\big\|\Big\}.
\end{eqnarray*}

Note here that since $a_1,\dots,a_r$ are self-adjoint, the mapping
$v\mapsto\break \sum_{i=1}^ra_iva_i: M_m(\C)\to M_m(\C)$ is
completely positive. Therefore, it attains its norm at the unit
$\unit_m$, and the norm is $\|\sum_{i=1}^ra_i^2\|$. Using
this in the estimates above, we find that
\begin{equation}
\begin{split}
\Big\|\sum_{i=1}^ra_iG_n(\lambda)a_iG_n(\lambda)+(a_0-\lambda)  
G_n(\lambda) +
\unit_m\Big\|
&\le\Big\|\sum_{i=1}^ra_i^2\Big\|\cdot
\E\Big\{\big\|K_n(\lambda)\big\|^2\Big\} \\[.2cm]
&\le\Big\|\sum_{i=1}^ra_i^2\Big\|\cdot
\E\Big\{\big\|K_n(\lambda)\big\|_{2,\Tr_m}^2\Big\},
\end{split}
\label{e3.1}
\end{equation}
where the last inequality uses that the operator norm of a matrix is
always dominated by the Hilbert-Schmidt norm.
It remains to estimate $\E\{\|K_n(\lambda)\|_{2,\Tr_m}^2\}$.
For this, let $H_{n,j,k}(\lambda)$, ($1\le j,k\le n$) denote the  
entries of
$H_n(\lambda)$; i.e.,
\begin{equation}
H_n(\lambda)=\sum_{j,k=1}^mH_{n,j,k}(\lambda)e(m,j,k),
\label{e3.2}
\end{equation}
where $e(m,j,k)$, ($1\le j,k\le m$) are the usual $m\times m$ matrix
units. Let, correspondingly, $K_{n,j,k}(\lambda)$ denote the entries of
$K_n(\lambda)$. Then
$K_{n,j,k}(\lambda)=H_{n,j,k}(\lambda)-\E\{H_{n,j,k}(\lambda)\}$, for  
all
$j,k$, so that $\V\{H_{n,j,k}(\lambda)\}=\E\{|K_{n,j,k}(\lambda)|^2\}$.  
Thus it
follows   that
\begin{equation}
\E\Big\{\big\|K_n(\lambda)\big\|_{2,\Tr_m}^2\Big\}
=\E\Big\{\sum_{j,k=1}^m|K_{n,j,k}(\lambda)|^2\Big\}
=\sum_{j,k=1}^m\V\big\{H_{n,j,k}(\lambda)\big\}.
\label{e3.2d}
\end{equation}

Note further that by \eqref{e3.2}
\[
\begin{split}
H_{n,j,k}(\lambda)&=\Tr_m\big(e(m,k,j)H_n(\lambda)\big) \\[.2cm]
&=m\cdot\tr_m\big(e(m,k,j)\cdot(\id_m\otimes\tr_n)
\big[(\lambda\otimes\unit_n-S_n)^{-1}\big]\big) \\[.2cm]
&=m\cdot\tr_m\otimes\tr_n\big[(e(m,j,k)\otimes\unit_n)
(\lambda\otimes\unit_n-S_n)^{-1}\big].
\end{split}
\]
For any $j,k$ in $\{1,2,\dots,m\}$, consider next the mapping
$f_{n,j,k}\colon\CE_{r,n}\to\C$ given by:
\begin{multline*}
f_{n,j,k}(v_1,\dots,v_r)\\
=m\cdot(\tr_m\otimes\tr_n)
\big[(e(m,k,j)\otimes\unit_n)(\lambda\otimes\unit_n
-a_0\otimes\unit_n -\textstyle{\sum_{i=1}^ra_i\otimes v_i})^{-1}\big],
\end{multline*}
for all $v_1,\dots,v_r$ in $\Mnsa$. Note then that
\[
H_{n,j,k}(\lambda)=f_{n,j,k}(X_1^{(n)},\dots,X_r^{(n)}),
\]
for all $j,k$. Using now the ``concentration estimate'' \eqref{e3.0a} in
Remark~\ref{matrix concentration}, it follows that for all $j,k$,
\begin{equation}
\V\big\{H_{n,j,k}(\lambda)\big\}\le \frac 1n \E\Big\{\big\|\grad
f_{n,j,k}(X_1^{(n)},\dots,X_r^{(n)})\big\|_e^2\Big\}.
\label{e3.2c}
\end{equation}
For fixed $j,k$ in $\{1,2,\dots,m\}$ and $v=(v_1,\dots,v_r)$ in
$\CE_{r,n}$, note that $\grad f_{n,j,k}(v)$ is the
vector in $\CE_{r,n}$, characterized by the property that
\[
\big\langle\grad f_{n,j,k}(v),w\big\rangle_e = \diff f_{n,j,k}(v+tw),
\]
for any vector $w=(w_1,\dots,w_r)$ in $\CE_{r,n}$.
It follows thus that
\begin{eqnarray}
\label{e3.2b}
\big\|\grad f_{n,j,k}(v)\big\|_e^2& =&
\max_{w\in S_1(\CE_{r,n})}\big|\big\langle\grad
f_{n,j,k}(v),w\big\rangle_e\big|^2
\\
&=&\max_{w\in S_1(\CE_{r,n})}\Big|\diff f_{n,j,k}(v+tw)\Big|^2.
\nonumber
\end{eqnarray}

Let $v=(v_1,\dots,v_n)$ be a fixed vector in $\CE_{r,n}$, and put
$\Sigma=a_0\otimes\unit_n + \sum_{i=1}^ra_i\otimes v_i$. Let further
$w=(w_1,\dots,w_n)$ be a fixed vector in $S_1(\CE_{r,n})$. It follows
then by Lemma~\ref{Kaplansky trick} that
\begin{eqnarray} &&\\
\label{e3.2a}
&&\hskip-12pt \diff f_{n,j,k}(v+tw)\nonumber \\[.1cm] 
&&= \diff\! m\cdot(\tr_m\otimes\tr_n)\nonumber\\[.1cm]
 &&\quad\cdot \big[(e(m,k,j)\otimes\unit_n)\big(\lambda\otimes\unit_n -  
a_0\otimes\unit_n
-\textstyle{\sum_{i=1}^ra_i\otimes(v_i+tw_i)}\big)^{-1}\big]  \nonumber\\[.1cm]
&&=m\cdot(\tr_m\otimes\tr_n)\nonumber\\[.1cm]
 &&\quad\cdot \Big[(e(m,k,j)\otimes\unit_n)\,\diff\!\!\big(\lambda\otimes\unit_n -  
a_0\otimes\unit_n
-\textstyle{\sum_{i=1}^ra_i\otimes(v_i+tw_i)}\big)^{-1}\Big]  \nonumber\\[.1cm]
&&=m\cdot(\tr_m\otimes\tr_n)\nonumber\\[.1cm]
 &&\quad\cdot \big[(e(m,k,j)\otimes\unit_n)\big(\lambda\otimes\unit_n-\Sigma\big)^{-1}
\big(\textstyle{\sum_{i=1}^ra_i\otimes w_i}\big)
\big(\lambda\otimes\unit_n-\Sigma\big)^{-1}\big].\nonumber
\end{eqnarray}

Using next the Cauchy-Schwartz inequality for $\Tr_n\otimes\Tr_m$, we  
find that
\begin{eqnarray}
&&m^2\big|(\tr_m\otimes\tr_n)
\big[e(m,k,j)\otimes\unit_n\label{e3.3}
\\[.1cm]
&&\qquad \cdot\big(\lambda\otimes\unit_n- 
\Sigma\big)^{-1}
\big(\textstyle{\sum_{i=1}^ra_i\otimes w_i}\big)
\big(\lambda\otimes\unit_n-\Sigma\big)^{-1}\big]\big|^2\nonumber \\[.1cm]
&&=\frac{1}{n^2}\big|(\Tr_m\otimes\Tr_n)
\big[e(m,k,j)\otimes\unit_n\nonumber\\[.1cm]
&&\qquad \cdot\big(\lambda\otimes\unit_n- 
\Sigma\big)^{-1}
\big(\textstyle{\sum_{i=1}^ra_i\otimes w_i}\big)
\big(\lambda\otimes\unit_n-\Sigma\big)^{-1}\big]\big|^2 \nonumber\\[.1cm]
&&\le\frac{1}{n^2}\big\|e(m,j,k)
\otimes\unit_n\big\|_{2,\Tr_m\otimes\Tr_n}^2\nonumber\\[.1cm]
&&\qquad\cdot
\big\|\big(\lambda\otimes\unit_n-\Sigma\big)^{-1}
\big(\textstyle{\sum_{i=1}^ra_i\otimes w_i}\big)
\big(\lambda\otimes\unit_n-\Sigma\big)^{ 
-1}\big\|_{2,\Tr_m\otimes\Tr_n}^2 \nonumber\\[.1cm]
&&=\frac{1}{n}\big\|\big(\lambda\otimes\unit_n-\Sigma\big)^{-1}
\big(\textstyle{\sum_{i=1}^ra_i\otimes w_i}\big)
\big(\lambda\otimes\unit_n-\Sigma\big)^{ 
-1}\big\|_{2,\Tr_m\otimes\Tr_n}^2.
\nonumber
\end{eqnarray}
Note here that
\begin{eqnarray*}
&&\hskip-12pt\big\|\big(\lambda\otimes\unit_n-\Sigma\big)^{-1}
\big(\textstyle{\sum_{i=1}^ra_i\otimes w_i}\big)
\big(\lambda\otimes\unit_n-\Sigma\big)^{ 
-1}\big\|_{2,\Tr_m\otimes\Tr_n}^2 \\[.2cm]
&&\qquad\le \big\|\big(\lambda\otimes\unit_n-\Sigma\big)^{-1}\big\|^2\cdot
\big\|\textstyle{\sum_{i=1}^ra_i\otimes  
w_i}\big\|_{2,\Tr_m\otimes\Tr_n}^2\cdot
\big\|\big(\lambda\otimes\unit_n-\Sigma\big)^{-1}\big\|^2 \\[.2cm]
&&\qquad\le \big\|\textstyle{\sum_{i=1}^ra_i\otimes  
w_i}\big\|_{2,\Tr_m\otimes\Tr_n}^2
\cdot\big\|\big({\rm Im}(\lambda)\big)^{-1}\big\|^4,
\end{eqnarray*}
where the last inequality uses Lemma~\ref{imaginaerdel} and the fact  
that
$\Sigma$ is self-adjoint:
\begin{eqnarray*}
\big\|\big(\lambda\otimes\unit_n-\Sigma\big)^{-1}\big\|
&\le& \big\|\big({\rm Im}(\lambda\otimes\unit_n-\Sigma\big)^{-1}\big\|
\\
&=&\big\|\big({\rm Im}(\lambda\otimes\unit_n)\big)^{-1}\big\|
=\big\|\big({\rm Im}(\lambda)\big)^{-1}\big\|.
\end{eqnarray*}
Note further that by Lemma~\ref{tensor estimat},
$\|\sum_{i=1}^ra_i\otimes w_i\|_{2,\Tr_m\otimes\Tr_n}
\!\le \! m^{1/2}\big\|\textstyle{\sum_{i=1}^ra_i^2}\big\|^{1/2}$,
since $w=(w_1,\dots,w_r)\in S_1(\CE_{r,n})$. We
conclude thus that
\begin{multline}
\big\|\big(\lambda\otimes\unit_n-\Sigma\big)^{-1}
\big(\textstyle{\sum_{i=1}^ra_i\otimes w_i}\big)
\big(\lambda\otimes\unit_n-\Sigma\big)^{ 
-1}\big\|_{2,\Tr_m\otimes\Tr_n}^2
\\
\le m\big\|\textstyle{\sum_{i=1}^ra_i^2}\big\|
\cdot\big\|\big({\rm Im}(\lambda)\big)^{-1}\big\|^4.
\label{e3.5}
\end{multline}

Combining now formulas
\eqref{e3.2a}--\eqref{e3.5}, it follows that for any $j,k$ in
$\{1,2,\dots,m\}$, any vector
$v=(v_1,\dots,v_r)$ in $\CE_{r,n}$ and any unit vector
$w=(w_1,\dots,w_r)$ in $\CE_{r,n}$, we have that
\[
\Big|\diff f_{n,j,k}(v+tw)\Big|^2
\le\frac{m}{n}\big\|\textstyle{\sum_{i=1}^ra_i^2}\big\|
\cdot\big\|\big({\rm Im}(\lambda)\big)^{-1}\big\|^4;
\]
  hence, by \eqref{e3.2b},
\[
\big\|\grad f_{n,j,k}(v)\big\|_e^2 \le
\frac{m}{n}\big\|\textstyle{\sum_{i=1}^ra_i^2}\big\|
\cdot\big\|\big({\rm Im}(\lambda)\big)^{-1}\big\|^4.
\]
Note that this estimate holds at any point $v=(v_1,\dots,v_r)$ in
$\CE_{r,n}$. Using this in conjunction with \eqref{e3.2c}, we may thus
conclude that
\[
\V\big\{H_{n,j,k}(\lambda)\big\}\le
\frac{m}{n^2}\big\|\textstyle{\sum_{i=1}^ra_i^2}\big\|
\cdot\big\|\big({\rm Im}(\lambda)\big)^{-1}\big\|^4,
\]
for any $j,k$ in $\{1,2\ldots,m\}$, and hence, by \eqref{e3.2d},
\begin{equation}
\E\Big\{\big\|K_n(\lambda)\big\|_{2,\Tr_m}^2\Big\}\le
\frac{m^3}{n^2}\big\|\textstyle{\sum_{i=1}^ra_i^2}\big\|
\cdot\big\|\big({\rm Im}(\lambda)\big)^{-1}\big\|^4.
\label{e3.6}
\end{equation}
Inserting finally \eqref{e3.6} into \eqref{e3.1}, we find that
\begin{multline*}
\Big\|\sum_{i=1}^ra_iG_n(\lambda)a_iG_n(\lambda)+(a_0-\lambda)  
G_n(\lambda) +
\unit_m\Big\|\\
\le \frac{m^3}{n^2}\big\|\textstyle{\sum_{i=1}^ra_i^2}\big\|^2
\cdot\big\|\big({\rm Im}(\lambda)\big)^{-1}\big\|^4,
\end{multline*}
and this is the desired estimate
\hfill\qed

\begin{lemma}\label{kaederegel}
Let $N$ be a positive integer{\rm ,} let $I$ be an open interval in $\R${\rm ,} and  
let
$t\mapsto a(t)\colon I\to M_N(\C)_{\rm sa}$ be a $C^1$\/{\rm -}\/function.  
Consider
further a function $\varphi$ in $C^1(\R)$. Then the function
$t\mapsto\tr_N[\varphi(a(t))]$ is $C^1$\/{\rm -}\/function on $I${\rm ,} and
\[
\frac{\rm d}{{\rm d}t}\tr_N\big[\varphi(a(t))\big]
=\tr_N\big[\varphi'(a(t))\cdot a'(t)\big].
\]
\end{lemma}
\vskip8pt
 
\Proof  
This is well known. For the reader's convenience we include a proof:
Note first that for any $k$ in $\N$,
\[
\frac{\rm d}{{\rm d}t}\big(a(t)^k\big) = \sum^{k-1}_{j=0} a(t)^j
a'(t)a(t)^{k-j-1}.
\]
Hence, by the trace property $\tr_N(xy)=\tr_N(yx)$, we get
\[
\frac{\rm d}{{\rm d}t} (\tr_N(a(t)^k) = \tr_N(k a(t)^{k-1}a'(t)).
\]
Therefore
\[
\frac{\rm d}{{\rm d}t} \tr_N(p(a(t))) = \tr_N(p'(a(t))a'(t))
\]
for all polynomials $p\in\C[X]$. The general case $\varphi\in
C^1(I)$ follows easily from this by choosing a sequence of
polynomials $p_n\in\C[X]$, such that $p_n\to\varphi$ and
$p'_n\to\varphi'$ uniformly on compact subsets of $I$, as $n\to\infty$.
\hfill\qed

\begin{proposition}\label{diff-kvotient trick}
Let $a_0,a_1,\dots,a_r$ be matrices in $\Mmsa$ and put as in  
\eqref{e2.0}
$$
S_n=a_0\otimes\unit_n+\sum_{i=1}^ra_i\otimes X_i^{(n)}.
$$
Let further $\varphi\colon\R\to\C$ be a $C^1$-function with compact
support{\rm ,} and consider the
random matrices $\varphi(S_n)$ and $\varphi'(S_n)$ obtained by
applying the spectral mapping associated to the self\/{\rm -}\/adjoint \/{\rm (}\/random\/{\rm )}\/  
matrix
$S_n$. We then have\/{\rm :}\/
$$
\V\big\{(\tr_m\otimes\tr_n)[\varphi(S_n)\big]\big\} \le
\frac{1}{n^2}\Big\|\sum_{i=1}^ra_i^2\Big\|^2
\E\big\{(\tr_m\otimes\tr_n)\big[|\varphi'|^2(S_n)\big]\big\}.
$$
\end{proposition}
\vskip8pt

\Proof   Consider the mappings $g\colon\CE_{r,n}\to M_{nm}(\C)_{\rm sa}$  
and
$f\colon\CE_{r,n}\to\C$ given by
\[
g(v_1,\dots,v_r)=a_0\otimes\unit_n+\sum_{i=1}^ra_i\otimes v_i, \qquad
(v_1,\dots,v_r\in\Mnsa),
\]
and
\[
f(v_1,\dots,v_r)=(\tr_m\otimes\tr_n)\big[\varphi(g(v_1,\dots,v_r))\big 
],
\qquad (v_1,\dots,v_r\in\Mmsa),
\]
Note then that
$S_n=g(X_1^{(n)},\dots,X_r^{(n)})$ and that
$$(\tr_m\otimes\tr_n)[\varphi(S_n)]=f(X_1^{(n)},\dots,X_r^{(n)}).
$$
Note also that $f$ is a bounded function on $\Mnsa$, and, by
Lemma~\ref{kaederegel},
it has bounded continuous partial derivatives. Hence, we obtain from
\eqref{e3.0a} in Remark~\ref{matrix concentration} that
\begin{equation}
\V\big\{(\tr_m\otimes\tr_n)[\varphi(S_n)]\big\}\le
\frac{1}{n}\E\Big\{\big\|\grad
f(X_1^{(n)},\dots,X_r^{(n)})\big\|_e^2\Big\}.
\label{e3.7}
\end{equation}

Recall next that for any $v$ in $\CE_{r,n}$, $\grad f(v)$ is the
vector in $\CE_{r,n}$, characterized by the property that
\[
\big\langle\grad f(v),w\big\rangle_e = \diff f(v+tw),
\]
for any vector $w=(w_1,\dots,w_r)$ in $\CE_{r,n}$.
It follows thus that
\begin{equation}
\big\|\grad f(v)\big\|_e^2 =
\max_{w\in S_1(\CE_{r,n})}\big|\big\langle\grad
f(v),w\big\rangle_e\big|^2
=\max_{w\in S_1(\CE_{r,n})}\Big|\diff f(v+tw)\Big|^2,
\label{e3.8}
\end{equation}
at any point $v=(v_1,\dots,v_r)$ of $\CE_{r,n}$. Now, let
$v=(v_1,\dots,v_r)$ be a fixed point in $\CE_{r,n}$ and let
$w=(w_1,\dots,w_r)$ be a fixed point in $S_1(\CE_{r,n})$.
By Lemma~\ref{kaederegel}, we have then that
\[
\begin{split}
\diff f(v+tw) &= \diff(\tr_m\otimes\tr_n)\big[\varphi(g(v+tw))\big]
\\[.2cm]
&=(\tr_m\otimes\tr_n)\Big[\varphi'(g(v))\cdot\diff g(v+tw)\Big] \\[.2cm]
&=(\tr_m\otimes\tr_n)\big[\varphi'(g(v))\cdot
\textstyle{\sum_{i=1}^ra_i\otimes w_i}\big].
\end{split}
\]
Using then the Cauchy-Schwartz inequality for $\Tr_m\otimes\Tr_n$, we  
find
that
\[
\begin{split}
\Big|\diff f(v+tw)\Big|^2 &=
\frac{1}{m^2n^2}\Big|(\Tr_m\otimes\Tr_n)\big[\varphi'(g(v))\cdot
\textstyle{\sum_{i=1}^ra_i\otimes w_i}\big]\Big|^2 \\[.2cm]
&=\frac{1}{n^2m^2}\big\|\overline{\varphi}'(g(v))\big\|_{2,\Tr_m\otimes\ 
Tr_n}^2\cdot
\big\|\textstyle{\sum_{i=1}^ra_i\otimes  
w_i}\big\|_{2,\Tr_m\otimes\Tr_n}^2.
\end{split}
\]
Note   that
\[
\big\|\overline{\varphi}'(g(v))\big\|_{2,\Tr_m\otimes\Tr_n}^2
=\Tr_m\otimes\Tr_n\big[|\varphi'|^2(g(v))\big]
=mn\cdot\tr_m\otimes\tr_n\big[|\varphi'|^2(g(v))\big],
\]
and, by Lemma~\ref{tensor estimat},
\[
\big\|\textstyle{\sum_{i=1}^ra_i\otimes  
w_i}\big\|_{2,\Tr_m\otimes\Tr_n}^2
\le m\big\|\textstyle{\sum_{i=1}^r}a_i^2\big\|,
\]
since $w$ is a unit vector with respect to  $\|\cdot\|_e$. We find thus that
\[
\Big|\diff f(v+tw)\Big|^2 \le
\frac{1}{n}\big\|\textstyle{\sum_{i=1}^r}a_i^2\big\|
\tr_m\otimes\tr_n\big[|\varphi'|^2(g(v))\big].
\]
Since this estimate holds for any unit vector $w$ in $\CE_{r,n}$, we
conclude, using \eqref{e3.8}, that
\[
\big\|\grad f(v)\big\|^2_e\le
\frac{1}{n}\big\|\textstyle{\sum_{i=1}^r}a_i^2\big\|
\tr_m\otimes\tr_n\big[|\varphi'|^2(g(v))\big],
\]
for any point $v$ in $\CE_{r,n}$. Combining this with \eqref{e3.7}, we
obtain the desired estimate.
\hfill\qed

\section{Estimation of $\|G_n(\lambda)-G(\lambda)\|$} \label{sec4}

\begin{lemma}\label{middelvaerdi af norm}
For each $n$ in $\N${\rm ,} let $X_n$ be a random matrix in  
$\SGRM(n,\frac{1}{n})$. Then
\begin{equation}
\E\big\{\|X_n\|\big\}\le 2+2\sqrt{\frac{\log(2n)}{2n}}, \qquad (n\in\N).
\label{e4.0}
\end{equation}
In particular{\rm ,} it follows that
\begin{equation}
\E\big\{\|X_n\|\big\}\le 4,
\label{e4.0a}
\end{equation}
for all $n$ in $\N$.
\end{lemma}

\Proof   In \cite[Proof of Lemma~3.3]{HT1} it was proved that for any $n$  
in
$\N$ and any positive number $t$, we have
\begin{equation}
\E\big\{\Tr_n(\exp(tX_n))\big\}\le
n\exp\big(2t+\textstyle{\frac{t^2}{2n}}\big).
\label{e4.1}
\end{equation}
Let $\lambda_{\max}(X_n)$ and $\lambda_{\min}(X_n)$ denote the largest
and smallest eigenvalue of $X_n$ as functions of $\omega\in\Omega$. Then
\[
\begin{split}
\exp(t\|X_n\|)&=\max\{\exp(t\lmax(X_n)),\exp(-t\lmin(X_n))\} \\[.2cm]
&\le\exp(t\lmax(X_n))+\exp(-t\lmin(X_n))\\[.2cm]
&\le
\Tr_n\big(\exp(tX_n)+\exp(-tX_n)\big).
\end{split}
\]
Using this in connection with Jensen's inequality, we find that
\begin{eqnarray}
\exp\big(t\E\{\|X_n\|\}\big)&\le& \E\big\{\exp(t\|X_n\|)\big\}
\label{e4.2}\\[.2cm]
&\le&\E\big\{\Tr_n(\exp(tX_n))\big\} + \E\big\{\Tr_n(\exp(-tX_n))\big\}
\nonumber\\[.2cm]
&=&2\E\big\{\Tr_n(\exp(tX_n))\big\},
\nonumber
\end{eqnarray}
where the last equality is due to the fact that
  $-X_n\in\SGRM(n,\frac{1}{n})$ too. Combining \eqref{e4.1} and
  \eqref{e4.2} we obtain the estimate
\[
\exp\big(t\E\{\|X_n\|\}\big)\le
2n\exp\big(2t+\textstyle{\frac{t^2}{2n}}\big),
\]
and hence, after taking logarithms and dividing by $t$,
\begin{equation}
\E\{\|X_n\|\}\le \frac{\log(2n)}{t}+2+\frac{t}{2n}.
\label{e4.3}
\end{equation}
This estimate holds for all positive numbers $t$. As a function of
$t$, the right-hand side of \eqref{e4.3} attains its minimal value at
$t_0=\sqrt{2n\log(2n)}$ and the minimal value is
$2+2\sqrt{\log(2n)/2n}$. Combining this with \eqref{e4.3} we
obtain \eqref{e4.0}. The estimate \eqref{e4.0a} follows subsequently
by noting that the function $t\mapsto\log(t)/t$ ($t>0$) attains
its maximal value at $t=\e$, and thus $2+2\sqrt{\log(t)/t}\le
2+2\sqrt{1/\e}\approx 3.21$ for all positive numbers $t$.
\Endproof\vskip4pt 

In the following we consider a fixed  positive integer $m$ and fixed
self-adjoint matrices $a_0,\dots,a_r$ in $\Mmsa$. We consider
further, for each positive integer $n$, independent random matrices
$X_1^{(n)},\dots,X_r^{(n)}$ in $\SGRM(n,\frac{1}{n})$. As in
Sections~\ref{sec2} and \ref{sec3}, we define
\[
S_n=a_0 + \sum_{i=1}^ra_i\otimes X_i^{(n)}.
\]
and, for any matrix $\lambda$ in $M_m(\C)$ such that
${\rm Im}(\lambda)$ is positive definite, we put
\[
H_n(\lambda)=(\id_m\otimes\tr_n)
\big[(\lambda\otimes\unit_n-S_n)^{-1}\big],
\]
and
\[
G_n(\lambda)=\E\{H_n(\lambda)\}.
\]

\begin{proposition}\label{estimat1}
Let $\lambda$ be a matrix in $M_m(\C)$ such that ${\rm Im}(\lambda)$ is
positive definite. Then $G_n(\lambda)$ is invertible and
\[
\big\|G_n(\lambda)^{-1}\big\|\le
\big(\|\lambda\|+K\big)^2\big\|(\im\lambda)^{-1}\big\|,
\]
where $K= \|a_0\|+ 4\sum_{i=1}^r\|a_i\|$.
\end{proposition}

\Proof   We note first that
\[
\begin{split}
\im&\big((\lambda\otimes\unit_n-S_n)^{-1}\big)\\[.2cm]
&=\frac{1}{2\i}\big((\lambda\otimes\unit_n-S_n)^{-1}
-(\lambda^*\otimes\unit_n-S_n)^{-1}\big) \\[.2cm]
&=\frac{1}{2\i}\big((\lambda\otimes\unit_n-S_n)^{-1}
\big((\lambda^*\otimes\unit_n-S_n)-(\lambda\otimes\unit_n-S_n)\big)
(\lambda^*\otimes\unit_n-S_n)^{-1}\big) \\[.2cm]
&=-(\lambda\otimes\unit_n-S_n)^{-1}({\rm Im}(\lambda)\otimes\unit_n)
(\lambda^*\otimes\unit_n-S_n)^{-1}.
\end{split}
\]
 From this it follows that $-{\rm Im}((\lambda\otimes\unit_n-S_n)^{-1})$ is
positive definite at any $\omega$ in~$\Omega$, and the inverse is
given by
\[
\big(-{\rm Im}((\lambda\otimes\unit_n-S_n)^{-1})\big)^{-1} =
(\lambda^*\otimes\unit_n-S_n)((\im\lambda)^{ 
-1}\otimes\unit_n)(\lambda\otimes\unit_n-S_n).
\]
In particular, it follows that
\[
0\le \big(-{\rm Im}((\lambda\otimes\unit_n-S_n)^{-1})\big)^{-1}
\le\big\|\lambda\otimes\unit_n-S_n\big\|^2\big\|(\im\lambda)^{-1}\big\|
\cdot\unit_m\otimes\unit_n,
\]
and this implies that
\[
-{\rm Im}\big((\lambda\otimes\unit_n-S_n)^{ 
-1}\big)\ge\frac{1}{\|\lambda\otimes\unit_n-S_n\|^2\|(\im\lambda)^{ 
-1}\|}
\cdot\unit_m\otimes\unit_n.
\]
Since the slice map $\id_m\otimes\tr_n$ is positive, we have thus
established that
\begin{eqnarray*}
-\im H_n(\lambda)&\ge&
\frac{1}{\|\lambda\otimes\unit_n-S_n\|^2\|(\im\lambda)^{-1}\|}
\cdot\unit_m
\\
&\ge& \frac{1}{(\|\lambda\|+\|S_n\|)^2\|(\im\lambda)^{-1}\|}
\cdot\unit_m,
\end{eqnarray*}
so that
\[
-\im G_n(\lambda) = \E\{-\im  
H_n(\lambda)\}\ge\frac{1}{\|(\im\lambda)^{-1}\|}
\E\Big\{\frac{1}{(\|\lambda\|+\|S_n\|)^2}\Big\}\unit_m.
\]

Note here that the function $t\mapsto\frac{1}{(\|\lambda\|+t)^2}$ is  
convex,
so applying Jensen's inequality to the random variable $\|S_n\|$,
yields the estimate
\[
\E\Big\{\frac{1}{(\|\lambda\|+\|S_n\|)^2}\Big\}\ge
\frac{1}{(\|\lambda\|+\E\{\|S_n\|\})^2},
\]
where
\begin{eqnarray*}
\E\{\|S_n\|\}&\le&\E\Big\{ \|a_0\| +  
\sum_{i=1}^r\|a_i\|\cdot\|X_i^{(n)}\|\Big\}\\
& =&
\|a_0\| + \sum_{i=1}^r\|a_i\|\cdot\E\big\{\|X_i^{(n)}\|\big\} \le
\|a_0\| + 4\sum_{i=1}^r\|a_i\|,
\end{eqnarray*}
by application of Lemma~\ref{middelvaerdi af norm}.
Putting $K=4\sum_{i=1}^r\|a_i\|$, we may thus conclude that
\[
-\im G_n(\lambda) \ge
  \frac{1}{\|(\im\lambda)^{-1}\|}\frac{1}{(\|\lambda\|+K)^2}\unit_m.
\]
By Lemma~\ref{imaginaerdel}, this implies that $G_n(\lambda)$ is
invertible and that
\[
\big\|G_n(\lambda)^{ 
-1}\big\|\le(\|\lambda\|+K)^2\cdot\big\|(\im\lambda)^{-1}\big\|,
\]
as desired.
\hfill\qed

\begin{corollary}\label{kor til estimat1}
Let $\lambda$ be a matrix in $M_m(\C)$ such that $\im\lambda$ is
positive definite. Then
\begin{equation}
\Big\| a_0 + \sum_{i=1}^ra_iG_n(\lambda)a_i + G_n(\lambda)^{-1} -
\lambda\Big\| \le
\frac{C}{n^2}(K+\|\lambda\|)^2\big\|(\im\lambda)^{-1}\big\|^5,
\label{e4.3c}
\end{equation}
where{\rm ,} as before{\rm ,} $C=m^3\|\sum_{i=1}^ra_i^2\|^2$ and
$K=\|a_0\| + 4\sum_{i=1}^r\|a_i\|$.
\end{corollary}

\Proof   Note that
\begin{multline*}
a_0 + \sum_{i=1}^ra_iG_n(\lambda)a_i + G_n(\lambda)^{-1} -\lambda\\
 =
\Big(\sum_{i=1}^ra_iG_n(\lambda)a_iG_n(\lambda)+(a_0 -\lambda)  
G_n(\lambda)
+\unit_m\Big)G_n(\lambda)^{-1}.
\end{multline*}
Hence, \eqref{e4.3c} follows by combining Theorem~\ref{master ineq}
with Proposition~\ref{estimat1}.
\Endproof\vskip4pt 

In addition to the given matrices $a_0,\dots,a_r$ in $\Mmsa$,
we consider next, as replacement for the random matrices
$X_1^{(n)},\dots,X_r^{(n)}$, free self-adjoint
operators $x_1,\dots,x_r$ in some $C^*$-probability space
$(\CB,\tau)$. We assume that $x_1,\dots,x_r$ are identically
semi-circular distributed, such that $\tau(x_i)=0$ and $\tau(x_i^2)=1$
for all $i$. Then put
\begin{equation}
\label{eq4-7a}
s=a_0\otimes \unit_{\CB} + \sum_{i=1}^ra_i\otimes x_i\in  
M_m(\C)\otimes\CB.
\end{equation}

Consider further the subset $\CO$ of $M_m(\C)$, given by
\begin{eqnarray}
\CO&=&\{\lambda\in M_m(\C)\mid {\rm Im}(\lambda) \ \textrm{is positive
   definite}\} \label{e4.3a}\\
&= &\{\lambda\in M_m(\C)\mid\lmin(\im\lambda)>0\}
\nonumber
\end{eqnarray}
and for each positive number $\delta$, put
\begin{equation}
\CO_{\delta}=\{\lambda\in\CO\mid \|(\im\lambda)^{-1}\|<\delta\}
=\{\lambda\in\CO\mid \lmin(\im\lambda)>\delta^{-1}\}.
\label{e4.3b}
\end{equation}
Note that $\CO$ and $\CO_{\delta}$ are open subsets of $M_m(\C)$.

If $\lambda\in\CO$, then it follows from Lemma~\ref{imaginaerdel} that
$\lambda\otimes\unit_{\CB}-s$ is invertible, since $s$ is
self-adjoint. Hence, for each $\lambda$ in $\CO$, we may define
\[
G(\lambda)=\id_m\otimes\tau\big[(\lambda\otimes\unit_{\CB}-s)^{-1}\big].
\]
As in the proof of Lemma~\ref{estimat1}, it follows that
$G(\lambda)$ is invertible for any $\lambda$ in $\CO$. Indeed, for
$\lambda$ in $\CO$, we have
\[
\begin{split}
\im&\big((\lambda\otimes\unit_{\CB}-s)^{-1}\big)\\[.2cm]
&=\frac{1}{2\i}\big((\lambda\otimes\unit_{\CB}-s)^{-1}
\big((\lambda^*\otimes\unit_{\CB}-s)-(\lambda\otimes\unit_{\CB}-s)\big)
(\lambda^*\otimes\unit_{\CB}-s)^{-1}\big) \\[.2cm]
&=-(\lambda\otimes\unit_{\CB}-s)^{-1}({\rm Im}(\lambda)\otimes\unit_{\CB})
(\lambda^*\otimes\unit_{\CB}-s)^{-1},
\end{split}
\]
which shows that $-{\rm Im}((\lambda\otimes\unit_{\CB}-s)^{-1})$ is
positive definite and that
\[
\begin{split}
0\le\big(-{\rm Im}((\lambda\otimes\unit_{\CB}-s)^{-1})\big)^{-1} &=
(\lambda^*\otimes\unit_{\CB}-s)((\im\lambda)^{ 
-1}\otimes\unit_{\CB})(\lambda\otimes\unit_{\CB}-s)
\\[.2cm]
&\le\big\|\lambda\otimes\unit_{\CB}-s\big\|^2\big\|(\im\lambda)^{ 
-1}\big\|
\cdot\unit_m\otimes\unit_{\CB}.
\end{split}
\]
Consequently,
\[
-{\rm Im}\big((\lambda\otimes\unit_{\CB}-s)^{-1}\big)\ge
\frac{1}{\|\lambda\otimes\unit_{\CB}-s\|^2\|(\im\lambda)^{-1}\|}
\cdot\unit_m\otimes\unit_{\CB},
\]
so that
\[
-\im G(\lambda)\ge
\frac{1}{\|\lambda\otimes\unit_{\CB}-s\|^2\|(\im\lambda)^{-1}\|}
\cdot\unit_m.
\]
By Lemma~\ref{imaginaerdel}, this implies that $G(\lambda)$ is
invertible and that
\[
\big\|G(\lambda)^{-1}\big\|\le
\big\|(\lambda\otimes\unit_{\CB}-s)\big\|^2\big\|(\im\lambda)^{ 
-1}\big\|.
\]

The following lemma shows that the estimate \eqref{e4.3c} in
Corollary~\ref{kor til estimat1} becomes an exact equation, when
$G_n(\lambda)$ is replaced by $G(\lambda)$.

\begin{lemma}\label{formel for G(lambda)}
With $\CO$ and $G(\lambda)$ defined as above{\rm ,} we have that
\[
a_0 + \sum_{i=1}^ra_iG(\lambda)a_i + G(\lambda)^{-1} = \lambda,
\]
for all $\lambda$ in $\CO$.
\end{lemma}

\Proof   We start by recalling the definition of the R-transform $\CR_s$
of (the distribution of) $s$ with amalgamation over $M_m(\C)$: It can
be shown (cf.\ \cite{V6}) that the expression
\[
G(\lambda)=\id_m\otimes\tau\big[(\lambda\otimes\unit_{\CB}-s)^{-1}\big],
\]
gives rise to a well-defined and bijective mapping on a region of the
form
\[
\CU_{\delta}=\big\{\lambda\in M_m(\C)\mid \lambda \ \textrm{is  
invertible
   and} \ \|\lambda^{-1}\|<\delta\big\},
\]
where $\delta$ is a (suitably small) positive number. Denoting by
$G^\brinv$ the inverse of the mapping $\lambda\mapsto G(\lambda)$
$(\lambda\in\CU_{\delta})$, the R-transform $\CR_s$ of $s$ with
amalgamation over $M_m(\C)$ is defined as
\[
\CR_s(\rho)=G^\brinv(\rho)-\rho^{-1}, \qquad (\rho\in
G(\CU_{\delta})).
\]
In \cite{Le} it was proved that
\[
\CR_s(\rho)=a_0 + \sum_{i=1}^ra_i\rho a_i,
\]
so that
\[
G^\brinv(\rho)=a_0 + \sum_{i=1}^ra_i\rho a_i + \rho^{-1}, \qquad  
(\rho\in
G(\CU_{\delta}));
\]
  hence
\begin{equation}
a_0 + \sum_{i=1}^ra_i G(\lambda)a_i + G(\lambda)^{-1} = \lambda, \qquad
(\lambda\in\CU_{\delta}).
\label{e4.4}
\end{equation}
Note now that by Lemma~\ref{imaginaerdel}, the set $\CO_{\delta}$,
defined in \eqref{e4.3b}, is a subset of $\CU_{\delta}$, and hence
\eqref{e4.4} holds, in particular, for $\lambda$ in $\CO_{\delta}$.
Since $\CO_{\delta}$ is an open, nonempty subset of $\CO$ (defined in
\eqref{e4.3a}) and since $\CO$ is a nonempty connected (even convex)  
subset of
$M_m(\C)$, it follows then from the
principle of uniqueness of analytic continuation (for analytical
functions in $m^2$ complex variables) that formula \eqref{e4.4}  
actually holds
for all $\lambda$ in $\CO$, as desired.
\Endproof\vskip4pt 

For $n$ in $\N$ and $\lambda$ in the set $\CO$ (defined in
\eqref{e4.3a}), we introduce further the following notation:
\begin{eqnarray}
\Lambda_n(\lambda)&=&a_0 +  
\sum_{i=1}^ra_iG_n(\lambda)a_i+G_n(\lambda)^{-1},
\label{e4.4d}
\\[.2cm]
\varepsilon(\lambda)&=&\frac{1}{\|(\im\lambda)^{-1}\|} \ = \
\lmin(\im\lambda),
\label{e4.4c}
\\[.2cm]
\CO_n'&=&\big\{\lambda\in\CO\bigm |
\textstyle{\frac{C}{n^2}(K+\|\lambda\|)^2\varepsilon(\lambda)^{ 
-6}<\frac{1}{2}}\big\},
\label{e4.4b}
\end{eqnarray}
where, as before, $C=\frac{\pi^2}{8}m^3\|\sum_{i=1}^ra_i^2\|^2$ and
$K=\|a_0\| +  4\sum_{i=1}^r\|a_i\|$. Note that $\CO'_n$ is an open  
subset of
$M_m(\C)$, since the mapping $\lambda\mapsto\varepsilon(\lambda)$ is
continuous on $\CO$. With the above notation we have the following

\begin{lemma}\label{formel for Gn(lambda)}
For any positive integer $n$ and any matrix $\lambda$ in $\CO_n',$ 
\begin{equation}
\im\Lambda_n(\lambda)\ge\frac{\varepsilon(\lambda)}{2}\unit_m.
\label{e4.4a}
\end{equation}
In particular{\rm ,} $\Lambda_n(\lambda)\in\CO$. Moreover
\begin{equation}
a_0 + \sum_{i=1}^ra_iG(\Lambda_n(\lambda))a_i +  
G(\Lambda_n(\lambda))^{-1} =
a_0 + \sum_{i=1}^ra_iG_n(\lambda)a_i + G_n(\lambda)^{-1},
\label{e4.5}
\end{equation}
for any $\lambda$ in $\CO_n'$.
\end{lemma}

\Proof   Note that the right-hand side of \eqref{e4.5} is nothing else
than $\Lambda_n(\lambda)$. Therefore, \eqref{e4.5} follows from
Lemma~\ref{formel for G(lambda)}, once we have established that
$\Lambda_n(\lambda)\in\CO$ for all $\lambda$ in $\CO_n'$. This, in
turn, is an immediate consequence of \eqref{e4.4a}. It suffices thus
to verify \eqref{e4.4a}.
Note first that for any $\lambda$ in $\CO$, we have by
Corollary~\ref{kor til estimat1} that
\[
\begin{split}
\big\|\im\Lambda_n(\lambda)-\im\lambda\big\| \le
\big\|\Lambda_n(\lambda)-\lambda\big\|
&=\Big\|a_0 + \sum_{i=1}^ra_iG_n(\lambda)a_i + G_n(\lambda)^{-1}
-\lambda\Big\| \\[.2cm]
&\le \frac{C}{n^2}(K+\|\lambda\|)^2\varepsilon(\lambda)^{-5}.
\end{split}
\]
\newpage
\noindent
In particular, $\im\Lambda_n(\lambda)-\im\lambda\ge
-\frac{C}{n^2}(K+\|\lambda\|)^2\varepsilon(\lambda)^{-5}\unit_m$, and
since also $\im\lambda\ge\varepsilon(\lambda)\unit_m$, by definition of
$\varepsilon(\lambda)$, we conclude that
\begin{equation}
\im\Lambda_n(\lambda) = \im\lambda + (\im\Lambda_n(\lambda)-\im\lambda)
\ge\big(\varepsilon(\lambda)
-\textstyle{\frac{C}{n^2}}(K+\|\lambda\|)^2\varepsilon(\lambda)^{ 
-5}\big)\unit_m,
\label{e4.6}
\end{equation}
for any $\lambda$ in $\CO$. Assume now that $\lambda\in\CO_n'$. Then
$\frac{C}{n^2}(K+\|\lambda\|)^2\varepsilon(\lambda)^{ 
-5}<\frac{1}{2}\varepsilon(\lambda)$,
and inserting this in \eqref{e4.6}, we find that
\[
\im\Lambda_n(\lambda)\ge\textstyle{\frac{1}{2}}\varepsilon(\lambda)\unit 
_m,
\]
as desired.
\hfill\qed

\begin{proposition}\label{sammenhaeng mellem G og Gn}
Let $n$ be a positive integer. Then with  $G${\rm ,} $G_n$ and $\CO_n'$ as
defined above{\rm ,} we have that
\[
G(\Lambda_n(\lambda))=G_n(\lambda),
\]
for all $\lambda$ in $\CO'_n$.
\end{proposition}

\Proof   Note first that the functions $\lambda\mapsto G_n(\lambda)$ and
$\lambda\mapsto G(\Lambda_n(\lambda))$ are both analytical functions
(of $m^2$ complex variables) defined on $\CO_n'$ and taking values in
$M_m(\C)$. Applying the principle of uniqueness of analytic
continuation, it suffices thus to prove the following two
assertions:

\begin{itemize}

\item[(a)] The set $\CO_n'$ is an open connected subset of $M_m(\C)$.

\item[(b)] The formula $G(\Lambda_n(\lambda))=G_n(\lambda)$ holds for
   all $\lambda$ in some open, nonempty subset $\CO_n''$ of
   $\CO_n'$.

\end{itemize}

{\it Proof of} (a).   We have already noted that $\CO_n'$ is
open. Consider the subset $I_n$ of $\R$ given by:
\[
I_n=\big\{t\in{}]0,\infty[ \bigm |
\textstyle{\frac{C}{n^2}}(K+t)^2t^{-6}<\frac{1}{2}\big\},
\]
with $C$ and $K$ as above. Note that since the function
$t\mapsto(K+t)^2t^{-6}$ $(t>0)$ is continuous and strictly
decreasing, $I_n$ has the form: $I_n={}]t_n,\infty[$, where $t_n$ is
uniquely determined by the equation:
$\textstyle{\frac{C}{n^2}}(K+t)^2t^{-6}=\frac{1}{2}$.
Note further that for any $t$ in $I_n$, $\i t\unit_m\in\CO_n'$,
and hence the set
\[
\CI_n=\{\i t\unit_m\mid t\in I_n\},
\]
is an arc-wise connected subset of $\CO_n'$. To prove (a), it suffices
then to show that any $\lambda$ in $\CO_n'$ is connected to some
point in $\CI_n$ via a continuous curve $\gamma_{\lambda}$, which is
entirely contained in $\CO_n'$. So let $\lambda$ from $\CO_n'$ be
given, and note that
$0\le\varepsilon(\lambda)=\lmin(\im\lambda)\le\|\lambda\|$. Thus,
\[
\frac{C}{n^2}(K+\varepsilon(\lambda))^2\varepsilon(\lambda)^{-6}\le
\frac{C}{n^2}(K+\|\lambda\|)^2\varepsilon(\lambda)^{-6}<\frac{1}{2},
\]
and therefore $\varepsilon(\lambda)\in I_n$ and
$\i\varepsilon(\lambda)\unit_m\in\CI_n$. Now, let
$\gamma_{\lambda}\colon[0,1]\to M_m(\C)$ be the straight line from
$\i\varepsilon(\lambda)\unit_m$ to $\lambda$, i.e.,
\[
\gamma_{\lambda}(t)=(1-t)\i\varepsilon(\lambda)\unit_m+t\lambda, \qquad
(t\in[0,1]).
\]
We show that $\gamma_{\lambda}(t)\in\CO_n'$ for all $t$ in
$[0,1]$. Note for this that
\[
\im\gamma_{\lambda}(t) = (1-t)\varepsilon(\lambda)\unit_m+t\im\lambda,
\qquad (t\in[0,1]),
\]
so obviously $\gamma_{\lambda}(t)\in\CO$ for all $t$ in $[0,1]$.  
Furthermore,
if $0\le r_1\le r_2\le\cdots\le r_m$ denote the eigenvalues of
${\rm Im}(\lambda)$, then, for each $t$ in $[0,1]$,
$(1-t)\varepsilon(\lambda)+tr_j$ ($j=1,2,\dots,m)$ are the eigenvalues
of $\im\gamma_{\lambda}(t)$. In particular, since  
$r_1=\varepsilon(\lambda)$,
$\varepsilon(\gamma_{\lambda}(t))=\lmin(\im\gamma_{\lambda}(t))
=\varepsilon(\lambda)$ for all $t$ in $[0,1]$. Note also that
\[
\|\gamma_{\lambda}(t)\|\le (1-t)\varepsilon(\lambda)+t\|\lambda\| \le
(1-t)\|\lambda\|+t\|\lambda\|=\|\lambda\|,
\]
for all $t$ in $[0,1]$. Altogether, we conclude that
\[
\frac{C}{n^2}(K+\|\gamma_{\lambda}(t)\|)^2\varepsilon(\gamma_{\lambda}(t 
))^{-6}
\le \frac{C}{n^2}(K+\|\lambda\|)^2\varepsilon(\lambda)^{-6} <
\frac{1}{2},
\]
and hence $\gamma_{\lambda}(t)\in\CO_n'$ for all $t$ in $[0,1]$, as
desired.

\demo{Proof of {\rm (b)}} Consider, for the moment, a fixed matrix $\lambda$
from $\CO_n'$, and put $\zeta=G_n(\lambda)$ and
$\upsilon=G(\Lambda_n(\lambda))$. Then Lemma~\ref{formel for
   Gn(lambda)} asserts that
\[
a_0 + \sum_{i=1}^ra_i\upsilon a_i + \upsilon^{-1} = a_0 +  
\sum_{i=1}^ra_i\zeta a_i
+ \zeta^{-1},
\]
so that
\[
\upsilon\Big(\sum_{i=1}^ra_i\upsilon a_i + \upsilon^{-1}\Big)\zeta
  =\upsilon\Big( \sum_{i=1}^ra_i\zeta a_i + \zeta^{-1}\Big)\zeta;
\]
  hence
\[
\sum_{i=1}^r\upsilon a_i(\upsilon-\zeta)a_i\zeta = \upsilon-\zeta.
\]
In particular, it follows that
\begin{equation}
\Big(\|\upsilon\|\|\zeta\|\sum_{i=1}^r\|a_i\|^2\Big)\|\upsilon-\zeta\|
\ge \|\upsilon-\zeta\|.
\label{e4.7}
\end{equation}
Note here that by Lemma~\ref{imaginaerdel},
\begin{eqnarray}
\|\zeta\|&=&\|G_n(\lambda)\|=
\big\|\id_m\otimes\tr_n\big[(\lambda\otimes\unit_n-S_n)^{-1}\big]\big\|
\label{e4.8}\\[.2cm]
&\le& \big\|(\lambda\otimes\unit_n-S_n)^{-1}\big\|
\le\big\|(\im\lambda)^{-1}\big\|=\frac{1}{\varepsilon(\lambda)}.
\nonumber
\end{eqnarray}
\newpage
\noindent
Similarly, it follows that
\begin{equation}
\|\upsilon\|=\|G(\Lambda_n(\lambda))\|\le
\big\|(\Lambda_n(\lambda)\otimes\unit_{\CB}-s)^{-1}\big\|
\le
\big\|(\im\Lambda_n(\lambda))^{ 
-1}\big\|\le\frac{2}{\varepsilon(\lambda)},
\label{e4.9}
\end{equation}
where the last inequality follows from \eqref{e4.4a} in
Lemma~\ref{formel for Gn(lambda)}. Combining
\eqref{e4.7}--\eqref{e4.9}, it follows that
\begin{equation}
\Big(\frac{2}{\varepsilon(\lambda)^2}\sum_{i=1}^r\|a_i\|^2\Big)\|\upsilon-\zeta\|
\ge \|\upsilon-\zeta\|.
\label{e4.10}
\end{equation}
This estimate holds for all $\lambda$ in $\CO_n'$. If
$\lambda$ satisfies, in addition, that
$\frac{2}{\varepsilon(\lambda)^2}\sum_{i=1}^r\|a_i\|^2\break<1$, then
\eqref{e4.10} implies that $\zeta=\upsilon$, i.e.,
$G_n(\lambda)=G(\Lambda_n(\lambda))$. Thus, if we put
\[
\CO_n''=\big\{\lambda\in\CO_n'\bigm |
\varepsilon(\lambda)>\textstyle{\sqrt{2\sum_{i=1}^r\|a_i\|^2}}\big\},
\]
we have established that $G_n(\lambda)=G(\Lambda_n(\lambda))$ for all
$\lambda$ in $\CO_n''$. Since $\varepsilon(\lambda)$ is a continuous
function of $\lambda$, $\CO_n''$ is clearly an open subset of
$\CO_n'$, and it remains to check that $\CO_n''$ is nonempty. Note,
however, that for any positive number $t$, the matrix $\i t\unit_m$
is in $\CO$ and it satisfies that $\|\i t\unit_m\|=\varepsilon(\i
t\unit_m)=t$. From this, it follows easily that $\i
t\unit_m\in\CO_n''$ for all sufficiently large positive numbers
$t$. This concludes the proof of (b) and hence the proof of
Proposition~\ref{sammenhaeng mellem G og Gn}.
\hfill\qed

\begin{theorem}\label{estimat af Gn(lambda)-G(lambda)}
Let $r,m$ be positive integers{\rm ,} let $a_1,\dots,a_r$ be self\/{\rm -}\/adjoint
matrices in $M_m(\C)$ and{\rm ,} for each positive integer $n${\rm ,} let
$X_1^{(n)},\dots,X_r^{(n)}$ be independent random matrices in
$\SGRM(n,\frac{1}{n})$. Consider further free self\/{\rm -}\/adjoint
identically semi\/{\rm -}\/circular distributed operators $x_1,\dots,x_r$ in
some $C^*$\/{\rm -}\/probability space $(\CB,\tau)${\rm ,} and normalized such that
$\tau(x_i)=0$ and $\tau(x_i^2)=1$ for all $i$. Then put as in
\eqref{e2.0a} and \eqref{eq4-7a}\/{\rm :}\/
\begin{eqnarray*}
s&=&a_0 \otimes\unit_{\CB} + \sum_{i=1}^ra_i\otimes x_i\in
  M_m(\C)\otimes\CB\\
S_n &=&a_0\otimes\unit_n + \sum_{i=1}^ra_i\otimes X_i^{(n)}\in
  M_m(\C)\otimes M_n(\C), \quad (n\in\N),
\end{eqnarray*}
and for $\lambda$ in $\CO=\{\lambda\in M_m(\C)\mid {\rm Im}(\lambda) \
\textrm{is positive definite}\}$ define
\begin{eqnarray*}
G_n(\lambda)&=&\E\big\{ (\id_m\otimes\tr_n)
\big[(\lambda\otimes\unit_n-S_n)^{-1}\big]\big\}\\
G(\lambda)&=&(\id_m\otimes\tau)
\big[(\lambda\otimes\unit_{\CB}-s)^{-1}\big].
\end{eqnarray*}
Then{\rm ,} for any $\lambda$ in $\CO$ and any positive integer $n${\rm ,} we have
\begin{equation}
\big\|G_n(\lambda)-G(\lambda)\big\|\le
\frac{4C}{n^2}(K+\|\lambda\|)^2\big\|(\im\lambda)^{-1}\big\|^7,
\label{e4.11}
\end{equation}
where $C=m^3\|\sum_{i=1}^r a_i^2\|^2$ and
$K=\|a_0\| + 4\sum_{i=1}^r\|a_i\|$.
\end{theorem}

\Proof   Let $n$ in $\N$ be fixed, and assume first that $\lambda$ is in
the set $\CO_n'$ defined in \eqref{e4.4b}. Then, by
Proposition~\ref{sammenhaeng mellem G og Gn}, we have
\[
\begin{split}
\big\|G_n(\lambda)-G(\lambda)\big\| &=
\big\|G(\Lambda_n(\lambda))-G(\lambda)\big\| \\[.2cm]
&=\big\|\id_m\otimes\tau\big[(\Lambda_n(\lambda)\otimes\unit_{\CB}- 
s)^{-1}
-(\lambda\otimes\unit_{\CB}-s)^{-1}\big]\big\| \\[.2cm]
&\le\big\|(\Lambda_n(\lambda)\otimes\unit_{\CB}-s)^{-1}
-(\lambda\otimes\unit_{\CB}-s)^{-1}\big\|.
\end{split}
\]
Note here that
\begin{multline*}
(\Lambda_n(\lambda)\otimes\unit_{\CB}-s)^{-1}
-(\lambda\otimes\unit_{\CB}-s)^{-1}
\\[.2cm]
=(\Lambda_n(\lambda)\otimes\unit_{\CB}-s)^{-1}
\big((\lambda-\Lambda_n(\lambda)\otimes\unit_n\big)
(\lambda\otimes\unit_{\CB}-s)^{-1},
\end{multline*}
and therefore, taking Lemma~\ref{imaginaerdel} into account,
\[
\begin{split}
\big\|G_n(\lambda)-G(\lambda)\big\|&\le
\big\|(\Lambda_n(\lambda)\otimes\unit_{\CB}-s)^{-1}\big\|\cdot
\big\|\lambda-\Lambda_n(\lambda)\big\|\cdot
\big\|(\lambda\otimes\unit_{\CB}-s)^{-1}\big\| \\[.2cm]
&\le \big\|(\im\Lambda_n(\lambda))^{-1}\big\|\cdot
\big\|\lambda-\Lambda_n(\lambda)\big\|\cdot
\big\|(\im\lambda)^{-1}\big\|.
\end{split}
\]
Now, $\|(\im\lambda)^{-1}\|=1/\varepsilon(\lambda)$ (cf.\  
\eqref{e4.4c}),
and hence, by \eqref{e4.4a} in Lemma~\ref{formel for Gn(lambda)},
$\|(\im\Lambda_n(\lambda))^{-1}\|\le2/\varepsilon(\lambda)=
2\|(\im\lambda)^{-1}\|$.
Furthermore, by \eqref{e4.4d} and Corollary~\ref{kor til estimat1},
\begin{eqnarray*}
\big\|\Lambda_n(\lambda)-\lambda\big\|& =&
\Big\|a_0+  
\sum_{i=1}^ra_iG_n(\lambda)a_i+G_n(\lambda)^{-1}-\lambda\Big\| \\
&\le&
\frac{C}{n^2}(K+\|\lambda\|)^2\big\|(\im\lambda)^{-1}\big\|^5.
\end{eqnarray*}
Thus, we conclude that
\[
\big\|G_n(\lambda)-G(\lambda)\big\|\le
\frac{2C}{n^2}(K+\|\lambda\|)^2\big\|(\im\lambda)^{-1}\big\|^7,
\]
which shows, in particular, that \eqref{e4.11} holds for all $\lambda$
in $\CO_n'$.

Assume next that $\lambda\in\CO\setminus\CO_n'$, so that
\begin{equation}
\frac{C}{n^2}(K+\|\lambda\|)^2\big\|(\im\lambda)^{-1}\big\|^6=
\frac{C}{n^2}(K+\|\lambda\|)^2\varepsilon(\lambda)^{-6}\ge\frac{1}{2}.
\label{e4.12}
\end{equation}
By application of Lemma~\ref{imaginaerdel}, it follows that
\begin{equation}
\big\|G(\lambda)\big\|\le
\big\|(\lambda\otimes\unit_{\CB}-s)^{-1}\big\|\le
\big\|(\im\lambda)^{-1}\big\|,
\label{e4.13}
\end{equation}
and similarly we find that
\[
\big\|\id_m\otimes\tr_n\big[(\lambda\otimes\unit_n-S_n(\omega))^{ 
-1}\big]\big\|
\le\big\|(\im\lambda)^{-1}\big\|,
\]
at all points $\omega$ in $\Omega$. Hence, after integrating with respect to 
$\omega$ and using Jensen's inequality,
\begin{equation}
\|G_n(\lambda)\|\le \E\big\{\big\|\id_m\otimes\tr_n
\big[(\lambda\otimes\unit_n-S_n)^{-1}\big]\big\|\big\}
\le\big\|(\im\lambda)^{-1}\big\|.
\label{e4.14}
\end{equation}

Combining \eqref{e4.12}--\eqref{e4.14}, we find that
\begin{multline*}
\big\|G_n(\lambda)-G(\lambda)\big\|\le 2\big\|(\im\lambda)^{-1}\big\|
\\
=\frac{1}{2}\cdot4\big\|(\im \lambda)^{-1}\big\|
\le \frac{4C}{n^2}(K+\|\lambda\|)^2\big\|(\im\lambda)^{-1}\big\|^7,
\end{multline*}
verifying that \eqref{e4.11} holds for $\lambda$ in
$\CO\setminus\CO_n'$ too.
\hfill\qed

\section{The spectrum of $S_n$}\label{sec5}

Let $r,m\in\N$, let $a_0,\dots,a_r\in M_m(\C)_{\rm sa}$ and for each
$n\in\N$, let $X_1^{(n)},\dots\break\dots,X_r^{(n)}$ be $r$ independent random
matrices in $\SGRM(n,\frac1n)$. Let further $x_1,\dots,x_r$ be a
semi-circular family in a $C^*$-probability space $(\CB,\tau)$, and
define $S_n$, $s$, $G_n(\lambda)$ and $G(\lambda)$ as in
Theorem~\ref{estimat af Gn(lambda)-G(lambda)}.

\begin{lemma} \label{cor5-2}
For $\lambda\in\C$ with $\im\lambda >0${\rm ,} put
\begin{equation}
\label{eq5-6}
g_n(\lambda) =  
\E\big\{(\tr_m\otimes\tr_n)[(\lambda\unit_{mn}-S_n)^{-1}]\big\}
\end{equation}
and
\begin{equation}
\label{eq5-7}
g(\lambda) =  
(\tr_m\otimes\tau)\big[(\lambda(\unit_{m}\otimes\unit_{\CB})-s)^{ 
-1}\big].
\end{equation}
Then
\begin{equation}
\label{eq5-8}
|g_n(\lambda)-g(\lambda)| \le \frac{4C}{n^2}
\big(K+|\lambda|\big)^2(\im\lambda)^{-7}
\end{equation}
where $C${\rm ,} $K$ are the constants defined in Theorem~{\rm \ref{estimat af  
Gn(lambda)-G(lambda)}.}
\end{lemma}

\Proof  
This is immediate from Theorem~\ref{estimat af Gn(lambda)-G(lambda)}  
because
\[
g_n(\lambda) = \tr_m(G_n(\lambda \unit_m))
\]
and
\vglue4pt
\hfill $
\displaystyle{g(\lambda) = \tr_m(G(\lambda \unit_m)). }
$\Endproof\vskip12pt

Let Prob$(\R)$ denote the set of Borel probability measures on $\R$. We  
equip
Prob$(\R)$ with the weak$^*$-topology given by $C_0(\R)$, i.e., a net
$(\mu_\alpha)_{\alpha\in A}$ in Prob$(\R)$ converges in  
weak$^*$-topology
to $\mu\in\mbox{Prob}(\R)$, if and only if
\[
\lim_\alpha\bigg(\int_\R\varphi\6\mu_\alpha\bigg) = \int_\R \varphi\6\mu
\]
for all $\varphi\in C_0(\R)$.

Since $S_n$ and $s$ are self-adjoint, there are, by Riesz'  
representation
theorem, unique probability measures
$\mu_n$, $n=1,2,\dots$ and $\mu$ on $\R$, such that
\begin{eqnarray}
\label{eq5-9}
\int_\R \varphi\6\mu_n &=&
\E\big\{(\tr_m\otimes\tr_n)\varphi(S_n)\big\}\\ 
\label{eq5-10}
\int_\R \varphi\6\mu &=& (\tr_m\otimes\tau)\varphi(s)
\end{eqnarray}
for all $\varphi\in C_0(\R)$. Note that $\mu$ is compactly supported
while $\mu_n$, in general, is not compactly supported.

\begin{theorem} \label{thm5-3}
Let $S_n$ and $s$ be given by \eqref{e2.0a} and \eqref{eq4-7a}{\rm ,} and let
$C=\frac{\pi^2}{8} m^3\|\sum^r_{i=1} a_i^2\|^2$ and
$K=\|a_0\|+4\sum^r_{i=1} \|a_i\|$. Then for all $\varphi\in
C_c^\infty(\R,\R)${\rm ,}
\begin{equation}
\label{eq5-11}
\E\big\{(\tr_m\otimes\tr_n)\varphi(S_n)\big\} =
(\tr_m\otimes\tau)\varphi(s)+R_n
\end{equation}
where
\begin{equation}
\label{eq5-12}
|R_n| \le\frac{4C}{315\pi
  n^2}\int_\R\big|((1+D)^8\varphi)(x)\big|\big(K+2+|x|\big)^2 \6x
\end{equation}
and $D=\frac{\rm d}{{\rm d}x}$. In particular $R_n=O(\frac{1}{n^2})$ for
$n\to\infty$.
\end{theorem}

\Proof  
Let $g_n,g,\mu_n,\mu$ be as in \eqref{eq5-6}, \eqref{eq5-7},
\eqref{eq5-9} and \eqref{eq5-10}. Then for any complex number
$\lambda$, such that ${\rm Im}(\lambda)>0$, we have
\begin{eqnarray}
\label{eq5-13}
g_n(\lambda) &=& \int_\R \frac{1}{\lambda -x}\6\mu_n(x)\\
\label{eq5-14}
g(\lambda) &=& \int_\R \frac{1}{\lambda-x}\6\mu(x).
\end{eqnarray}
Hence $g_n$ and $g$ are the Stieltjes transforms (or Cauchy transforms,  
in the
terminology of \cite{vdn}) of $\mu_n$ and $\mu$ in the half plane
$\im\lambda>0$. Hence, by the inverse Stieltjes transform,
\[
\mu_n = \lim_{y\to 0^+}\Big(-\frac{1}{\pi}{\rm Im}(g_n(x+\i{}y))\6x\Big)
\]
where the limit is taken in the weak$^*$-topology on Prob$(\R)$. In
particular, for all $\varphi$ in $C_c^\infty(\R,\R)$:
\begin{equation}
\label{eq5-15}
\int_\R \varphi(x) \6\mu_n(x) = \lim_{y\to 0^+} \Big[-\frac1\pi
\im\Big(\int_\R \varphi(x)g_n(x+\i{}y)\6x\Big)\Big].
\end{equation}
In the same way we get for $\varphi\in C_c^\infty(\R,\R)$:
\begin{equation}
\label{eq5-16}
\int_\R \varphi(x) \6\mu(x) = \lim_{y\to 0^+} \Big[-\frac1\pi
\im\int_\R \varphi(x)g(x+\i{}y)\6x\Big].
\end{equation}

In the rest of the proof, $n\in\N$ is fixed, and we put
$h(\lambda)=g_n(\lambda)-g(\lambda)$.
Then by \eqref{eq5-15} and \eqref{eq5-16}
\begin{equation}
\label{eq5-17}
\Big| \int_\R \varphi(x)\6\mu_n(x)-\int_{\R}\varphi(x)\6\mu(x)\Big|
\le\frac1\pi \limsup_{y\to 0^+} \Big|\int_\R  
\varphi(x)h(x+\i{}y)\6x\Big|.
\end{equation}
For $\im\lambda >0$ and $p\in\N$, put
\begin{equation}
\label{eq5-18}
I_p(\lambda) = \frac{1}{(p-1)!} \int^\infty_0 h(\lambda
+t)t^{p-1}\e^{-t}\6t.
\end{equation}
\newpage
\noindent
Note that $I_p(\lambda)$ is well defined because, by \eqref{eq5-13} and
\eqref{eq5-14}, $h(\lambda)$ is uniformly bounded in any half-plane of
the form $\im\lambda\ge\varepsilon$, where $\varepsilon >0$. Also, it is
easy to check that $I_p(\lambda)$ is an analytic function of $\lambda$,
and its first derivative is given by
\begin{equation}
\label{eq5-19}
I_p'(\lambda) = \frac{1}{(p-1)!} \int^\infty_0 h'(\lambda +t)t^{p-1}
\e^{-t} \6t
\end{equation}
where $h'=\frac{{\rm d}h}{{\rm d}\lambda}$. We claim that
\begin{eqnarray}
\label{eq5-20}
I_1(\lambda)-I'_1(\lambda) &=& h(\lambda)\\[.2cm]
\label{eq5-21}
I_p(\lambda)-I'_p(\lambda) &=& I_{p-1}(\lambda),\quad p\ge 2.
\end{eqnarray}
Indeed, by \eqref{eq5-19} and partial integration we get
\begin{eqnarray*}
I_1'(\lambda) &=& \big[h(\lambda +t)\e^{-t}\big]^\infty_0 +  
\int^\infty_0
h(\lambda +t)\e^{-t}\6t\\[.2cm]
&=& -h(\lambda)+I_1(\lambda),
\end{eqnarray*}
which proves \eqref{eq5-20} and in the same way we get for $p\ge 2$,
\begin{eqnarray*}
I'_p(\lambda) &=& \frac{1}{(p-1)!} \int^\infty_0 h'(\lambda
+t)t^{p-1}\e^{-t} \6t\\[.2cm]
&=& -\frac{1}{(p-1)!} \int^\infty_0 h(\lambda
+t)((p-1)t^{p-2}-t^{p-1})\e^{-t}\6t\\[.2cm]
&=& -I_{p-1} (\lambda)+I_p(\lambda),
\end{eqnarray*}
which proves \eqref{eq5-21}. Assume now that $\varphi\in  
C_c^\infty(\R,\R)$
and that $y>0$. Then, by \eqref{eq5-20} and partial integration, we have
\begin{eqnarray*}
\int_\R\varphi(x)h(x+\i{}y)\6x &=& \int_\R \varphi(x)I_1(x+\i{}y)\6x -  
\int_\R
\varphi(x)I'_1(x+\i{}y)\6x\\[.2cm]
&=& \int_\R \varphi(x)I_1(x+\i{}y)\6x + \int_\R
\varphi'(x)I_1(x+\i{}y)\6x\\[.2cm]
&=& \int_\R ((1+D)\varphi)(x)\cdot I_1(x+\i{}y)\6x,
\end{eqnarray*}
where $D=\frac{\rm d}{{\rm d}x}$. Using \eqref{eq5-21}, we can continue  
to
perform partial integrations, and after $p$ steps we obtain
\[
\int_\R \varphi(x)h(x+\i{}y)\6x = \int_\R ((1+D)^p\varphi)(x)\cdot
I_p(x+\i{}y)\6x.
\]
Hence, by \eqref{eq5-17}, we have for all $p\in\N$:
\begin{multline}
\label{eq5-22}
\Big|\int_\R \varphi(x)\6\mu_n(x)-\int_{\R}\varphi(x)\6\mu(x)\Big|
\\
\le\frac1\pi \limsup_{y\to 0^+} \Big|\int_\R
((1+D)^p\varphi)(x)\cdot I_p(x+\i{}y)\6x\Big|.
\end{multline}

Next, we use \eqref{eq5-8} to show that for $p=8$ and $\im\lambda >0$  
one
has
\begin{equation}
\label{eq5-23}
|I_8(\lambda)| \le \frac{4C(K+2+|\lambda|)^2}{315n^2}.
\end{equation}
To prove \eqref{eq5-23}, we apply Cauchy's integral theorem to the
function
\[
F(z) = \frac{1}{7!} h(\lambda +z)z^7\e^{-z},
\]
which is analytic in the half-plane $\im z>-\im\lambda$. Hence for $r>0$
\[
\int_{[0,r]}F(z)\6z + \int_{[r,r+\i{}r]} F(z)\6z + \int_{[r+\i{}r,0]}
F(z)\6z=0
\]
where $[\alpha,\beta]$ denotes the line segment connecting $\alpha$ and
$\beta$ in $\C$ oriented from $\alpha$ to $\beta$. Put
\[
M(\lambda) = \sup\big\{ |h(w)| \bigm | \im w \ge \im\lambda\big\}.
\]
Then by \eqref{eq5-13} and \eqref{eq5-14}, $M(\lambda)
\le\frac{2}{|\im\lambda|}<\infty$. Hence
\begin{eqnarray*}
\Big|\int_{[r,r+\i{}r]} F(z)\6z\Big| &\le & \frac{M(\lambda)}{7!}
\int^r_0|r+\i{}t|^7 \e^{-r}\6t\\[.2cm]
&\le & \frac{M(\lambda)}{7!} (2r)^7 r\cdot \e^{-r}\\[.2cm]
&\to & 0, \qquad\mbox{for $r\to\infty$}.
\end{eqnarray*}
Therefore,
\begin{eqnarray}
I_8(\lambda) &=& \frac{1}{7!}\int^\infty_0 h(\lambda +t)t^7 \e^{-t}\6t
 \\[.2cm]
&=& \lim_{r\to\infty} \int_{[0,r]} F(z)\6z \nonumber \\[.2cm]
&=& \lim_{r\to\infty} \int_{[0,r+\i{}r]} F(z)\6z \nonumber \\[.2cm]
\label{eq5-24}
&=& \frac{1}{7!} \int^\infty_0 h(\lambda +(1+\i)t)((1+\i)t)^7
\e^{-(1+\i)t}(1+\i)\6t.\nonumber
\end{eqnarray}
By \eqref{eq5-8},
\[
|h(w)| \le\frac{4C}{n^2} (K+|w|)^2(\im w)^{-7},\qquad \im w > 0.
\]
Inserting this in \eqref{eq5-24} we get
\begin{eqnarray*}
|I_8(\lambda)| &\le & \frac{4C}{7!n^2}\int^\infty_0
\frac{\big(K+|\lambda|+\sqrt{2} t\big)^2}{(\im\lambda +t)^7} (\sqrt{2}  
t)^7
\e^{-t}\sqrt{2}\6t\\[.2cm]
&\le & \frac{2^6 C}{7! n^2}\int^\infty_0\big(K+|\lambda|+\sqrt{2}  
t\big)^2
\e^{-t}\6t\\
&=& \frac{4C}{315
   n^2}\big((K+|\lambda|)^2+2\sqrt{2}(K+|\lambda|)+4\big)\\[.2cm]
&\le & \frac{4C}{315n^2} (K+|\lambda|+2)^2.
\end{eqnarray*}
This proves \eqref{eq5-23}. Now, combining \eqref{eq5-22} and
\eqref{eq5-23}, we have
\[
\begin{split}
\Big|\int_\R\varphi(x)\6\mu_n(x)&-\int_{\R}\varphi(x)\6\mu(x)\Big|
\\[.2cm]
&\le\frac{4C}{315\pi n^2} \limsup_{y\to 0^+}\int_\R
\big|((1+D)^8\varphi)(x)\big|\big(K+2+|x+\i{}y|\big)^2\6x\\[.2cm]
&=\frac{4C}{315\pi n^2}
\int_\R\big|((1+D)^8\varphi)(x)\big|\big(K+2+|x|\big)^2\6x
\end{split}
\]
for all $\varphi\in C_c^\infty(\R,\R)$. Together with \eqref{eq5-9} and
\eqref{eq5-10} this proves Theorem \ref{thm5-3}.
\phantom{endofline}\hfill\qed

\begin{lemma}\label{lemma5-4}
Let $S_n$ and $s$ be given by \eqref{e2.0a} and \eqref{eq4-7a}{\rm ,} and let
$\varphi\colon\R\to\R$ be a $C^\infty$-function which is constant
outside a compact subset of $\R$. Assume further that
\begin{equation}
\label{eq5-25}
\supp(\varphi)\cap \spe(s)=\emptyset .
\end{equation}
Then
\begin{eqnarray}
\label{eq5-26}
\E\big\{(\tr_m\otimes\tr_n)\varphi(S_n)\big\} &=&
O\big(\textstyle{\frac{1}{n^2}}\big),\qquad\mbox{for  
$n\to\infty$}\\[.2cm]
\label{eq5-27}
\V\big\{(\tr_m\otimes\tr_n)\varphi(S_n)\big\} &=&
O\big(\textstyle{\frac{1}{n^4}}\big),\qquad\mbox{for $n\to\infty$}
\end{eqnarray}
where $\V$ is the absolute variance of a complex random variable
\/{\rm (}\/cf.~{\rm \S\ref{sec3}).} Moreover
\begin{equation}
\label{eq5-28}
(\tr_m\otimes\tr_n)\varphi(S_n(\omega)) = O(n^{-4/3})
\end{equation}
for almost all $\omega$ in the underlying probability space $\Omega$.
\end{lemma}

\Proof  
By the assumptions, $\varphi=\psi+c$, for some $\psi$ in  
$C_c^\infty(\R,\R)$
and some constant $c$ in $\R$. By Theorem \ref{thm5-3}
\[
\E\big\{(\tr_m\otimes\tr_n)\psi(S_n)\big\} =
(\tr_m\otimes\tau)\psi(s)+O\big(\textstyle{\frac{1}{n^2}}\big), \qquad
\mbox{for $n\to\infty$},
\]
and hence also
\[
\E\big\{(\tr_m\otimes\tr_n)\varphi(S_n)\big\} =
(\tr_m\otimes\tau)\varphi(s)+O\big(\textstyle{\frac{1}{n^2}}\big),
\qquad\mbox{for $n\to\infty$}.
\]
\enlargethispage*{1000pt}
\newpage
\noindent
But since $\varphi$ vanishes on $\spe(s)$, we have $\varphi(s)=0$. This
proves \eqref{eq5-26}. Moreover, applying Proposition~\ref{diff-kvotient
   trick} to $\psi\in C_c^\infty(\R)$, we have
\begin{equation}
\label{eq5-29}
\V\big\{(\tr_m\otimes\tr_n)\psi(S_n)\big\}\le\frac{1}{n^2}
\Big\|\sum^r_{i=1}
a_i^2\Big\|^2\E\big\{(\tr_m\otimes\tr_n)(\psi'(S_n))^2\big\}.
\end{equation}
By \eqref{eq5-25}, $\psi'=\varphi'$ also vanishes on $\spe(s)$. Hence,
by Theorem \ref{thm5-3}
\[
\E\big\{(\tr_m\otimes\tr_n)|\psi'(S_n)|^2\big\} =
O\big(\textstyle{\frac{1}{n^2}}\big),\qquad\mbox{as $n\to\infty$}.
\]
Therefore, by \eqref{eq5-29}
\[
\V\big\{(\tr_m\otimes\tr_n)\psi(S_n)\big\} =
O\big(\textstyle{\frac{1}{n^4}}\big),\qquad\mbox{as $n\to\infty$.}
\]
Since $\varphi(S_n)=\psi(S_n)+c\unit_{mn}$,
$\V\big\{(\tr_m\otimes\tr_n)\varphi(S_n)\big\}
=\V\big\{(\tr_m\otimes\tr_n)\psi(S_n)\big\}$.
This proves \eqref{eq5-27}. Now put
\begin{eqnarray*}
Z_n &=& (\tr_m\otimes\tr_n)\varphi(S_n)\\[.2cm]
\Omega_n &=& \big\{\omega\in\Omega\bigm | |Z_n(\omega)|\ge  
n^{-4/3}\big\}.
\end{eqnarray*}
By \eqref{eq5-26} and \eqref{eq5-27}
\[
\E\big\{|Z_n|^2\big\} = |\E\{Z_n\}|^2 +\V\{Z_n\}
=O\big(\textstyle{\frac{1}{n^4}}\big), \qquad \mbox{for $n\to\infty$}.
\]
Hence
\begin{equation}
P(\Omega_n)=\int_{\Omega_n}\6P(\omega)\le
\int_{\Omega_n}\big|n^{4/3}Z_n(\omega)\big|^2\6P(\omega)
\le n^{8/3}\E\big\{|Z_n|^2\big\}=O(n^{-4/3}),
\end{equation}
for $n\to\infty$. In particular $\sum^\infty_{n=1}
P(\Omega_n)<\infty$. Therefore, by the
Borel-Cantelli lemma (see e.g.\ \cite{bre}), $\omega\notin\Omega_n$  
eventually, as
$n\to\infty$, for almost all $\omega\in\Omega$; i.e., $|Z_n(\omega)| <
n^{-4 /3}$ eventually, as $n\to\infty$, for almost all
$\omega\in\Omega$. This proves \eqref{eq5-28}.
\Endproof\vskip4pt

\begin{theorem}
\label{thm5-5}
Let $m\in\N$ and let $a_0,\dots,a_r\in M_m(\C)_{\rm sa}${\rm ,} $S_n$ and
$s$ be as in Theorem~{\rm \ref{estimat af Gn(lambda)-G(lambda)}.}
Then for any $\varepsilon>0$ and for almost all $\omega\in\Omega${\rm ,}
\[
\spe(S_n(\omega))\subseteq \spe(s) \ + \ ]-\varepsilon,\varepsilon[,
\]
eventually as $n\to\infty$.
\end{theorem}

\Proof  
Put
\begin{eqnarray*}
K &=& \spe(s)+\big[-\textstyle{\frac{\varepsilon}{2},
\frac{\varepsilon}{2}}\big]\\[.2cm]
F &=& \big\{t\in\R\mid d(t,\spe(s))\ge\varepsilon\big\}.
\end{eqnarray*}
Then $K$ is compact, $F$ is closed and $K\cap F=\emptyset$. Hence there
exists $\varphi\in C^\infty(\R)$, such that $0\le\varphi\le 1$,
$\varphi(t)=0$ for $t\in K$ and $\varphi(t)=1$ for $t\in F$
(cf.~\cite[(8.18) p.~237]{F}).
Since $\C\backslash F$ is a bounded set, $\varphi$ satisfies the
requirements
of lemma \ref{lemma5-4}. Hence by \eqref{eq5-28}, there exists a
$P$-null set $N\subseteq\Omega$, such that for all
$\omega\in\Omega\backslash N$:
\[
(\tr_m\otimes\tr_n)\varphi(S_n(\omega)) = O(n^{-4/3}),\qquad\mbox{as
$n\to\infty$}.
\]
Since $\varphi\ge 1_F$, it follows that
\[
(\tr_m\otimes\tr_n)1_F(S_n(\omega))=O(n^{-4/3}),\qquad\mbox{as
$n\to\infty$}.
\]
But for fixed $\omega\in\Omega\backslash N$, the number of eigenvalues
(counted with multiplicity) of the matrix $S_n(\omega)$ in the set $F$
is equal to $mn(\tr_m\otimes\tr_n)1_F(S_n(\omega))$, which is
$O(n^{-1/3})$ as $n\to\infty$. However, for each $n\in\N$ the above
number is an integer. Hence, the number of eigenvalues of $S_n(\omega)$
in $F$ is zero eventually as $n\to\infty$. This shows that
\[
\spe(S_n(\omega)) \subseteq\C\backslash F =
\spe(s) \ + \ ]-\varepsilon,\varepsilon[
\]
eventually as $n\to\infty$, when $\omega\in\Omega\backslash N$.
\hfill\qed

\section{Proof of the main theorem} \label{sec6}

Throughout this section, $r\in\N\cup\{\infty\}$, and, for each $n$ in
$\N$, we let $(X_i^{(n)})_{i=1}^r$ denote a finite or countable set
of independent random matrices from $\SGRM(n,\frac{1}{n})$, defined on
the same probability space $(\Omega,\CF,P)$. In addition, we let
$(x_i)_{i=1}^r$ denote a corresponding semi-circular family in a
$C^*$-probability space $(\CB,\tau)$, where $\tau$ is a faithful state
on $\CB$. Furthermore, as in \cite{vdn}, we let $\C\<(X_i)_{i=1}^r\>$
denote the algebra of all polynomials in $r$ noncommuting variables.  
Note
that $\C\<(X_i)_{i=1}^r\>$ is a unital $*$-algebra, with the
$*$-operation given by:
\[
(cX_{i_1}X_{i_2}\cdots X_{i_k})^* =
\overline{c}X_{i_k}X_{i_{k-1}}\cdots X_{i_2}X_{i_1},
\]
for $c$ in $\C$, $k$ in $\N$ and $i_1,i_2,\dots,i_k$ in
$\{1,2,\dots,r\}$, when $r$ is finite, and in $\N$ when $r=\infty$.
The purpose of this section is to conclude the proof of the main theorem
(Theorem~\ref{thm6-1} below) by combining the results of the previous
sections.

\begin{theorem}\label{thm6-1}
Let $r$ be in $\N\cup\{\infty\}$. Then there exists a $P$\/{\rm -}\/null\/{\rm -}\/set\break
$N\subseteq\Omega${\rm ,} such that for all $p$ in $\C\<(X_i)_{i=1}^r\>$ and
all $\omega$ in $\Omega\setminus N${\rm ,} we have
\[
\lim_{n\to\infty}\big\|p\big((X_i^{(n)}(\omega))_{i=1}^r\big)\big\|
=\big\|p\big((x_i)_{i=1}^r\big)\big\|.
\]
\end{theorem}
\vglue8pt

We start by proving the following

\begin{lemma}
\label{lemma6-2} Assume that $r\in\N$.
Then there exists a $P$\/{\rm -}\/null set $N_1\subseteq\Omega${\rm ,} such that for
all $p$ in $\C\<(X_i)_{i=1}^r\>$ and all
$\omega$ in $\Omega\backslash N_1$\/{\rm :}\/
\begin{equation}
\label{eq6-2}
\liminf_{n\to\infty}\big\|
p\big(X_1^{(n)}(\omega),\dots,X_r^{(n)}(\omega)\big)\big\|
\ge \|p(x_1,\dots,x_r)\|.
\end{equation}
\end{lemma}
\vskip8pt

\Proof   We first prove that for each $p$ in $\C\<X_1,\dots,X_r\>$,
there exists a $P$-null-set $N(p)$, depending on $p$, such that
\eqref{eq6-2} holds for all $\omega$ in $\Omega\setminus N(p)$. This
assertion is actually a special case of \cite[Prop.~4.5]{T}, but for
the readers convenience, we include a more direct proof:
Consider first a fixed $p\in\C\<X_1,\dots,X_r\>$.
Let $k\in\N$ and put $q=(p^*p)^k$. By \cite[Cor.~3.9]{T} or \cite{HP2},
\begin{equation}
\label{eq6-3}
\lim_{n\to\infty}
\tr_n\big(q(X_1^{(n)}(\omega),\dots,X_r^{(n)}(\omega))\big)
=\tau\big(q(x_1,\dots,x_r)\big),
\end{equation}
for almost all $\omega\in\Omega$. For $s\ge 1$, $Z\in M_n(\C)$ and
$z\in\CB$, put
$\|Z\|_s=\tr_n(|Z|^s)^{1/s}$ and $\|z\|_s=\tau(|z|^s)^{1/s}$.
Then \eqref{eq6-3} can be rewritten as
\begin{equation}
\label{eq6-4}
\lim_{n\to\infty}
\big\|p\big(X_1^{(n)}(\omega),\dots,X_r^{(n)}(\omega)\big)\big\|_{2k}^{2 
k}
=\big\|p(x_1,\dots,x_r)\big\|^{2k}_{2k}
\end{equation}
for $\omega\in\Omega\backslash N(p)$, where $N(p)$ is a $P$-null-set.  
Since
$\N$ is a countable set, we can assume that $N(p)$ does not depend
on $k\in\N$.
For every bounded Borel function $f$ on a probability space, one has
\begin{equation}
\|f\|_{\infty}=\lim_{k\to\infty}\|f\|_k,
\label{eq6-4a}
\end{equation}
(cf.\ \cite[Exercise~7, p.~179]{F}). Put $a=p(x_1,\dots,x_r)$, and
let $\Gamma\colon \CD\to C(\hat{\CD})$ be the Gelfand transform of the
Abelian $C^*$-algebra $\CD$ generated by $a^*a$ and $\unit_{\CB}$, and
let $\mu$ be the probability measure on $\hat{\CD}$ corresponding to
$\tau_{\mid\CD}$. Since $\tau$ is faithful,
$\supp(\mu)=\hat{\CD}$. Hence,
$\|\Gamma(a^*a)\|_{\infty}=\|\Gamma(a^*a)\|_{\sup}=\|a^*a\|$.
Applying then \eqref{eq6-4a} to the function $f=\Gamma(a^*a)$, we find
that
\begin{equation}
\|a\|=\|a^*a\|^{1/2}=\lim_{k\to\infty}\|a^*a\|_k^{1/2}=
\lim_{k\to\infty}\|a\|_{2k}.
\label{eq6-4b}
\end{equation}

Let $\varepsilon>0$. By \eqref{eq6-4b}, we can choose $k$ in $\N$, such  
that
\[
\|p(x_1,\dots,x_r)\|_{2k} > \|p(x_1,\dots,x_r)\|-\varepsilon.
\]
Since $\|Z\|_s\le\|Z\|$ for all $s\ge 1$ and all $Z\in M_n(\C)$, we have
by \eqref{eq6-4}
\[
\liminf_{n\to\infty}
\big\|p\big(X_1^{(n)}(\omega),\dots,X_r^{(n)}(\omega)\big)\big\|
\ge\|p(x_1,\dots,x_r)\|_{2k}
> \|p(x_1,\dots,x_r)\|-\varepsilon,
\]
for all $\omega\in\Omega\backslash N(p)$, and since $N(p)$ does not  
depend on
$\varepsilon$, it follows that \eqref{eq6-2} holds for all
$\omega\in\Omega\backslash N(p)$. Now put $N'=
\bigcup_{p\in\CP}N(p)$, where $\CP$ is the set of polynomials from
$\C\<X_1,\dots,X_r\>$
with coefficients in $\Q+\i\Q$. Then $N'$ is again a null set,
and \eqref{eq6-2} holds for all $p\in\CP$ and all
$\omega\in\Omega\backslash N'$.

By \cite[Thm.~2.12]{Ba} or \cite[Thm.~3.1]{HT1},
$\lim_{n\to\infty}\|X_i^{(n)}(\omega)\| = 2$, $i=1,\dots,r$, for almost  
all
$\omega\in\Omega$. In particular
\begin{equation}
\label{eq6-6}
\sup_{n\in\N} \|X_i^{(n)}(\omega)\| < \infty,\quad i=1,\dots,r,
\end{equation}
for almost all $\omega\in\Omega$. Let $N''\subseteq\Omega$ be the set of
$\omega\in\Omega$ for which \eqref{eq6-6} fails for some
$i\in\{1,\dots,r\}$. Then $N_1=N'\cup N''$ is a null set, and a simple
approximation argument shows that \eqref{eq6-2} holds for all $p$ in
$\C\< X_1,\dots,X_r\>$, when $\omega\in\Omega\backslash N_1$.
\Endproof\vskip4pt 

In order to complete the proof of Theorem~\ref{thm6-1}, we have to prove

\begin{proposition}\label{prop6-3}
Assume that $r\in\N$. Then there
is a $P$\/{\rm -}\/null set\break $N_2\subseteq\Omega${\rm ,} such that for all
polynomials $p$ in $r$ noncommuting variables and all
$\omega\in\Omega\backslash N_2${\rm ,}
\[
\limsup_{n\to\infty}
\big\|p\big(X_1^{(n)}(\omega),\dots,X_r^{(n)}(\omega)\big)\big\|
\le\|p(x_1,\dots,x_r)\|.
\]
\end{proposition}

The proof of Proposition \ref{prop6-3} relies on Theorem \ref{thm5-5}
combined with the linearization trick in the form of  
Theorem~\ref{thm1-4}.
Following the notation of \cite{BK} we put
\[
\prod_n M_n(\C) = \big\{(Z_n)^\infty_{n=1} \bigm | Z_n\in M_n(\C), \
\textstyle{\sup_{n\in\N}} \|Z_n\| <\infty\big\}
\]
and
\[
\sum_n M_n(\C) = \big\{(Z_n)^\infty_{n=1} \bigm | Z_n\in M_n(\C), \
\textstyle{\lim_{n\to\infty}}\|Z_n\|=0\big\},
\]
and we let $\CC$ denote the quotient $C^*$-algebra
\begin{equation}
\label{eq6-7}
\CC = \prod_n M_n(\C)\Big/ \sum_n M_n(\C).
\end{equation}
Moreover, we let $\rho: \prod_n M_n(\C)\to \CC$ denote the quotient
map. By \cite[Lemma~6.13]{RLL},
the quotient norm in $\CC$ is given by
\begin{equation}
\label{eq6-8}
\big\|\rho\big((Z_n)^\infty_{n=1}\big)\big\| = \limsup_{n\to\infty}
\|Z_n\|,
\end{equation}
for $(Z_n)^\infty_{n=1}\in\prod M_n(\C)$.

Let $m\in\N$. Then we can identify $M_m(\C)\otimes \CC$ with $$\prod_n
M_{mn}(\C)\ / \ \sum M_{mn}(\C),$$ where $\prod_n M_{mn}(\C)$ and $\sum_n
M_{mn}(\C)$ are defined as $\prod_n M_n(\C)$ and $\sum_n M_n(\C)$, but
with  $Z_n\in M_{mn}(\C)$ instead of $Z_n\in M_n(\C)$.  
Moreover, for
$(Z_n)^\infty_{n=1}\in\prod_n M_{mn}(\C)$, we have, again by
\cite[Lemma~6.13]{RLL},
\begin{equation}
\label{eq6-9}
\big\|(\id_m\otimes\rho)\big((Z_n)^\infty_{n=1}\big)\big\| =
\limsup_{n\to\infty}\|Z_n\|.
\end{equation}

\begin{lemma}
\label{lemma6-4}
Let $m\in\N$ and let $Z=(Z_n)^\infty_{n=1}\in\prod_n M_{mn}(\C)${\rm ,} such
that each $Z_n$ is normal. Then for all $k\in\N$
\[
\spe\big((\id_m\otimes\rho)(Z)\big)\subseteq\overline{\bigcup^\infty_{n= 
k}
\spe(Z_n)}.
\]
\end{lemma}

\Proof  
Assume $\lambda\in\C$ is not in the closure of $\bigcup^\infty_{n=k}
\spe(Z_n)$. Then there exists an $\varepsilon >0$, such that
$d(\lambda,\spe(Z_n))\ge\varepsilon$ for all $n\ge k$. Since $Z_n$ is
normal, it follows that
$\|(\lambda\unit_{mn}-Z_n)^{-1}\|\le\frac{1}{\varepsilon}$ for all $n\ge
k$. Now put
\[
y_n =
\begin{cases}
0, & \quad \textrm{if} \ 1\le n\le k-1,\\
(\lambda\unit_{mn}-Z_n)^{-1}, &\quad \textrm{if} \ n\ge k.
\end{cases}
\]
Then $y=(y_n)^\infty_{n=1}\in\prod_n M_{mn}(\C)$, and one checks easily
that $\lambda\unit_{M_m(\C)\otimes \CC}-(\id_m\otimes\rho)(Z)$ is
invertible in
$M_m(\C)\otimes \CC = \prod_n M_{mn}(\C)\ / \ \sum_n M_{mn}(\C)$ with
inverse $(\id_m\otimes\rho)y$. Hence $\lambda\notin
\spe((\id_m\otimes\rho)(Z))$.
\Endproof\vskip4pt

{\it Proof of Proposition {\rm \ref{prop6-3}} and Theorem
   {\rm \ref{thm6-1}}}.
Assume first that $r\in\N$. Put
\[
\Omega_0 = \big\{\omega\in\Omega\bigm |
\textstyle{\sup_{n\in\N}}\|X_i^{(n)}(\omega)\|
<\infty, \ i=1,\dots,r\big\}.
\]
By \eqref{eq6-6}, $\Omega\backslash\Omega_0$ is a $P$-null set. For
every $\omega\in\Omega_0$, we define
\[
y_i(\omega)\in \CC=\prod_n M_n(\C)\Big/\sum_n M_n(\C)
\]
by
\begin{equation}
\label{eq6-9a}
y_i(\omega) = \rho\big((X_i^{(n)}(\omega))^\infty_{n=1}\big),\quad
i=1,\dots,r.
\end{equation}
Then for every noncommutative polynomial $p\in\C\< X_1,\dots,X_r\>$ and
every $\omega$ in $\Omega_0$, we
get by \eqref{eq6-8} that
\begin{equation}
\label{eq6-10}
\big\|p(y_1(\omega),\dots,y_r(\omega))\big\| = \limsup_{n\to\infty}
\big\|p\big(X_1^{(n)}(\omega),\dots,X_r^{(n)}(\omega)\big)\big\|.
\end{equation}

Let $j\in\N$ and $a_0,a_1,\dots,a_r\in\Mmsa$. Then by Theorem  
\ref{thm5-5}
there exists a null set $N(m,j,a_0,\dots,a_r)$, such that for
\[
\spe\big(a_0\otimes\unit_n +
\textstyle{\sum^r_{i=1}}a_i\otimes X_i^{(n)}(\omega)\big)
\subseteq \spe\big(a_0\otimes\unit_{\CB}  
+\textstyle{\sum^r_{i=1}}a_i\otimes
x_i\big) \ + \ \big]\textstyle{-\frac1j,\frac1j}\big[,
\]
eventually, as $n\to\infty$, for all $\omega\in\Omega\backslash
N(m,j,a_0,\dots,a_r)$. Let $N_0=\break
\bigcup N(m,j,a_0,\dots,a_r)$, where the
union is taken over all $m,j\in\N$ and $a_0,\dots,a_r\in M_n(\Q +
\i\Q)_{\rm sa}$. This is a countable union. Hence $N_0$ is again a  
$P$-null
set, and by Lemma \ref{lemma6-4}
\[
\spe\big(a_0\otimes\unit_n + \textstyle{\sum^r_{i=1}}a_i\otimes
y_i(\omega)\big)\subseteq\spe\big(a_0\otimes\unit_{\CB} +
\textstyle{\sum^r_{i=1}}a_i\otimes x_i\big)+
\big[\textstyle{-\frac1j,\frac1j}\big],
\]
for all $\omega\in\Omega_0\backslash N_0$, all $m,j\in\N$ and all
$a_0,\dots,a_r\in M_n(\Q + \i\Q)_{\rm sa}$. Taking intersection over
$j\in\N$ on the right-hand side, we get
\[
\spe\big(a_0\otimes\unit_n + \textstyle{\sum^r_{i=1}}a_i\otimes
y_i(\omega)\big)\subseteq\spe\big(a_0\otimes\unit_{\CB} +
\textstyle{\sum^r_{i=1}}a_i\otimes x_i\big),
\]
for $\omega\in\Omega_0\backslash N_0$, $m\in\N$ and $a_0,\dots,a_r\in
M_n(\Q+ \i\Q)_{\rm sa}$. Hence, by Theorem~\ref{thm1-4},
\[
\big\|p(y_1(\omega),\dots,y_r(\omega))\big\|\le\|p(x_1,\dots,x_r)\|,
\]
for all $p\in\C\< X_1,\dots,X_r \>$ and all
$\omega\in\Omega_0\backslash N_0$, which, by \eqref{eq6-10}, implies  
that
\[
\limsup_{n\to\infty}
\big\|p\big(X_1^{(n)}(\omega),\dots,X_r^{(n)}(\omega)\big)\big\|
\le\|p(x_1,\dots,x_r)\|,
\]
for all $p\in\C \< X_1,\dots,X_r\>$ and all $\omega\in\Omega_0\backslash
N_0$. This proves Proposition \ref{prop6-3}, which, together with Lemma
\ref{lemma6-2}, proves Theorem \ref{thm6-1} in the case $r\in\N$. The
case $r=\infty$ follows from the case $r\in\N$, because
$\C\<(X_i)_{i=1}^\infty\>=\cup_{r=1}^\infty\C\<(X_i)_{i=1}^r\>$.
\phantom{endofline}\hfill\qed

\section{$\Ext(\Cred(F_r))$ is not a group}\label{ext-resultat}

We start this section by translating Theorem~\ref{thm6-1} into a
corresponding result, where the self-adjoint Gaussian random matrices  
are
replaced by random unitary matrices and the semi-circular system is
replaced by a free family of Haar-unitaries.

Define $C^1$-functions $\varphi\colon\R\to\R$ and $\psi\colon\R\to\C$ by
\begin{equation}
\varphi(t)=
\begin{cases}
-\pi, &\textrm{if} \quad t\le-2, \\
\int_0^t \sqrt{4-s^2}\6s, &\textrm{if} \quad -2<t<2, \\
\pi, &\textrm{if} \quad t\ge2.
\end{cases}
\label{e7.1}
\end{equation}
and
\begin{equation}
\psi(t)=\e^{\i\varphi(t)}, \qquad (t\in\R).
\label{e7.2}
\end{equation}
Let $\mu$ be the standard semi-circle distribution on $\R$:
\[
\6\mu(t)=\frac{1}{2\pi}\sqrt{4-t^2}\cdot 1_{[-2,2]}(t)\6t,
\]
and let $\varphi(\mu)$ denote the push-forward measure of $\mu$ by  
$\varphi$,
i.e., $\varphi(\mu)(B)=\mu(\varphi^{-1}(B))$ for any Borel subset $B$  
of $\R$.
Since $\varphi'(t)=\sqrt{4-t^2}\cdot1_{[-2,2]}(t)$ for all $t$ in $\R$,  
it
follows that $\varphi(\mu)$ is the uniform distribution on  
$[-\pi,\pi]$, and,
hence, $\psi(\mu)$ is the Haar measure on the unit circle $\TT$ in $\C$.

The following lemma is a simple application of Voiculescu's results in
\cite{V2}.

\begin{lemma}\label{C*-iso}
Let $r\in\N\cup\{\infty\}$ and let $(x_i)_{i=1}^r$ be a
semi\/{\rm -}\/circular system in a $C^*$\/{\rm -}\/probability space $(\CB,\tau)${\rm ,} where
$\tau$ is a faithful state on $\CB$. Let $\psi\colon\R\to\TT$ be the
function defined in \eqref{e7.2}{\rm ,} and then put
\[
u_i=\psi(x_i), \qquad (i=1,\dots,r).
\]
Then there is a \/{\rm (}\/surjective\/{\rm )}\/ $*$\/{\rm -}\/isomorphism $\Phi\colon\Cred(F_r)\to
C^*((u_i)_{i=1}^r)${\rm ,} such that
\[
\Phi\big(\lambda(g_i)\big)=u_i, \qquad (i=1,\dots,r),
\]
where $g_1,\dots,g_r$ are the generators of the free group $F_r${\rm ,} and
$\lambda\colon F_r\to\CB(\ell^2(F_r))$ is the left regular  
representation
of $F_r$ on $\ell^2(F_r)$.
\end{lemma}

\Proof   Recall that $\Cred(F_r)$ is, by definition, the $C^*$-algebra in
$\CB(\ell^2(F_r))$ generated by $\lambda(g_1),\dots,\lambda(g_r)$. Let  
$e$
denote the unit in $F_r$ and let $\delta_e\in\ell^2(F_r)$ denote the
indicator function for $\{e\}$. Recall then that the vector state
$\eta=\<\cdot\delta_e,\delta_e\>\colon\CB(\ell^2(F_r))\to\C$,  
corresponding
to $\delta_e$, is faithful on $\Cred(F_r)$.
We recall further from \cite{V2} that $\lambda(g_1),\dots,\lambda(g_r)$
are $*$-free operators with respect to  $\eta$, and that each
$\lambda(g_i)$ is a Haar unitary, i.e.,
\[
\eta(\lambda(g_i)^n)=
\begin{cases}
1, &\textrm{if} \quad n=0,\\
0, &\textrm{if} \quad n\in\Z\setminus\{0\}.
\end{cases}
\]
Now, since $(x_i)_{i=1}^r$ are free self-adjoint operators in
$(\CB,\tau)$, $(u_i)_{i=1}^r$ are $*$-free unitaries in
$(\CB,\tau)$, and since, as noted above, $\psi(\mu)$ is the Haar  
measure on
$\TT$, all the $u_i$'s are Haar unitaries as well. Thus, the
$*$-distribution of
$(\lambda(g_i))_{i=1}^r$ with respect to  $\eta$ (in the sense of
\cite{V2}) equals that of $(u_i)_{i=1}^r$ with respect to~$\tau$. Since $\eta$
and $\tau$ are both faithful, the existence of a $*$-isomorphism  
$\Phi$, with
the properties set out in the lemma, follows from \cite[Remark~1.8]{V2}.
\Endproof\vskip4pt 

Let $r\in\N\cup\{\infty\}$.
As in Theorem~\ref{thm6-1}, we consider next, for each $n$ in $\N$,
independent random matrices $(X_i^{(n)})_{i=1}^r$ in
$\SGRM(n,\frac{1}{n})$. We then define, for each $n$, random unitary
$n\times n$ matrices $(U_i^{(n)})_{i=1}^r$, by setting
\begin{equation}
U_i^{(n)}(\omega)=\psi(X_i^{(n)}(\omega)), \qquad (i=1,2,\dots,r),
\label{e7.3}
\end{equation}
where $\psi\colon\R\to\TT$ is the function defined in \eqref{e7.2}.
Consider further the (free) generators $(g_i)_{i=1}^r$ of $F_r$. Then,
by the universal property of a free group, there exists, for each $n$ in
$\N$ and each $\omega$ in $\Omega$, a unique group homomorphism:
\[
\pi_{n,\omega}\colon F_r\to\CU(n)=\CU(M_n(\C)),
\]
satisfying
\begin{equation}
\pi_{n,\omega}(g_i)=U_i^{(n)}(\omega), \qquad (i=1,2,\dots,r).
\label{e7.4}
\end{equation}

\begin{theorem}\label{unitaer version}
Let $r\in\N\cup\{\infty\}$ and let{\rm ,} for each $n$ in $\N${\rm ,}
$(U_i^{(n)})_{i=1}^r$ be the random unitaries given by
\eqref{e7.3}. Let further for each $n$ in $\N$ and each $\omega$ in
$\Omega${\rm ,} $\pi_{n,\omega}\colon F_r\to\CU(n)$ be the group homomorphism
given by \eqref{e7.4}.

Then there exists a $P$-null set $N\subseteq\Omega${\rm ,} such that for all
$\omega$ in $\Omega\setminus N$ and all functions $f\colon F_r\to\C$  
with
finite support{\rm ,} we have
\[
\lim_{n\to\infty}\Big\|\sum_{\gamma\in
   F_r}f(\gamma)\pi_{n,\omega}(\gamma)\Big\|
= \Big\|\sum_{\gamma\in F_r}f(\gamma)\lambda(\gamma)\Big\|,
\]
where{\rm ,} as above{\rm ,} $\lambda$ is the left regular representation of $F_r$  
on
$\ell^2(F_r)$.
\end{theorem}

\Proof   In the proof we shall need the following simple observation: If
$a_1,\dots,a_s$, $b_1,\dots,b_s$ are $2s$ operators on a Hilbert space
$\CK$, such that $\|a_i\|$,\break $\|b_i\|\le1$ for all $i$ in  
$\{1,2,\dots,s\}$,
then
\begin{equation}
\|a_1a_2\cdots a_s - b_1b_2\cdots b_s\|\le \sum_{i=1}^s\|a_i-b_i\|.
\label{e7.5}
\end{equation}
We shall need further that for any positive $\varepsilon$ there exists
a polynomial $q$ in one variable, such that
\begin{equation}
|q(t)|\le1, \qquad (t\in[-3,3]),
\label{e7.6}
\end{equation}
and
\begin{equation}
|\psi(t)-q(t)|\le\varepsilon, \qquad (t\in[-3,3]).
\label{e7.7}
\end{equation}
Indeed, by Weierstrass' approximation theorem we may choose a
polynomial $q_0$ in one variable, such that
\begin{equation}
|\psi(t)-q_0(t)|\le\varepsilon/2, \qquad (t\in[-3,3]).
\label{e7.7a}
\end{equation}
Then put $q=(1+\varepsilon/2)^{-1}q_0$ and note that since $|\psi(t)|=1$
for all $t$ in $\R$, it follows from \eqref{e7.7a} that \eqref{e7.6}
holds. Furthermore,
\[
|q_0(t)- 
q(t)|\le\textstyle{\frac{\varepsilon}{2}}|q(t)|\le\frac{\varepsilon}{2},
  \qquad (t\in[-3,3]),
\]
which, combined with \eqref{e7.7a}, shows that \eqref{e7.7} holds.

After these preparations, we
start by proving the theorem in the case $r\in\N$.
For each $n$ in $\N$, let $X_1^{(n)},\dots,X_r^{(n)}$ be
independent random matrices in $\SGRM(n,\frac{1}{n})$ defined on
$(\Omega,\CF,P)$, and define the random unitaries
$U_1^{(n)},\dots,U_r^{(n)}$ as in \eqref{e7.3}. Then let $N$ be a
$P$-null set as in the main theorem (Theorem~\ref{thm6-1}). By
considering, for each $i$ in $\{1,2,\dots,r\}$, the polynomial
$p(X_1,\dots,X_r)=X_i$, it follows then from the main theorem that
\[
\lim_{n\to\infty}\big\|X_i^{(n)}(\omega)\big\|=2,
\]
for all $\omega$ in $\Omega\setminus N$. In particular, for each  
$\omega$ in
$\Omega\setminus N$, there exists an $n_{\omega}$ in $\N$, such that
\[
\big\|X_i^{(n)}(\omega)\big\|\le 3, \qquad
\mbox{whenever $n\ge n_{\omega}$ and $i\in\{1,2,\dots,r\}$.}
\]

Considering then the polynomial $q$ introduced above, it follows
from \eqref{e7.6} and \eqref{e7.7} that for all $\omega$ in  
$\Omega\setminus N$, we
have
\begin{equation}
\big\|q\big(X_i^{(n)}(\omega)\big)\big\|\le1, \qquad \mbox{whenever
   $n\ge n_{\omega}$ and $i\in\{1,2,\dots,r\}$,}
\label{e7.8}
\end{equation}
and
\begin{equation}
\big\|U_i^{(n)}(\omega)- 
q\big(X_i^{(n)}(\omega)\big)\big\|\le\varepsilon,
   \qquad \mbox{whenever $n\ge n_{\omega}$ and $i\in\{1,2,\dots,r\}$.}
\label{e7.9}
\end{equation}
Next, if $\gamma\in F_r\setminus\{e\}$, then $\gamma$ can be written
(unambiguesly) as a
reduced word: $\gamma=\gamma_1\gamma_2\cdots\gamma_s$, where
$\gamma_j\in\{g_1,g_2,\dots,g_r,g_1^{-1},g_2^{-1},\dots,g_r^{-1}\}$
for each $j$ in $\{1,2,\dots,s\}$, and where $s=|\gamma|$ is the
length of the reduced word for $\gamma$. It follows then, by
\eqref{e7.4}, that $\pi_{n,\omega}(\gamma)=a_1a_2\cdots a_s$, where
\begin{eqnarray*}
a_j&=&\pi_{n,\omega}(\gamma_j)\\
&\in&\big\{U_1^{(n)}(\omega),\dots,U_r^{(n)}( 
\omega),
U_1^{(n)}(\omega)^*,\dots,U_r^{(n)}(\omega)^*\big\}, \qquad  
(j=1,2,\dots,s).
\end{eqnarray*}

Combining now \eqref{e7.5}, \eqref{e7.8} and \eqref{e7.9}, it follows
that for any $\gamma$ in $F_r\setminus\{e\}$, there exists a
polynomial $p_{\gamma}$ in $\C\<X_1,\dots,X_r\>$, such that
\begin{multline}
\big\|\pi_{n,\omega}(\gamma)-
p_{\gamma}\big(X_1^{(n)}(\omega),\dots,X_r^{(n)}(\omega)\big)\big\|\le
|\gamma|\varepsilon,\\
 \mbox{whenever $n\ge n_{\omega}$ and
   $\omega\in \Omega\setminus N$}.
\label{e7.10}
\end{multline}
Now, let $\{x_1,\dots,x_r\}$ be a semi-circular system in a
$C^*$-probability space $(\CB,\tau)$, and put $u_i=\psi(x_i)$,
$i=1,2,\dots,r$. Then, by Lemma~\ref{C*-iso}, there is a surjective
$*$-isomorphism $\Phi\colon\Cred(F_r)\to C^*(u_1,\dots,u_r)$, such that
$(\Phi\circ\lambda)(g_i)=u_i$, $i=1,2,\dots,r$. Since $\|x_i\|\le 3$,
$i=1,2,\dots,r$, the arguments that lead to \eqref{e7.10} show also  
that
for any $\gamma$ in $F_r\setminus\{e\}$,
\begin{equation}
\big\|(\Phi\circ\lambda)(\gamma)-p_{\gamma}(x_1,\dots,x_r)\big\|
\le|\gamma|\varepsilon,
\label{e7.11}
\end{equation}
where $p_{\gamma}$ is the same polynomial as in \eqref{e7.10}. Note that
\eqref{e7.10} and \eqref{e7.11} also hold in the case $\gamma=e$, if we  
put
$p_e(X_1,\dots,X_r)=1$, and $|e|=0$.

Consider now an arbitrary function $f\colon F_r\to\C$ with finite  
support,
and then define the polynomial $p$ in $\C\<X_1,\dots,X_r\>$, by:
$p=\sum_{\gamma\in F_r}f(\gamma)p_{\gamma}$. Then, for any $\omega$ in
$\Omega\setminus N$ and any $n\ge n_{\omega}$, we have
\begin{equation}
\Big\|\sum_{\gamma\in F_r}f(\gamma)\pi_{n,\omega}(\gamma)
-p\big(X_1^{(n)}(\omega),\dots,X_r^{(n)}(\omega)\big)\Big\|
\le \Big(\sum_{\gamma\in F_r}|f(\gamma)|\cdot | \gamma |  
\Big)\varepsilon,
\label{e7.12}
\end{equation}
and
\begin{equation}
\Big\|\sum_{\gamma\in F_r}f(\gamma)\cdot(\Phi\circ\lambda)(\gamma)
-p(x_1,\dots,x_r)\Big\|
\le \Big(\sum_{\gamma\in F_r}|f(\gamma)|\cdot | \gamma |  
\Big)\varepsilon,
\label{e7.13}
\end{equation}
Taking also Theorem~\ref{thm6-1} into account, we may, on the basis of
\eqref{e7.12} and \eqref{e7.13}, conclude that for any $\omega$ in
$\Omega\setminus N$, we have
\[
\limsup_{n\to\infty}\Bigg|\Big\|\sum_{\gamma\in
   F_r}f(\gamma)\pi_{n,\omega}(\gamma)\Big\| - \Big\|\sum_{\gamma\in
   F_r}f(\gamma)\cdot(\Phi\circ\lambda)(\gamma)\Big\|\Bigg| \le
   2\varepsilon\Big(\sum_{\gamma\in F_r}|f(\gamma)|\cdot | \gamma |  
\Big).
\]
Since $\varepsilon>0$ is arbitrary, it follows that for any $\omega$ in
$\Omega\setminus N$,
\[
\lim_{n\to\infty}\Big\|\sum_{\gamma\in
   F_r}f(\gamma)\pi_{n,\omega}(\gamma)\Big\|=
\Big\|\sum_{\gamma\in F_r}f(\gamma)\cdot(\Phi\circ\lambda)(\gamma)\Big\|
=\Big\|\sum_{\gamma\in F_r}f(\gamma)\lambda(\gamma)\Big\|,
\]
where the last equation follows from the fact that $\Phi$ is a
$*$-isomorphism. This proves Theorem~\ref{unitaer version} in the case
where $r\in\N$. The case $r=\infty$ follows by trivial modifications
of the above argument.
\hfill\qed

\begin{remark}\label{aabent spoergsmaal}
The distributions of the random unitaries
$U_1^{(n)},\dots,U_r^{(n)}$ in Theorem~\ref{unitaer version} are quite
complicated. For instance, it is easily seen that for all $n$ in $\N$,
\[
P\big(\big\{\omega\in\Omega\bigm |  
U_1^{(n)}(\omega)=-\unit_n\big\}\big)>0.
\]
It would be interesting to know whether Theorem~\ref{unitaer version}  
also
holds, if, for each $n$ in~$\N$, $U_1^{(n)},\dots,U_r^{(n)}$ are  
replaced
be stochastically independent random unitaries
$V_1^{(n)},\dots,V_r^{(n)}$, which are all distributed according to the
normalized Haar measure on $\CU(n)$.
\end{remark}

\begin{corollary}\label{MF-algebra}
For any $r$ in $\N\cup\{\infty\}${\rm ,} the $C^*$\/{\rm -}\/algebra
$\Cred(F_r)$ has a unital embedding into the quotient $C^*$\/{\rm -}\/algebra
\[
\CC = \prod_n M_n(\C)\ \Big/ \sum_n M_n(\C),
\]
introduced in Section~{\rm \ref{sec6}.} In particular{\rm ,} $\Cred(F_r)$ is an
MF-algebra in the sense of Blackadar and Kirchberg \/{\rm (}\/cf.\ {\rm \cite{BK})}.
\end{corollary}

\Proof   This follows immediately from
   Theorem~\ref{unitaer
   version} and formula \eqref{eq6-8}. In fact, one only needs the  
existence
   of one $\omega$ in
$\Omega$ for which the convergence in Theorem~\ref{unitaer version}  
holds!
\Endproof\vskip4pt 

We remark that Corollary~\ref{MF-algebra} could also have been proved
directly from the main theorem (Theorem~\ref{thm6-1}) together with
Lemma~\ref{C*-iso}.

\begin{corollary}\label{ext is not a group}\quad
For any $r$ in $\{2,3,\ldots\}\cup\{\infty\}${\rm ,} the semi\/{\rm -}\/group\break
$\Ext(\Cred(F_r))$ is {\rm not} a group.
\end{corollary}

\Proof\quad   In Section 5.14 of Voiculescu's paper \cite{V5}, it is proved  
that\break
$\Ext(\Cred(F_r))$ cannot be a group, if there exists a sequence
$(\pi_n)_{n\in\N}$ of unitary representations $\pi_n\colon
F_r\to \CU(n)$, with the property that
\begin{equation}
\lim_{n\to\infty}\Big\|\sum_{\gamma\in F_r}f(\gamma)\pi_n(\gamma)\Big\|
=\Big\|\sum_{\gamma\in F_r}f(\gamma)\lambda(\gamma)\Big\|,
\label{e7.14}
\end{equation}
for any function $f\colon F_r\to\C$ with finite support.

For any $r\in\{2,3,\ldots\}\cup\{\infty\}$, the existence of such a  
sequence
$(\pi_n)_{n\in\N}$ follows
immediately from Theorem~\ref{unitaer version}, by considering one  
single
$\omega$ from the sure event $\Omega\setminus N$ appearing in that  
theorem.
\hfill\qed

\begin{remark}\label{voiculescus argument}
Let us briefly outline Voiculescu's  argument in \cite{V5} for the fact  
that
\eqref{e7.14} implies Corollary~\ref{ext is not a group}. It is
obtained by combining the following
two results of Rosenberg \cite{Ro} and Voiculescu \cite{V4},  
respectively:

\begin{enumerate}

\item If $\Gamma$ is a discrete countable nonamenable group, then
$\Cred(\Gamma)$ is not quasi-diagonal (\cite{Ro}).

\item A separable unital $C^*$-algebra $\CA$ is quasi-diagonal if and
only if there exists a sequence of natural numbers $(n_k)_{k\in\N}$
and a sequence $(\varphi_k)_{k\in\N}$\break of completely positive unital
maps $\varphi_k\colon\CA\to M_{n_k}(\C)$,
such that\break $\lim_{k\to\infty}\|\varphi_k(a)\|=\|a\|$
and $\lim_{k\to\infty}\|\varphi_k(ab)-\varphi_k(a)\varphi_k(b)\|=0$ for  
all
$a,b\in\CA$ (\cite{V4}).

\end{enumerate}

Let $\CA$ be a separable unital $C^*$-algebra. Then, as mentioned in
the introduction, $\Ext(\CA)$ is the set of equivalence classes
$[\pi]$ of one-to-one unital
$*$-homomorphisms $\pi$ of $\CA$ into the Calkin algebra
$\CC(\CH)=\CB(\CH)/\CK(\CH)$ over a separable infinite dimensional
Hilbert space $\CH$. Two
such $*$-homomorphisms are equivalent if they are equal up to a unitary
transformation of $\CH$. $\Ext(\CA)$ has a natural semi-group structure  
and
$[\pi]$ is invertible in $\Ext(\CA)$ if and only if $\pi$ has a unital
completely positive lifting: $\psi\colon\CA\to \CB(\CH)$
(cf.~\cite{Arv}). Let now $\CA=\Cred(F_r)$, where $r\in\{2,3,\dots\}\cup
\{\infty\}$. Moreover, let $\pi_n\colon F_r\to\CU_n$, $n\in\N$, be a  
sequence of
unitary representations satisfying \eqref{e7.14} and let $\CH$ be the  
Hilbert
space $\CH=\bigoplus^\infty_{n=1} \C^n$. Clearly, $\prod_n
M_n(\C)/\sum_n M_n(\C)$ embeds naturally into the Calkin algebra
$\CC(\CH)=\CB(\CH)/\CK(\CH)$. Hence, there
exists a one-to-one $*$-homomorphism $\pi\colon\CA\to \CC(\CH)$, such  
that
\[
\pi(\lambda(h)) = \rho\begin{pmatrix} \pi_1(h) & \ & 0 \\ \ & \pi_2(h)
   & \ \\ 0 & \ & \ddots \end{pmatrix},
\]
for all $h\in F_r$ (here $\rho$ denotes the quotient map from $\CB(\CH)$
to $\CC(\CH)$). Assume
$[\pi]$ is invertible in $\Ext(\CA)$. Then $\pi$ has a unital completely
positive lifting $\varphi\colon\CA\to \CB(\CH)$. Put
$\varphi_n(a)=p_n\varphi(a)p_n$, $a\in\CA$, where $p_n\in \CB(\CH)$ is  
the
orthogonal projection onto the component $\C^n$ of $\CH$. Then each
$\varphi_n$ is a unital completely positive map from $\CA$ to
$M_n(\C)$, and it is easy to check that
\[
\lim_{n\to\infty} \|\varphi_n(\lambda(h))-\pi_n(h)\|=0,\qquad (h\in  
F_r).
\]
 
From this it follows that
\[
\lim_{n\to\infty}\|\varphi_n(a)\|=\|a\| \quad \mbox{and} \quad
\lim_{n\to\infty} \|\varphi_n(ab)-\varphi_n(a)\varphi_n(b)\|=0,\qquad
(a,b\in \CA)
\]
so by (ii), $\CA=\Cred(F_r)$ is quasi-diagonal. But since $F_r$ is not
amenable for $r\ge2$, this contradicts (i). Hence $[\pi]$ is not  
invertible in
$\Ext(\CA)$.
\end{remark}

\begin{remark}
\label{rem7-7}
let $\CA$ be a separable unital $C^*$-algebra and let $\pi\colon\CA\to
\CC(\CH)=\CB(\CH)/\CK(\CH)$ be a one-to-one *-homomorphism. Then $\pi$  
gives rise to
an extension of $\CA$ by the compact operators $\CK=\CK(\CH)$, i.e., a
$C^*$-algebra $\CB$ together with a short exact sequence of
*-homomorphisms
\[
0\to \CK\stackrel{\iota}{\rightarrow} \CB
\stackrel{q}{\rightarrow}\CA\to0.
\]
Specifically, with $\rho\colon\CB(\CH)\to\CC(\CH)$ the quotient map,
$\CB=\rho^{-1}(\pi(\CA))$, $\iota$ is the inclusion map of $\CK$
into $\CB$ and $q=\pi^{-1}\circ\rho$.
Let now $\CA=\Cred(F_r)$, let $\pi\colon\CA\to \CC(\CH)$
be the one-to-one unital *-homomorphism from Remark~\ref{voiculescus
   argument}, and let $\CB$ be
the compact extension of $\CA$ constructed above. We then have
\begin{itemize}
\item[a)] $\CA=\Cred(F_r)$ is an exact $C^*$-algebra, but the compact
extension $\CB$ of $\CA$ is not exact.
\item[b)] $\CA=\Cred(F_r)$ is not quasi-diagonal but the compact
extension $\CB$ of $\CA$ is quasi-diagonal.
\end{itemize}
To prove a), note that $\Cred(F_r)$ is exact by \cite[Cor.~3.12]{DH} or
\cite[p.~453, l.~1--3]{Ki2}. Assume $\CB$ is also exact. Then, in
particular, $\CB$ is locally reflexive (cf.~\cite{Ki2}). Hence by the
lifting theorem in \cite{EH} and the nuclearity of $\CK$, the identity
map $\CA\to\CA$ has a unital completely positive lifting
$\varphi\colon\CA\to\CB$. If we consider $\varphi$ as a map from $\CA$  
to
$\CB(\CH)$, it is a unital completely positive lifting of  
$\pi\colon\CA\to
\CC(\CH)$, which contradicts that $[\pi]$ is not invertible in
$\Ext(\CA)$. To prove b), note that by Rosenberg's result, quoted in
(i) above, $\Cred(F_r)$ is not quasi-diagonal.
On the other hand, by the definition of $\pi$ in
Remark~\ref{voiculescus argument}, every $x\in\CB$ is
a compact perturbation of an operator of the form
\[
y= \begin{pmatrix} y_1 && 0\\ & y_2 & \\ 0 && \ddots\end{pmatrix},
\]
where $y_n\in M_n(\C)$, $n\in\N$. Hence $\CB$ is quasi-diagonal.
\end{remark}

\section{Other applications}\label{sec8}

Recall that a $C^*$-algebra $\CA$ is called exact if, for every pair  
$(\CB,\CJ)$
consisting of a $C^*$-algebra $\CB$ and closed two-sided ideal $\CJ$
in $\CB$, the sequence
\begin{equation}
\label{eq8-1}
0\to \CA\mathop{\otimes}_{\min} \CJ\to \CA\mathop{\otimes}_{\min}\CB \to
\CA\mathop{\otimes}_{\min}(\CB/\CJ)\to0
\end{equation}
is exact (cf.\ \cite{Ki1}). In generalization of the construction  
described in
the\break paragraph preceding Lemma~\ref{lemma6-4}, we may, for any sequence
$(\CA_n)^\infty_{n=1}$ of\break $C^*$-algebras, define two $C^*$-algebras
\begin{eqnarray*}
\prod_n \CA_n &=& \big\{(a_n)^\infty_{n=1} \mid a_n\in \CA_n,\
\textstyle{\sup_{n\in\N}}\|a_n\| < \infty\big\}\\
\sum_n \CA_n &=& \big\{(a_n)^\infty_{n=1} \mid a_n\in \CA_n,\
\textstyle{\lim_{n\to\infty}}\|a_n\|=0\big\}.
\end{eqnarray*}
The latter $C^*$-algebra is a closed two-sided ideal in the first, and
the norm in the quotient $C^*$-algebra $\prod_n\CA_n/\sum_n\CA_n$ is  
given
by
\begin{equation}
\label{eq8-2}
\big\|\rho\big((x_n)^\infty_{n=1}\big)\big\| = \limsup_{n\to\infty}
\|x_n\|,
\end{equation}
where $\rho$ is the quotient map (cf.\ \cite[Lemma 6.13]{RLL}) . In the
following we let $\CA$ denote an exact $C^*$-algebra. By \eqref{eq8-1}  
we
have the following natural identification of $C^*$-algebras
\[
\CA \mathop{\otimes}_{\min} \Big(\prod_n M_n(\C)\Big/  
\sum_nM_n(\C)\Big) =
\Big(\CA\mathop{\otimes}_{\min}\prod_n M_n(\C)\Big)\Big/
\Big(\CA\mathop{\otimes}_{\min}
\sum_nM_n(\C)\Big).
\]
Moreover, we have (without assuming exactness) the following natural
identification
\[
\CA\mathop{\otimes}_{\min}\sum_n M_n(\C)=\sum_n M_n(\CA)
\]
and the natural inclusion
\[
\CA\mathop{\otimes}_{\min}\prod_n M_n(\C)\subseteq\prod_n M_n(\CA).
\]
If $\dim(\CA)<\infty$, the inclusion becomes an identity, but in  
general the
inclusion is proper. Altogether we have for all exact $C^*$-algebras  
$\CA$
a natural inclusion
\begin{equation}
\label{eq8-3}
\CA\mathop{\otimes}_{\min}\Big(\prod_n M_n(\C)\Big/
\sum_nM_n(\C)\Big)\subseteq \prod_nM_n(\CA)\Big/ \sum_n M_n(\CA).
\end{equation}
Similarly, if $n_1<n_2<n_3<\cdots$, are natural numbers, then
\begin{equation}
\label{eq8-4}
\CA\mathop{\otimes}_{\min}\Big(\prod_kM_{n_k}(\C)\Big/\sum_k  
M_{n_k}(\C)\Big)\subseteq
\prod_kM_{n_k}(\CA)\Big/\sum_kM_{n_k}(\CA).
\end{equation}
After these preparations we can now prove the following generalizations
of Theorems~\ref{thm6-1} and \ref{unitaer version}.

\begin{theorem}
\label{thm8-1}
Let $(\Omega,\CF,P)${\rm ,} $N${\rm ,}
$(X_i^{(n)})_{i=1}^r$ and $(x_i)_{i=1}^r$ be as in
Theorem~{\rm \ref{thm6-1},}
and let $\CA$ be a unital exact $C^*$\/{\rm -}\/algebra. Then for all
polynomials $p$ in $r$ noncommuting variables and with coefficients in
$\CA$ \/{\rm (}\/i.e.{\rm ,} $p$ is in the algebraic tensor product
$\CA\otimes\C\<(X_i)_{i=1}^r\>${\rm ),} and all $\omega\in\Omega\backslash N${\rm ,}
\begin{equation}
\label{eq8-5}
\lim_{n\to\infty}\big\|p\big((X_i^{(n)}(w))_{i=1}^r\big)\big\|_{M_n(\CA) 
}
=\big\|p\big((x_i)_{i=1}^r\big)\big\|
_{\CA\otimes_{\min}C^*((x_i)_{i=1}^r,\unit_{\CB})}.
\end{equation}
\end{theorem}
\vskip8pt

\Proof   We consider only the case $r\in\N$. The case $r=\infty$ is
proved similarly.
By Theorem~\ref{thm6-1} we can for each $\omega\in\Omega\backslash N$  
define a
unital embedding $\pi_\omega$ of $C^*(x_1,\dots,x_r,\unit_{\CB})$ into
$\prod_n M_n(\C)/\sum_n M_n(\C)$, such that
\[
\pi_\omega(x_i) =  
\rho\big(\big(X_i^{(n)}(\omega)\big)^\infty_{n=1}\big),\quad
i=1,\dots,r,
\]
where $\rho\colon\prod_n M_n(\C)\to\prod_n M_n(\C)/\sum_n M_n(\C)$ is  
the
quotient map. Since $\CA$ is exact, we can, by \eqref{eq8-3}, consider
$\id_\CA\otimes\pi_\omega$ as a unital embedding of
$\CA\mathop{\otimes}_{\min}C^*(x_1,\dots,x_r,\unit_{\CB})$ into  
$\prod_n
M_n(\CA)/\sum_n M_n(\CA)$, for which
\[
(\id_\CA\otimes\pi_\omega)(a\otimes x_i)=\tilde{\rho}\big(\big(a\otimes
X_i^{(n)}(\omega)\big)^\infty_{n=1}\big), \quad i=1,\dots,r,
\]
where $\tilde{\rho}\colon\prod_n M_n(\CA)\to\prod_n M_n(\CA)/\sum  
M_n(\CA)$
is the quotient map. Hence, for every $p$ in  
$\CA\otimes\C\<X_1,\dots,X_r\>$,
\[
(\id_\CA\otimes\pi_\omega)\big(p(x_1,\dots,x_r)\big)
=\tilde{\rho}\big(\big(
p(X_1^{(n)}(\omega),\dots,X_r^{(n)}(\omega))\big)^\infty_{n=1}\big).
\]

By \eqref{eq8-2} it follows that for all $\omega\in\Omega/N$, and
every $p$ in $\CA\otimes\C\<X_1,\dots,X_r\>$,
\[
\big\|p(x_1,\dots,x_r)\big\|
_{\CA\otimes_{\min}C^*(x_1,\dots,x_r,\unit_{\CB})} =
\limsup_{n\to\infty}\big\|p\big(X_1^{(n)}(\omega),
\dots,X_r^{(n)}(\omega)\big)\big\|_{M_n(\CA)}.
\]

Consider now a fixed $\omega\in\Omega\backslash N$. Put
\[
\alpha =
\liminf_{n\to\infty}\big\|p\big(X_1^{(n)}(\omega),
\dots,X_r^{(n)}(\omega)\big)\big\|_{M_n(\CA)},
\]
and choose natural numbers $n_1<n_2<n_3<\cdots$, such that
\[
\alpha = \lim_{k\to\infty}\big\|p\big(X_1^{(n_k)}(\omega),
\dots,X_r^{(n_k)}(\omega)\big)\big\|_{M_n(\CA)}.
\]
By Theorem~\ref{thm6-1} there is a unital embedding $\pi'_\omega$ of
$C^*(x_1,\dots,x_r,\unit_{\CB})$ into the quotient $\prod_k
M_{n_k}(\C)/\sum_k M_{n_k}(\C)$, such that
\[
\pi'_\omega(x_i)=\rho'\big(\big(X_i^{(n_k)}(\omega)\big)^\infty_{k=1}\big),
\quad i=1,\dots,r,
\]
where $\rho'\colon\prod_k M_{n_k}(\C)\to \prod_k M_{n_k}(\C)/\sum_k
M_{n_k}(\C)$ is the quotient map. Using \eqref{eq8-4} instead of
\eqref{eq8-3}, we get, as above, that
\begin{eqnarray*}
\|p(x_1,\dots,x_r)\|_{\CA\otimes_{\min}C^*(x_1,\dots,x_r,\unit_{\CB})}
&=& \limsup_{k\to\infty}
\big\|p\big(X_1^{(n_k)}(\omega),
\dots,X_r^{(n_k)}(\omega)\big)\big\|_{M_n(\CA)}\\[.2cm]
&=& \alpha\\[.2cm]
&=& \liminf_{n\to\infty}
\big\|p\big(X_1^{(n)}(\omega),
\dots,X_r^{(n)}(\omega)\big)\big\|_{M_n(\CA)}.
\end{eqnarray*}
This completes the proof of \eqref{eq8-5}.
\hfill\qed 

\begin{theorem}
\label{thm8-2}
Let $(\Omega,\CF,P)${\rm ,} $(U_i^{(n)})_{i=1}^r${\rm ,}
$\pi_{n,\omega},\lambda$ and $N$ be as in
Theorem~{\rm \ref{unitaer version}.}
Then for every unital exact $C^*$\/{\rm -}\/algebra $\CA${\rm ,} every function
$f\colon F_r\to \CA$ with finite support \/{\rm (}\/i.e.\ $f$ is in the algebraic
tensor product $\CA\otimes\C F_r$\/{\rm ),}\/ and for every
$\omega\in\Omega\backslash N$
\[
\lim_{n\to\infty}\Big\|\sum_{\gamma\in F_r}f(\gamma)\otimes
\pi_{n,\omega}(\gamma)\Big\|_{M_n(\CA)}
=\Big\|\sum_{\gamma\in F_r}f(\gamma)\otimes\lambda(\gamma)\Big\|
_{\CA\otimes_{\min}\Cred(F_r)}.
\]
\end{theorem}
\vskip8pt

\Proof  
This follows from Theorem~\ref{unitaer version} in the same way as  
Theorem
\ref{thm8-1} follows from Theorem~\ref{thm6-1}, so we leave the details  
of
the proof to the reader.
\Endproof\vskip4pt 

In Corollary \ref{cor8-3} below we use Theorem \ref{thm8-1} to give new
proofs of two key results from our previous paper \cite{HT2}. As in
\cite{HT2} we denote by $\GRM(n,n,\sigma^2)$ or $\GRM(n,\sigma^2)$ the
class of $n\times n$ random matrices $Y=(y_{ij})_{1\le i,j\le n}$, whose
entries $y_{ij}$, $1\le i,j\le n$, are $n^2$ independent and identically distributed  complex Gaussian
random variables
with density $(\pi\sigma^2)^{-1}\exp(-|z|^2/\sigma^2)$, $z\in\C$. It is
elementary to check that $Y$ is in $\GRM(n,\sigma^2)$, if and only if
$Y=\frac{1}{\sqrt{2}}(X_1+\i X_2)$, where
\[
X_1 = \frac{1}{\sqrt{2}}(Y+Y^*),\quad X_2 = \frac{1}{\i\sqrt{2}}(Y-Y^*)
\]
are two stochastically independent self-adjoint random matrices from the
class $\SGRM(n,\sigma^2)$.

\begin{corollary}[{[}HT2, Thm.~4.5~and~Thm.~8.7{]}]\label{cor8-3}
Let $\CH,\CK$ be Hilbert\break spaces{\rm ,} let $c>0${\rm ,} let $r\in\N$ and let
$a_1,\dots,a_r\in\CB(\CH,\CK)$ such that
\[
\Big\|\sum_{i=1}^r a^*_ia_i\Big\|\le c\quad\mbox{and}\quad
\Big\|\sum_{i=1}^r a_ia^*_i\Big\|\le 1,
\]
and such that $\{a_i^*a_j\mid i,j=1,\dots,r\}\cup\{\unit_{\CB(\CH)}\}$
generates an
exact $C^*$\/{\rm -}\/algebra $\CA\subseteq \CB(\CH)$. Assume further that
$Y_1^{(n)},\dots,Y_r^{(n)}$ are stochastically independent random
matrices from the class $\GRM(n,\frac 1n)${\rm ,} and put $S_n=\sum^r_{i=1}
a_i\otimes Y_i^{(n)}$. Then for almost all $\omega$ in the underlying
probability space $\Omega${\rm ,}
\begin{equation}
\label{eq8-6}
\limsup_{n\to\infty}\max\big\{\spe(S_n(\omega)^*S_n(\omega))\big\}\le
(\sqrt{c}+1)^2.
\end{equation}
If{\rm ,} furthermore{\rm ,} $c>1$ and $\sum^r_{i=1} a_i^*a_i=c\unit_{\CB(\CH)}${\rm ,}  
then
\begin{equation}
\label{eq8-7}
\liminf_{n\to\infty}\min\big\{\spe(S_n(\omega)^*S_n(\omega))\big\}
\ge(\sqrt{c}-1)^2.
\end{equation}
\end{corollary}
\vskip8pt

\Proof  
By the comments preceding Corollary \ref{cor8-3}, we can write
\[
Y_i^{(n)} = \frac{1}{\sqrt{2}}(X_{2i-1}^{(n)}+\i X_{2i}^{(n)}),\quad
i=1,\dots,r,
\]
where $X_1^{(n)},\dots,X_{2r}^{(n)}$ are independent self-adjoint  
random matrices
from\break $\SGRM(n,\frac 1n)$. Hence $S^*_nS_n$ is a second order polynomial
in $(X_1^{(n)},\dots,X_{2r}^{(n)})$ with coefficient in the exact unital
$C^*$-algebra $\CA$ generated by $\{a_i^*a_j\mid i,j=1,\dots,r\}\cup
\{\unit_{\CB(\CH)}\}$. Hence, by Theorem \ref{thm8-1}, there is a  
$P$-null
set $N$ in the underlying probability space $(\Omega,\CF,P)$ such that
\[
\lim_{n\to\infty}\|S^*_n(\omega)S_n(\omega)\| = \Big\|\Big(\sum^r_{i=1}
a_i\otimes y_i\Big)^*\Big(\sum^r_{i=1} a_i\otimes y_i\Big)\Big\|,
\]
where $y_i=\frac{1}{\sqrt{2}}(x_{2i-1}+\i x_{2i})$ and
$(x_1,\dots,x_{2r})$ is any semicircular system in a $C^*$-probability
space $(\CB,\tau)$ with $\tau$ faithful. Hence, in the terminology of
\cite{V2}, $(y_1,\dots,y_r)$ is a circular system with the normalization
$\tau(y^*_iy_i)=1$, $i=1,\dots,r$. By \cite{V2}, a concrete model for  
such
a circular system is
\[
y_i = \ell_{2i-1}+\ell^*_{2i},\quad i=1,\dots,r
\]
where $\ell_1,\dots,\ell_{2r}$ are the creation operators on the full
Fock space
\[
\CT=\CT(\CL)=\C\oplus \CL\oplus \CL^{\otimes 2}\oplus\cdots
\]
over a Hilbert space $\CL$ of dimension $2r$, and $\tau$ is the vector
state given by the unit vector $1\in\C\subseteq\CT(\CL)$. Moreover,  
$\tau$
is a faithful trace on the $C^*$-algebra
$\CB=C^*(y_1,\dots,y_{2r},\unit_{\CB(\CT(\CL))})$. The creation  
operators
$\ell_1,\dots,\ell_{2r}$ satisfy
\[
\ell_i^*\ell_j =
\begin{cases}
1, &\textrm{if} \ i=j, \\
0, &\textrm{if} \ i\ne j.
\end{cases}
\]
Hence, we get
\[
\sum^r_{i=1} a_i\otimes y_i = \Big(\sum^r_{i=1}
a_i\otimes\ell_{2i-1}\Big)+\Big(\sum^r_{i=1}
a_i\otimes\ell_{2i}^*\Big)=z+w,
\]
where
\[
z^*z = \Big(\sum^r_{i=1} a_i^*a_i \Big)\otimes \unit_{\CB(\CT)}
\quad\mbox{and}\quad
ww^*=\Big(\sum_{i=1}^ra_ia^*_i\Big)\otimes\unit_{\CB(\CT)}.
\]
Thus,
\[
\Big\| \sum^r_{i=1} a_i\otimes y_i\Big\| \le \|z\|+\|w\|
\le  \Big\|\sum^r_{i=1} a_i^*a_i\Big\|^\frac12 + \Big\|\sum^r_{i=1}
a_ia_i^*\Big\|^\frac12
\le  \sqrt{c}+1.
\]
This proves \eqref{eq8-5}. If, furthermore, $c>1$ and $\sum^r_{i=1}
a^*_ia_i=c\cdot\unit_{\CB(\CH)}$, then  
$z^*z=c\unit_{\CA\otimes\CB(\CT)}$
and, as before, $\|w\|\le1$. Thus,
for all $\xi\in \CH\otimes\CT$, $\|z\xi\|=\sqrt{c}\|\xi\|$ and
$\|w\xi\|\le\|\xi\|$. Hence
\[
(\sqrt{c}-1)\|\xi\|\le\|(z+w)\xi\|\le(\sqrt{c}+1)\|\xi\|,\quad (\xi\in
\CH\otimes\CT),
\]
which is equivalent to
\[
(\sqrt{c}-1)^2\unit_{\CB(\CH\otimes\CT)}\le(z+w)^*(z+w)\le(\sqrt{c}+1)^2
\unit_{\CB(\CH\otimes\CT)},
\]
\[
-2\sqrt{c} \unit_{\CB(\CH\otimes\CT)}\le (z+w)^*(z+w)-(c+1)
\unit_{\CB(\CH\otimes\CT)}\le 2\sqrt{c} \unit_{\CB(\CH\otimes\CT)},
\]
and therefore
\begin{equation}
\label{eq8-8}
\big\|(z+w)^*(z+w)-(c+1)\unit_{\CB(\CH\otimes\CT)}\big\|\le 2\sqrt{c}.
\end{equation}
Since $S_n^*S_n$ is a second order polynomial in
$(X_1^{(n)},\dots,X_{2r}^{(n)})$ with coefficients in~$\CA$, the same
holds for $S_n^*S_n-(c+1)\unit_{M_n(\CA)}$. Hence, by Theorem  
\ref{thm8-1} and~\eqref{eq8-8},
\begin{eqnarray*}
&&\hskip-.75in \lim_{n\to\infty}\big\|S_n(\omega)^*  
S_n(\omega)-(c+1)\unit_{M_n(\CA)}\big\| \\[.2cm]
&&= \Big\|\Big(\sum_{i=1}^ra_i\otimes  
y_i\Big)^*\Big(\sum_{i=1}^ra_i\otimes
y_i\Big)-(c+1)\unit_{\CB(\CH\otimes\CT)}\Big\|\\[.2cm]
 &&\le 2\sqrt{c}.
\end{eqnarray*}
\newpage
\noindent
Therefore, $\liminf_{n\to\infty}
\min\{\spe(S_n(\omega)^*S_n(\omega))\}\ge(c+1)-2\sqrt{c}$, which proves
\eqref{eq8-7}.
\hfill\qed

\begin{remark}
\label{rem8-4}
The condition that $\{a_i^*a_j\mid  
i,j=1,\dots,r\}\cup\{\unit_{\CB(\CH)}\}$
generates
an exact $C^*$-algebra is essential for Corollary \ref{cor8-3} and hence
also for Theorem \ref{thm8-1}. Both \eqref{eq8-6} and \eqref{eq8-7} are
false in the general nonexact case (cf.\ \cite[Prop.~4.9]{HT2} and
\cite{HT3}).
\end{remark}

We turn next to a result about the constants $C(r)$, $r\in\N$,  
introduced by
Junge and Pisier in connection with their proof of
\begin{equation}
\label{eq8-9}
\CB(\CH)\mathop{\otimes}_{\max}\CB(\CH)\ne
\CB(\CH)\mathop{\otimes}_{\min} \CB(\CH).
\end{equation}

\begin{definition}[\cite{JP}]\label{def8-5}
For $r\in\N$, let $C(r)$ denote the infimum of all $C\in\R_+$ for which
there exists a sequence of natural numbers $(n(m))_{m=1}^\infty$ and a
sequence of $r$-tuples of $n(m)\times n(m)$ unitary matrices
\[
(u_1^{(m)},\dots,u_r^{(m)})^\infty_{m=1}
\]
such that for all $m,m'\in\N$, $m\ne m'$
\begin{equation}
\label{eq8-10}
\big\|\sum^r_{i=1} u_i^{(m)}\otimes\bar{u}_i^{(m')}\big\| \le C,
\end{equation}
where $\bar{u}_i^{(m')}$ is the unitary matrix obtained by complex
conjugation of the entries of $u_i^{(m')}$.
\end{definition}

To obtain \eqref{eq8-9}, Junge and Pisier proved that $\lim_{r\to\infty}
\frac{C(r)}{r}=0$. Subsequently, Pisier \cite{P3} proved that $C(r) \ge
2\sqrt{r-1}$ for all $r\ge 2$. Moreover, using Ramanujan graphs, Valette
\cite{V} proved that $C(r)=2\sqrt{r-1}$, when $r=p+1$ for an odd prime
number $p$. From Theorem \ref{thm8-2} we obtain

\begin{corollary}
\label{cor8-6}
$C(r)=2\sqrt{r-1}$ for all $r\in\N${\rm ,} $r\ge 2$.
\end{corollary}

\Proof  
Let $r\ge 2$, and let $g_1,\dots,g_r$ be the free generators of $F_r$
and let $\lambda$ denote the left regular representation of $F_r$ on
$\ell^2(F_r)$. Recall from \cite[Formulas (4) and (7)]{P3} that
\begin{equation}
\label{eq8-11}
\Big\|\sum^r_{i=1} \lambda(g_i)\otimes v_i\Big\| = 2\sqrt{r-1}
\end{equation}
for all unitaries $v_1,\dots,v_r$ on a Hilbert space $\CH$. Let $C>
2\sqrt{r-1}$. We will construct natural numbers $(n(m))^\infty_{m=1}$
and $r$-tuples of $n(m)\times n(m)$ unitary matrices
\[
(u_1^{(m)},\dots,u_r^{(m)})^\infty_{m=1}
\]
such that \eqref{eq8-10} holds for $m,m'\in\N$, $m\ne m'$. Note that by
symmetry it is sufficient to check \eqref{eq8-10} for $m'<m$. Put first
\[
n(1)=1\quad\mbox{and}\quad u_1^{(1)}=\cdots = u_r^{(1)}=1.
\]

Proceeding by induction,
let $M\in\N$ and assume that we have found $n(m)\in\N$ and $r$-tuples of
$n(m)\times m(n)$ unitaries $(u_1^{(m)},\dots,u_r^{(m)})$ for $2\le m\break\le
M$, such that \eqref{eq8-10} holds for $1\le m' < m\le M$. By
\eqref{eq8-11},
\[
\Big\|\sum^r_{i=1} \lambda(g_i)\otimes\bar{u}_i^{(m)}\Big\| =
2\sqrt{r-1},
\]
for $m=1,2,\dots,M$. Applying Theorem \ref{thm8-2} to the exact
$C^*$-algebras $\CA_{m'} = M_{n(m')}(\C)$, $m'=1,\dots,M$, we have
\[
\lim_{n\to\infty} \Big\|\sum^r_{i=1}
\pi_{n,\omega}(g_i)\otimes\bar{u}_i^{(m')}\Big\|=2\sqrt{r-1}<C, \quad
(m'=1,2,\dots,M),
\]
where $\pi_{n,\omega}\colon F_r\to\CU(n)$ are the group homomorphisms
given by \eqref{e7.4}.
Hence, we can choose $n\in\N$ such that
\[
\Big\|\sum^r_{i=1} \pi_{n,\omega}(g_i)\otimes\bar{u}_i^{(m')}\Big\|<
C,\quad m'=1,\dots,M.
\]
Put $n(M+1)=n$ and $u_i^{(M+1)}=\pi_{n,\omega}(g_i)$, $i=1,\dots,r$.  
Then
\eqref{eq8-10} is satisfied for all $m,m'$ for which $1\le m'<m\le
M+1$. Hence, by induction we get the desired sequence of numbers $n(m)$
and $r$-tuples of $n(m)\times n(m)$ unitary matrices.
\Endproof\vskip4pt 

We close this section with an application of Theorem~\ref{thm6-1} to  
powers
of random matrices:

\begin{corollary}
\label{cor8-7}
Let for each $n\in\N$ $Y_n$ be a random matrix in the class
$\GRM(n,\frac 1n)${\rm ,} i.e.{\rm ,} the entries of $Y_n$ are $n^2$
independent and identically distributed  complex Gaussian variables with density $\frac{n}{\pi}
\e^{-n|z|^2}${\rm ,} $z\in\C$. Then for all $p\in\N$
\[
\lim_{n\to\infty} \|Y_n(\omega)^p\| =
\bigg(\frac{(p+1)^{p+1}}{p^p}\bigg)^\frac12,
\]
for almost all $\omega$ in the underlying probability space $\Omega$.
\end{corollary}

\Proof   By the remarks preceding Corollary~\ref{cor8-3}, we have
\[
(Y_n)^p =
\bigg(\frac{1}{\sqrt{2}}\big(X_1^{(n)}+\i X_2^{(n)}\big)\bigg)^p,
\]
\newpage
\noindent
where, for each $n\in\N$, $X_1^{(n)},X_2^{(n)}$ are two independent
random matrices from $\SGRM(n,\frac 1n)$. Hence, by Theorem  
\ref{thm6-1},
we have for almost all $\omega\in\Omega$:
\[
\lim_{n\to\infty}\|Y_n(\omega)^p\|=\|y^p\|,
\]
where $y=\frac{1}{\sqrt{2}}(x_1+\i x_2)$, and $\{x_1,x_2\}$ is a  
semicircular
system in a $C^*$-prob\-ability space $(\CB,\tau)$ with $\tau$
faithful. Hence, $y$ is a circular element in $\CB$ with the standard
normalization $\tau(y^*y)=1$. By \cite[Prop.~4.1]{La}, we
therefore have $\|y^p\|=((p+1)^{p+1}/p^p)^\frac12$.
\hfill\qed

\begin{remark}
\label{rem8-7}
For $p=1$, Corollary~\ref{cor8-7} is just the complex version of\break Geman's
result \cite{Ge} for square matrices (see \cite[Thm.~2.16]{Ba} or
\cite[Thm.~7.1]{HT1}), but
for $p\ge 2$ the result is new. In \cite[Example 1, p.125]{We}, Wegmann
proved that the empirical eigenvalue distribution of $(Y^p_n)^*Y^p_n$
converges almost surely to a probability measure $\mu_p$ on $\R$ with
\[
\max(\supp(\mu_p)) = \frac{(p+1)^{p+1}}{p^p}.
\]
This implies that for all $\varepsilon>0$, the number of eigenvalues of
$(Y^p_n)^*Y^p_n$, which are larger than $(p+1)^{p+1}/p^p+\varepsilon$,
grows slower than $n$, as $n\to\infty$ (almost
surely). Corollary~\ref{cor8-7} shows that this number is, in fact,
eventually 0 as $n\to\infty$ (almost surely).
\end{remark}


\references {999999}

\bibitem[An]{An} \name{J.~Anderson}, A $C^*$-algebra $A$ for which
$\mathrm{Ext}(A$) is not a group, {\it Ann.\ of  Math}.\ {\bf 107}  
(1978), 455--458.
 
\bibitem[Ar]{Ar} \name{L.~Arnold}, On the asymptotic distribution of the
eigenvalues of random matrices, {\it J.\ Math.\ Anal.\ Appl}.\ {\bf 20}
(1967), 262--268.

\bibitem[Arv]{Arv} \name{W.~Arveson}, Notes on extensions of $C^*$-algebras,  
{\it
     Duke Math.\ J.}\ {\bf 44}  (1977), 329--355.

\bibitem[Ba]{Ba} \name{Z.\ D.~Bai}, Methodologies in spectral analysis of large
dimensional random matrices, A review, {\it Statistica Sinica} {\bf 9}   
(1999),
611--677.
 

\bibitem[BDF1]{BDF1} \name{L.\ G.~Brown,  R.\ G.~Douglas}, and P.\  
\name{A.~Fillmore}, Unitary
equivalence modulo the compact operators and extensions of
$C^*$-algebras, {\it Proc.~Conf.~on Operator Theory\/},  {\it Lecture  
Notes in
Math\/}.\ {\bf 34} (1973),  58--128,  {\it Springer-Verlag}, New York.

\bibitem[BDF2]{BDF2} \bibline,  
Extensions of
$C^*$-algebras and $K$-homology, {\it Ann.\ of Math}.\ {\bf 105}   
(1977), 265--324.


\bibitem[Be]{Be} \name{W.~Beckner}, A generalized Poincar\'e inequality for
Gaussian measures, {\it Proc.\ Amer.\ Math.\ Soc}.\ {\bf 105}  (1989),  
397--400.

\bibitem[BK]{BK} \name{B.~Blackadar} and  \name{E.~Kirchberg}, Generalized inductive
limits of finite dimensional $C^*$-algebras, {\it Math.\ Ann}.\ {\bf  
307}
(1997), 343--380.
 

\bibitem[Bre]{bre} \name{L.~Breiman}, {\it Probability\/}, {\it  Classics in  
Applied
     Mathematics} {\bf 7}, Society for Industrial and Applied  
Mathematics (SIAM), Philadelphia, PA (1992).

\bibitem[BY]{BY} \name{Z.\ D.~Bai} and \name{Y.\ Q.~Yin}, Neccesary and sufficient
conditions for almost sure convergence of the largest eigenvalue of a
Wigner matrix, {\it Ann.\ of Probab}.\ {\bf 16}  (1988), 1729--1741.

\bibitem[CE]{CE} \name{M.\ D.~Choi} and \name{E.~Effros}, The completely positive
lifting problem for $C^*$-algebras, {\it Ann.\ of Math}.\ {\bf 104}  
(1976),
585--609.

\bibitem[Cf]{Cf} \name{H.~Chernoff}, A note on an inequality involving the
normal distribution, {\it Ann.\ Probab}.\ {\bf 9}  (1981),
533--535.

\bibitem[Cn]{Cn} \name{L.\ H.\ Y.~Chen}, An inequality for the multivariate  
normal
distribution, {\it J.~Multivariate Anal}.\ {\bf 12}  (1982),
306--315.
 

\bibitem[Da]{Da} \name{K.H.~Davidson}, The distance between unitary orbits of
normal operators, {\it Acta Sci.~Math}.\ {\bf 50}  (1986), 213--223.

\bibitem[DH]{DH} \name{J.~De Canni\`ere} and \name{U.~Haagerup}, Multipliers of the  
Fourier
   algebras of some simple Lie groups and their discrete subgroups, {\it
   Amer.\ J.\ Math}.\ {\bf 107} (1984), 455--500.

 
\bibitem[EH]{EH} \name{E.~Effros} and \name{U.~Haagerup}, Lifting problems and local
reflexivity for $C^*$-algebras, {\it Duke Math.\ J}.\ {\bf 52}  (1985),  
103--128.

\bibitem[F]{F} \name{G.\ B.~Folland}, {\it Real Analysis\/}, {\it Modern  
Techniques and their
Applications\/}, {\it Pure and Applied Mathematics\/}, {\it John Wiley  
and Sons}, New York (1984).

 
\bibitem[Ge]{Ge} \name{S.~Geman}, A limit theorem for the norm of random  
matrices,
{\it Annals Prob}.\ {\bf 8} (1980) 252--261.
 
\bibitem[HP]{HP2} \name{F.~Hiai} and \name{D.~Petz}, Asymptotic freeness almost
everywhere for random matrices, {\it Acta Sci.\ Math.} (Szeged) {\bf 66}
(2000), 809--834.

\bibitem[HT1]{HT1} \name{U.~Haagerup} and \name{S.~Thorbj{\o}rnsen}, Random matrices  
with
complex Gaussian entries, {\it Expo.\ Math}.\ {\bf 21} (2003), 293-337.

\bibitem[HT2]{HT2} \bibline, Random matrices  
and
K-theory for exact $C^*$-algebras, {\it Documenta Math}.\ {\bf 4}  
(1999), 330--441.

\bibitem[HT3]{HT3} \bibline, Random matrices  
and
nonexact $C^*$-algebras, in {\it $C^*$-algebras} (J.\ Cuntz,\break S.\  
Echterhoff
eds.), 71--91, Springer-Verlag, New York (2000).

 

\bibitem[JP]{JP} \name{M.~Junge} and \name{G.~Pisier}, Bilinear forms on exact  
operator
spaces and $B(H)\otimes B(H)$, {\it Geom.\ Funct.\ Anal\/}.\ {\bf 5}  
(1995),
329--363.

\bibitem[Ki1]{Ki1} \name{E.~Kirchberg}, The Fubini theorem for Exact  
$C^*$-algebras,
{\it J.\ Operator Theory} {\bf 10} (1983), 3--8.

\bibitem[Ki2]{Ki2} \bibline, On nonsemisplit extensions, tensor  
products
and exactness of group\break $C^*$-algebras, {\it Invent.\ Math}.\ {\bf 112}  
(1993),
449--489.

\bibitem[KR]{KR} \name{R.V.~Kadison} and \name{J.~Ringrose}, {\it Fundamentals of the
Theory of Operator Algebras\/}, Vol.~II, {\it Academic Press}, New York  
(1986).

 \bibitem[La]{La} \name{F.~Larsen}, Powers of $R$-diagonal elements, {\it
     J.~Operator Theory} {\bf 47} (2002), 197--212.

\bibitem[Le]{Le} \name{F.~Lehner}, Computing norms of free operators with  
matrix
coefficients, {\it Amer.\ J.\ Math}.\ {\bf 121}  (1999), 453--486.

\bibitem[Me]{Me} \name{M.\ L.~Mehta}, {\it Random Matrices\/}, second edition,  
  Academic
Press, New York (1991).
 
\bibitem[Pas]{pas} \name{L.~Pastur}, A simple approach to global regime of  
random
   matriix theory, in {\it Mathematical Results in Statistical  
Mechanics} (Marseilles, 1998)
   (S.~Miracle-Sole, J.~Ruiz and V.~Zagrebnov, eds.), 429--454, World  
Sci.\ Publishing, River Edge, NJ (1999).

\bibitem[Pau]{pa} V\name{.~Paulsen}, {\it Completely Bounded Maps and
     Dilations\/}, {\it Pitman Research Notes in Mathematic} {\bf 146},   
Longman
     Scientific \& Technical, New York (1986).

\bibitem[Pe]{pe} \name{G.\ K.~Pedersen}, {\it Analysis Now\/},  {\it Grad.\  
Texts in Math\/}.\ {\bf 118}, Springer-Verlag, New York
   (1989).

\bibitem[P1]{P1} \name{G.~Pisier}, {\it The Volume of Convex Bodies and Banach  
Space
Geometry\/},  Cambridge Univ.\  Press, Cambridge  (1989).

\bibitem[P2]{P2} \bibline, A simple proof of a Theorem of Kirchberg and
related results on  $C^*$-norms, {\it J.\ Operator Theory} {\bf 35}
(1996), 317--335.

\bibitem[P3]{P3} \bibline, Quadratic forms in unitary operators, {\it  
Linear
Algebra and its Appl}.\ {\bf 267} (1997), 125--137.
 
\bibitem[RLL]{RLL} M.~R{\o}rdam, F.~Larsen, and N.J.~Laustsen, {\it An
Introduction to $K$-theory for\break $C^*$-algebras\/}, {Cambridge  
Univ.~Press}, Cambridge  (2000).

\bibitem[Ro]{Ro} \name{J.~Rosenberg}, Quasidiagonality and inuclearity   
(appendix to
strongly quasidiagonal operators by D.~Hadwin),  {\it J.\ Operator
Theory} {\bf 18}  (1987), 15--18.
 
 

\bibitem[T]{T} \name{S.~Thorbj{\o}rnsen}, Mixed moments of Voiculescu's  
Gaussian random
matrices, {\it J.\ Funct.\ Anal}.\ {\bf 176} (2000), 213--246.
  
\bibitem[V]{V} \name{A.~Valette}, An application of Ramanujan graphs to
$C^*$-algebra tensor products, {\it Discrete Math}.\ {\bf 167} (1997),  
597--603.

\bibitem[V1]{V0} \name{D.~Voiculescu}, A non commutative Weyl-von Neumann  
Theorem,
{\it Rev.\ Roum.\ Pures et Appl}.\ {\bf 21} (1976), 97--113.

\bibitem[V2]{V1} \bibline, Symmetries of some reduced free group
$C^*$-algebras, in {\it Operator Algebras and Their Connections with  
Topology
and Ergodic Theory}, {\it Lecture Notes in Math}.\ {\bf 1132},
  Springer-Verlag, New York  (1985), 556--588.

\bibitem[V3]{V2} \bibline, Circular and semicircular systems and
free product factors, in {\it Operator Algebras},  {\it Unitary  
Representations\/},
{\it Algebras, and Invariant Theory} (Paris, 1989), {\it Progress in  
Math}.\ {\bf 92},
  Birkh\"auser Boston, Boston, MA (1990), 45--60.

\bibitem[V4]{V3} \bibline, Limit laws for random matrices and free
products, {\it Invent.\   Math}.\ {\bf 104}  (1991), 202--220.

\bibitem[V5]{V4} \bibline, A note on quasi-diagonal $C^*$-algebras  
and
homotopy, {\it Duke Math.\ J}.\ {\bf 62} (1991), 267--271.

\bibitem[V6]{V5} \bibline, Around quasidiagonal operators, in {\it
Integral Equations and Operator Theory} {\bf 17}  (1993), 137--149.

\bibitem[V7]{V6} \bibline, Operations on certain noncommutative  
random
variables, in {\it Recent\break Advances in Operator Algebras\/} (Or\'leans  
1992),
{\it Ast\'erisque} {\bf 232} (1995), 243--275.

\bibitem[V8]{V8} \bibline, Free probability theory: Random matrices
and von Neumann algebras, {\it Proc.\  of the International Congress
of Mathematicians} (Z\"urich 1994), Vol.~1,  Birkh\"auser, Basel (1995),
227--241.

\bibitem[V9]{V9} \bibline, Free entropy, {\it Bull.~London  
Math.~Soc.}
   {\bf 34} (2002), 257--278.

\bibitem[VDN]{vdn} \name{D.~Voiculescu, K.~Dykema}, and \name{A.~Nica}, {\it Free  
Random
Variables}, {\it CMR Monograph Series} {\bf 1},  A.\ M.\ S.,  
Providence, RI (1992).
 

\bibitem[We]{We} \name{R.~Wegmann}, The asymptotic eigenvalue-distribution for  
a
certain class of random matrices, {\it J.~Math.~Anal.\ Appl}.\ {\bf 56}  
(1976), 113--132.

\bibitem[Wi]{Wi} \name{E.~Wigner}, Characterictic vectors of boardered matrices
with infinite dimensions, {\it Ann.\ of Math\/}.\ {\bf 62}  (1955),   
548--564.
 
\Endrefs
\end{document}